\documentclass[reqno,a4paper]{amsart}

\usepackage{ifpdf}
\usepackage{graphicx}
\usepackage{epsfig}
\usepackage{color}

\ifpdf
\usepackage[pdftex=true,a4paper]{hyperref}
\hypersetup{
bookmarks=true,
colorlinks=flase,
breaklinks=true,
pdfstartview={FitH},
unicode=true
}
\else
\usepackage[a4paper,hypertex]{hyperref}
\fi

\usepackage{a4wide}
\usepackage[utf8]{inputenc}
\usepackage[english,ngerman]{babel}
\usepackage{amssymb}
\usepackage{amsmath}
\usepackage{amsthm}
\usepackage{mathptmx} 
\usepackage[scaled=.90]{helvet}
\usepackage{mathrsfs}
\usepackage{bbold}
\usepackage{natbib} 
\usepackage{url} 
\usepackage{fancyhdr}

\numberwithin{equation}{section}

\newcommand{\titlestring}{The Aizenman-Sims-Starr and Guerra's schemes for
the SK model with multidimensional spins} 
\newcommand{\authorstring}{Anton Bovier and Anton Klimovsky}
\newcommand{\subjectstring}{Primary 60K35, 82B44, 60F10}
\newcommand{\keywordsstring}{Sherrington-Kirkpatrick model, multidimensional
spins, quenched large deviations, concentration of measure, Gaussian spins,
convexity, Parisi functional, Parisi formula}

\ifpdf
\hypersetup{
pdftitle = {\titlestring},
pdfauthor = {\authorstring},
pdfsubject = {\subjectstring},
pdfkeywords = {\keywordsstring},
pdfcreator = {},
pdfproducer = {}
}
\fi

\renewcommand{\P}{\mathbb {P}}
\newcommand{\E}{\mathbb {E}}

\newcommand{\R}{\mathbb {R}}
\newcommand{\Z}{\mathbb {Z}}
\newcommand{\N}{\mathbb {N}}

\newcommand{\I}{\mathbb{1}} 
\newcommand{\dd}{{\rm d}}
\newcommand{\ee}{{\rm e}}

\DeclareMathOperator{\card}{card}
\DeclareMathOperator{\supp}{supp}
\DeclareMathOperator{\var}{Var}
\DeclareMathOperator{\cov}{Cov}

\DeclareMathOperator{\interior}{int}
\DeclareMathOperator{\conv}{conv}

\DeclareMathOperator{\diam}{diam}
\DeclareMathOperator{\diag}{diag}

\DeclareMathOperator{\symmetric}{Sym}

\renewcommand{\Im}{\operatorname{Im}}

\newcommand{\eps}{\varepsilon}

\newtheorem{theorem}{Theorem}[section]

\newtheorem{corollary}{Corollary}[section]
\newtheorem{definition}{Definition}[section]

\newtheorem{lemma}{Lemma}[section]

\newtheorem{proposition}{Proposition}[section]

\newtheorem{remark}{Remark}[section] 
\newtheorem{assumption}{Assumption}[section]

\ifx\hypersetup\undefined\else
\fi

\begin{document}

\selectlanguage{english}

\begin{center}
\LARGE{\titlestring}
\end{center}
\vskip1cm
\begin{center}
\large{
Anton Bovier
\\
Institut für Angewandte Mathematik
\\
Rheinische Friedrich-Wilhelms-Universität Bonn
\\
Wegelerstrasse 6
\\
53115 Bonn
\\
Germany
\\
e-mail: bovier@uni-bonn.de
}
\end{center}
\vskip1cm
\begin{center}
\large{
Anton~Klimovsky\footnote{Supported in part by the DFG
Research Training Group ``Stochastic Processes and Probabilistic
Analysis'' and by the Helmholz-Gemeinschaft.}
\\
Department Mathematik
\\
Friedrich-Alexander-Universität Erlangen-Nürnberg
\\
Bismarckstrasse $1\frac{1}{2}$
\\
91054 Erlangen
\\
Germany
\\
e-mail: klimovsk@math.tu-berlin.de
}
\end{center}
\vskip1cm
\begin{center}
Abstract
\end{center}
We prove upper and lower bounds on the free energy of the
Sherrington-Kirkpatrick model with multidimensional spins in terms of
variational inequalities. The bounds are based on a multidimensional
extension of the Parisi functional. We generalise and unify the comparison scheme of
Aizenman, Sims and Starr and the one of Guerra involving the GREM-inspired
processes and Ruelle's probability cascades. For this purpose, an abstract
quenched large deviations principle of the Gärtner-Ellis type is obtained. We
derive Talagrand's representation of Guerra's remainder term for the
Sherrington-Kirkpatrick model with multidimensional spins. The derivation is 
based on well-known properties of Ruelle's probability cascades and the
Bolthausen-Sznitman coalescent. We study the properties of the
multidimensional Parisi functional by establishing a link with a certain class
of semi-linear partial differential equations. We embed the problem of strict
convexity of the Parisi functional in a more general setting and prove the
convexity in some particular cases which shed some light on the
original convexity problem of Talagrand. Finally, we prove the Parisi formula
for the local free energy in the case of multidimensional Gaussian a priori
distribution of spins using Talagrand's methodology of a priori estimates.
\vskip0.5cm
\noindent
\textbf{Key words:} \keywordsstring.
\vskip0.5cm
\noindent
\textbf{AMS 2000 Subject Classification:} \subjectstring.
\vskip0.5cm
\noindent

\newpage
\tableofcontents

\section{Introduction}
The Sherrington-Kirkpatrick (SK) model of a mean-field spin-glass has long
been one of the most enigmatic models of statistical mechanics. The
recent rigorous proof of the celebrated \emph{Parisi formula} for its
free energy, due to Talagrand \cite{TalagrandParisiFormula2006}, based on the
ingenious interpolation schemes of Guerra \cite{Guerra2003a} and Aizenman, Sims,
and Starr \cite{AizenmanSimsStarr2003} constitutes one of the major recent
achievements of probability theory. Recently, these results have been generalised to spherical
SK-models \cite{TalagrandSphericalSK} and to models with spins taking values in
a bounded subset of $\R$ \cite{panchenko-free-energy-generalized-sk-2005}.

In this paper, we are mainly concerned with the question of the validity of the
Parisi formula in the case where spins take values in a $d$-dimensional Riemannian manifold. We address the
issue of extending the approach of Aizenman, Sims and Starr, and the one of
Guerra to the multidimensional spins.
We study the properties of the multidimensional Parisi functional.
Motivated by a
problem posed by \cite{TalagrandParisiMeasures2006a}, we show the strict convexity of the local
Parisi functional in some cases.

We partially extend Talagrand's methodology of estimating the
remainder term to the multidimensional setting. In the case of
the multidimensional Gaussian a priori distribution of spins we prove the
validity of the Parisi formula in the low temperature regime.

\subsection*{Definition of the model}

Let $\Sigma \subset \R^d$ and denote $\Sigma_N \equiv \Sigma^N$. We
define  a family of  Gaussian processes
$
X
\equiv
\{
X(\sigma)
\}_{
\sigma \in \Sigma_N
}
$
as follows 
\begin{align}
\label{eq:as2:sk-with-multidimensional-spins-process}
X(\sigma) 
= 
X_N(\sigma)
\equiv
\frac{1}{N}
\sum_{
i,j=1
}^{N}
g_{i,j}
\langle
\sigma_i
,
\sigma_j
\rangle
,
\end{align}
where the \emph{interaction matrix}
$
G
\equiv
\{ 
g_{i,j} 
\}_{i,j=1}^N
$
consists of i.i.d. standard normal random variables and, for 
$
x, y \in \R^d
$,
$
\langle
x
,
y 
\rangle
\equiv
\sum_{u=1}^d
x_u
y_u
$
is the standard Euclidean scalar product. 
In what follows all random variables and processes are assumed to be centred.
We shall call $H_N(\sigma)\equiv -\sqrt{N} X_N(\sigma)$ a \emph{random
Hamiltonian} of our model.

Throughout the paper, we assume that we are given a large enough probability
space $(\Omega, \mathcal{F}, \P)$ such that all random variables under
consideration are defined on it. Without further notice, we
shall assume that all Gaussian random variables (vectors and processes) are
centred.

We shall be interested mainly in 
the \emph{free energy}
\begin{align}
\label{eq:as2:finite-volume-pressure}
p_N(\beta)
\equiv
\frac{1}{N}
\log
\int_{
\Sigma_N
}
\exp
\left(
\beta\sqrt{N}X(\sigma)
\right)
\dd 
\mu^{
\otimes N
}
(\sigma)
,
\end{align}
where $\beta \geq 0$ is the \emph{inverse temperature} and 
$
\mu \in \mathcal{M}_\text{f}(\Sigma)
$
is some arbitrary (not necessarily uniform or discrete) finite
\emph{a priori measure}. We assume that the a priori measure $\mu$ is such
that \eqref{eq:as2:finite-volume-pressure} is finite.
We shall be interested in proving bounds on the thermodynamic limits
of these quantities, e.g., on 

\begin{align}
\label{eq:as2:thermodynamic-pressure}
p(\beta)
\equiv
\lim_{
N \uparrow +\infty
}
p_N(\beta)
.
\end{align}
\begin{remark}
Note that there is no need to include the additional external field terms into 
the Hamiltonian \eqref{eq:as2:sk-with-multidimensional-spins-process},
since they could be absorbed into the a priori measure $\mu$.
\end{remark}

Mean-field spin-glass models (see, e.g., \cite{BovierBook2006}) with
multidimensional \emph{(Heisenberg)} spins were considered in
the theoretical physics 
literature, see, e.g., \cite{Sherrington2007} and references therein. Rigorous
results are, however, rather scarce. An early example is
\citep{Froehlich-Zegarlinski-Generalized-SK-1987}, where the authors
get 
bounds on the free energy in the high temperature regime. Methods of stochastic
analysis and large deviations are used in
\cite{Toubol1998} to identify the limiting distribution
of the partition function and also to obtain some information about
the geometry of the Gibbs measure for 
small  $\beta$. More recent treatments of the high temperature regime
using the very different methods are due to Talagrand
\cite{Talagrand2000}, see also
\cite[Subsection~2.13]{TalagrandSpinGlassesBook2003}. The importance of the SK model with multidimensional spins for understanding the ultrametricity of the
original SK model \cite{SherringtonKirkpatrick1975} 
(which corresponds to $d=1$ and $\mu$ being the Rademacher measure in the above
notations) was
emphasised in \cite{Talagrand2007}. 

For the SK model, Guerra's scheme gave historically the
first way to obtain the variational upper bound on the free energy in
terms of the Parisi functional. The scheme is based on the
comparison between two Gaussian processes: the first one being the original SK
Hamiltonian \eqref{eq:as2:sk-with-multidimensional-spins-process} and the other
one being a carefully chosen GREM
 inspired process indexed by $\sigma \in \Sigma_N$. The second
important ingredient is a recursively defined non-linear comparison functional
acting on the Gaussian processes indexed by $\sigma \in \Sigma_N$.

The Aizenman-Sims-Starr ($\text{AS}^2$) scheme
\cite{AizenmanSimsStarr2003,AizenmanSimsStarr2006} gives an intrinsic way to
obtain variational upper bounds on the free energy in the SK model. The
scheme is also based on a comparison between two Gaussian processes. The
first process is the sum of the original SK Hamiltonian $X$ and a GREM-inspired
process indexed by additional index space $\mathcal{A} \equiv \N^n$. The second 
one is another GREM-inspired process indexed by the extended configuration
space $\Sigma_N \times \mathcal{A}$. The scheme uses a comparison functional
defined on Gaussian processes indexed by the extended configuration space
equipped with the product measure between the original a priori measure and 
Ruelle's probability cascade (RPC) \cite{Ruelle1987}. The role of the comparison functional in the
$\text{AS}^2$ scheme is played by a free energy functional acting on the
Gaussian processes indexed by the extended configuration space. In \cite{Panchenko2007} Panchenko and Talagrand have reexpressed
Guerra's scheme for the SK model using the RPC.

Talagrand \cite{TalagrandParisiFormula2006} using Guerra's
scheme and the wealth of other ingenious analytical insights showed that the
variational upper bound is also the lower bound for the free energy in the SK
model. This established, hence, the remaining half of the Parisi formula.

A particular case ($d=1$, $\mu$ with bounded support) of the model we are
considering here was treated by Panchenko in
\cite{panchenko-free-energy-generalized-sk-2005}. He used the techniques of
\cite{TalagrandParisiFormula2006} to prove that in the case $d=1$ upper and lower bounds on the free energy coincide (cf.
\eqref{eq:as2:pressure-upper-bound} and
\eqref{vektor-sk-lower-bound-plus-remainder} in this chapter). However, the
results of \citep[Section 5 and the proofs of Theorems 2, 5 and
9]{panchenko-free-energy-generalized-sk-2005} are based on relatively detailed
differential properties of the optimal Lagrange multipliers in the saddle
point optimisation problem of interest. These properties are harder to obtain
in multidimensional situations such as that we are dealing with here. In fact,
as we show in Theorems~\ref{thm:as2:pressure-upper-bound} and
\ref{thm:as2:pressure-lower-bound}, one can obtain the same saddle point
variational principles without invoking the detailed properties of the
optimal Lagrange multipliers. This is achieved using a quenched large
deviations principle (LDP) of the Gärtner-Ellis type.

The most advanced recent study of spin-glass models with multidimensional
spins was attempted by Panchenko and Talagrand in \cite{PanchenkoTalagrandMultipleSKModels2006}, where the
multidimensional spherical spin-glass model was considered. The authors 
combined the techniques of
\cite{TalagrandParisiFormula2006,panchenko-free-energy-generalized-sk-2005} to
obtain partial results on the ultrametricity and also get some information on
the local free energy for their model.

\subsection*{Main results}
\label{sec:as2:main-results}
In this paper, we prove
upper and lower bounds on the free energy in the SK model 
with multidimensional spins  in terms
of variational inequalities involving the corresponding multidimensional
generalisation of the Parisi functional
(Theorems~\ref{thm:as2:pressure-upper-bound}, \ref{thm:as2:pressure-lower-bound},
\ref{thm:multidimensional-sk:free-energy-upper-bound},
\ref{eq:multidimensional-sk:vektor-sk-lower-bound-plus-remainder}). For this
purpose, we generalise and unify the $\text{AS}^2$ and
Guerra's schemes for the case of multidimensional spins, and employ a quenched
LDP which may be of independent interest
(Theorems~\ref{thm:as2:generic-quenched-ldp-upper-bound} and
\ref{vektor-sk-generic-quenched-ldp-lower-bound}). Both schemes are formulated
in a unifying framework based on the same comparison functional. The
functional acts on Gaussian processes indexed by an extended configuration space
as in the original $\text{AS}^2$ scheme. As a by-product, we provide also a short
derivation  of the remainder term in multidimensional Guerra's scheme
(Theorem~\ref{thm:lectire-03-an-analytic-projection-the-v-term-without-rpc})
using well-known properties of the RPC and the
Bolthausen-Sznitman coalescent. This gives a clear meaning to the remainder in
terms of averages with respect to a measure changed disorder. The change of
measure is induced by a reweighting of the RPC using the exponentials of the
GREM-inspired process\footnote{In $d=1$ the latter fact was also known to the author of \cite{Arguin2006}, private communication.}. See \cite{Panchenko2007} for
another approach in the case of the SK model ($d=1$).

We study the properties of the multidimensional Parisi functional by
establishing a link between the functional and
a certain class of non-linear partial differential equations (PDEs), see
Propositions~\ref{proposition:vektor-sk-parisis-pde-descrete-case},
\ref{eq:remainder:parisi-functional-lipschitzianity} and
Theorem~\ref{thm:remainder:piece-wise-viscosity-solutions-to-parisis-pde}. We
extend the Parisi functional to a continuous functional on a compact space
(Theorems~\ref{thm:remainder:extension-by-continuity-of-the-parisi-functional},
\ref{thm:remainder:piece-wise-viscosity-solutions-to-parisis-pde}). We show
that the class of PDEs corresponds to the Hamilton-Jacobi-Bellman (HJB)
equations induced by a linear problem of diffusion control
(Proposition~\ref{prop:remainder:hjb-for-value-function}). Motivated by a
problem posed by \cite{TalagrandParisiMeasures2006a}, we show the strict convexity of the local
Parisi functional in some cases
(Theorem~\ref{thm:remainder:strict-convexity-parisi-functional-1d}).

We partially extend Talagrand's methodology of estimating the
remainder term to the multidimensional setting 
(Theorem~\ref{thm:lectire-03-an-analytic-projection-the-v-term-without-rpc},
Proposition~\ref{eq:remainder:phi-2-upper-bound},
Theorem~\ref{thm:remainder:the-conditional-parisi-formula}).
In the case of multidimensional Gaussian a priori distribution of spins we 
prove the validity of the Parisi formula 
(Theorem~\ref{thm:gaussian:the-local-low-temperature-parisi-formula}).

We partially extend Talagrand's methodology of estimating the
remainder term to the multidimensional setting 
(Theorem~\ref{thm:lectire-03-an-analytic-projection-the-v-term-without-rpc},
Proposition~\ref{eq:remainder:phi-2-upper-bound},
Theorem~\ref{thm:remainder:the-conditional-parisi-formula}). Though the main
technical problem of the methodology in the general multidimensional setting
remains
(Remark~\ref{rem:remainder:the-main-thechnical-problem-of-talagrands-approach}).
In the case of the multidimensional Gaussian a priori distribution of spins we 
prove the validity of the Parisi formula 
(Theorem~\ref{thm:gaussian:the-local-low-temperature-parisi-formula}).

Below we introduce the notations, assumptions and formulate our main
results. The other results (mentioned above) are formulated and proved in the
subsequent sections.

\begin{assumption}
Suppose that the configuration space $\Sigma$ is bounded and such that
$0 \in \interior \conv \Sigma$, where $\conv \Sigma$ denotes the
convex hull of $\Sigma$.
\end{assumption}
The examples listed below verify this assumption: 
\begin{enumerate}
\item
Multicomponent Ising spins. $\Sigma = \{-1;1\}^d$ -- the discrete hypercube.
\item
Heisenberg spins.
$
\Sigma 
= 
\left\{
\sigma \in \R^d 
: 
\Vert \sigma \Vert_2 
= 
1 
\right\}
$
-- the unit Euclidean sphere.
\item
$
\Sigma 
= 
\left\{
\sigma \in \R^d 
: 
\Vert \sigma \Vert_2 
\leq 
1 
\right\}
$
-- the unit Euclidean ball. 
\end{enumerate}
\begin{remark}
The boundedness assumption can be relaxed and replaced by
concentration properties of the a priori measure. In Section~\ref{sec:gaussian-spins:gaussian-spins} we will
exemplify this in the case of a Gaussian a priori distribution. In general a
subgaussian distribution will suffice.
\end{remark}
Consider the space of all \emph{symmetric matrices}
$
\symmetric(d)
\equiv
\left\{
\Lambda
\in
\R^{d\times d}
\mid
\Lambda = \Lambda^{*}
\right\}
$. 
Denote 
\begin{align*}
\symmetric^+(d)
\equiv
\left\{
\Lambda
\in
\symmetric(d)
\mid
\Lambda
\succeq
0
\right\}
,
\end{align*}
where the notation $\Lambda \succeq 0$ means that the matrix $\Lambda$ is
non-negative definite.
We equip the space 
$
\symmetric(d)
$
with the Frobenius (Hilbert-Schmidt) norm
\begin{align*}
\Vert
M
\Vert_\text{F}^2
\equiv
\sum_{u,v=1}^d
M_{u,v}^2
,
\quad
M \in \symmetric(d)
.
\end{align*}
We shall also denote the corresponding (tracial) scalar product by $\langle
\cdot, \cdot \rangle$. 
For 
$
r 
>
\max
\{
\Vert \sigma \Vert_2^2
:
\sigma
\in
\Sigma
\}
$, 
define 
\begin{align*}
\mathcal{U} 
\equiv
\left\{
U
\in
\symmetric(d)
\mid
U \succeq 0
,
\Vert
U
\Vert_2
\leq
r
\right\}
.
\end{align*}
We will call the set $\mathcal{U}$ the \emph{space of the
admissible self-overlaps}. In analogy to the usual overlap in the standard SK
model, we define, for two configurations,
$
\sigma^{(i)}
= 
(\sigma^{(i)}_1, \sigma^{(i)}_2,\dots,\sigma^{(i)}_N)
\in
\Sigma_N
$, $i=1,2$, the (mutual) \emph{overlap matrix} 
$
R_N(\sigma^{(1)},\sigma^{(2)})
\in
\R^{
d \times d
}
$
whose entries are given by 
\begin{align}
\label{eq:as2:overlap-matrix-entries}
R_N(\sigma^{(1)},\sigma^{(2)})_{u,v}
\equiv 
\frac{1}{N}\sum_{i=1}^N
\sigma^{(1)}_{i,u}\sigma^{(2)}_{i,v},
\quad
u,v
\in
[1;d] \cap \N
.
\end{align}
Fix an \emph{overlap matrix} 
$
U \in \mathcal{U}
$.
Given 
a subset
$
\mathcal{V}
\subset
\mathcal{U}
$,
define the set of the \emph{local configurations},
\begin{align*}
\Sigma_N(\mathcal{V})
\equiv
\left\{
\sigma \in \Sigma_N
\mid
R_N(\sigma, \sigma)
\in 
\mathcal{V}
\right\}
.
\end{align*}
Next, define the \emph{local free energy}
\begin{align}
\label{eq:as2:local-free-energy}
p_N(\mathcal{V})
\equiv
\frac{1}{N}
\log
\int
_{\Sigma_N(\mathcal{V})}
\ee^
{
\beta\sqrt{N}X(\sigma)
}
\dd 
\mu^{
\otimes N
}
(\sigma)
.
\end{align}
We also define
\begin{align}
\label{eq:as2:local-limiting-free-energy}
p(\mathcal{V})
\equiv
p(\beta,\mathcal{V})
\equiv
\lim_{
N \uparrow +\infty
}
p_N(\mathcal{V})
,
\end{align}
where the existence of the limit follows from a result of Guerra and
Toninelli \cite[Theorem~1]{Guerra-Toninelli-Generalized-SK-2003}. 
Consider a sequence of matrices 
$
\mathcal{Q}
\equiv
\{
Q^{(k)}
\in
\symmetric(d)
\}_{
k=0
}^{n+1}
$
such that
\begin{align}
\label{eq:multidimesnional-sk:monotonicity-of-q-matrix-sequence}
0 
\equiv
Q^{(0)}
\prec
Q^{(1)}
\prec
\ldots
\prec
Q^{(n+1)}
\equiv
U
,
\end{align}
where the ordering is understood in the sense of the corresponding quadratic
forms. Consider in addition  a partition of the unit interval
$
x
\equiv
\{
x_k
\}_{k=0}^{n+1}
$,
i.e.,
\begin{align}
\label{eq:partition-of-the-unit-interval}
0 
\equiv
x_0
<
x_1
<
\ldots
<
x_{n+1}
\equiv
1
.
\end{align}
Let 
$
\{z^{(k)}\}_{k=0}^n
$ 
be a sequence of independent Gaussian $d$-dimensional vectors with
\begin{align*}
\cov
\left[
z^{(k)}
\right]
=
Q^{(k+1)}
-
Q^{(k)}
.
\end{align*}
Given 
$
\Lambda \in \symmetric(d)
$, 
define
\begin{align}
\label{eq:as2:x-end-condition}
X_{n+1}(x, \mathcal{Q}, U, \Lambda)
\equiv
\log
\int
_\Sigma
\exp
\left(
\sqrt{2}
\beta
\left\langle
\vphantom{\sum}
\smash{
\sum_{k=0}^{n}
}
z_k
,
\sigma
\right\rangle
+
\langle
\Lambda
\sigma
,
\sigma
\rangle
\right)
\dd \mu(\sigma)
.
\end{align}
Define, for $k \in \{n,\dots,0\}$, by a descending recursion,
\begin{align}
\label{eq:vektor-sk-x-k-recursive-definition}
X_k(x, \mathcal{Q}, U, \Lambda)
\equiv
\frac{1}{x_k}
\log
\E_{z^{(k)}}
\left[
\exp
\left(
x_k X_{k+1}(x, \mathcal{Q}, U, \Lambda)
\right)
\right]
\end{align}
with
\begin{align}
\label{eq:multidimensional-sk:x-0}
X_0(x, \mathcal{Q}, U, \Lambda)
\equiv
\E_{z^{(0)}}
\left[
X_{1}(x, \mathcal{Q}, U, \Lambda)
\right]
,
\end{align}
where 
$
\E_{
z^{(k)}
}
\left[
\cdot
\right]
$
denotes the expectation with respect to the $\sigma$-algebra generated by the
random vector $z^{(k)}$.
\begin{remark}
Section~\ref{sec:remainder:the-filtered-d-dimensional-grem} contains the more
general framework of dealing with the recursive quantities
\eqref{eq:multidimensional-sk:x-0} which in particular brings to light the
links with certain non-linear parabolic PDEs. For these PDEs the recursion
\eqref{asm:multidimensional-sk:hadamard-squares} is closely related to an
iterative application of the well-known Hopf-Cole transformation, see, e.g.,
\cite{Evans1998}.
\end{remark}
Define the \emph{local Parisi functional} as
\begin{align}
\label{eq:as2:local-parisi-functional}
f(x, \mathcal{Q}, U, \Lambda)
\equiv
-
\langle
\Lambda
,
U
\rangle
-
\frac{\beta^2}{2}
\sum_{k=1}^{n}
x_k
\left(
\Vert
Q^{(k+1)}
\Vert_\text{F}^2
-
\Vert
Q^{(k)}
\Vert_\text{F}^2
\right)
+
X_0(x, \mathcal{Q}, U, \Lambda)
.
\end{align}
\begin{assumption}[Hadamard squares]
\label{asm:multidimensional-sk:hadamard-squares} 
We shall say that a
sequence, 
$
\{
Q^{(i)}
\}_{
i=1
}^n
$, 
of matrices satisfies Assumption 
\ref{asm:multidimensional-sk:hadamard-squares}, if 
\begin{align}
\label{eq:multidimensional-sk:hadamard-squares}
\left(
Q^{(1)}
\right)^{
\odot
2
}
\prec
\ldots
\prec
\left(
Q^{(n)}
\right)^{
\odot
2
}
\prec
\left(
Q^{(n+1)}
\right)^{
\odot
2
}
.
\end{align}
\end{assumption}
\begin{remark}
\label{rem:as2:hadamard-squares-assumption}
The above assumption on the matrix order parameters $\mathcal{Q}$ is necessary
only to employ the $\text{AS}^2$ scheme. In contrast, Guerra's scheme
(Theorems~\ref{thm:multidimensional-sk:free-energy-upper-bound} and
\ref{eq:multidimensional-sk:vektor-sk-lower-bound-plus-remainder}) does not require the above assumption.
\end{remark}
One may verify that  the matrices $q$ and $\rho$ in
\cite[Theorems~2.13.1 and 2.13.2]{TalagrandSpinGlassesBook2003}
correspond to the matrices $Q^{(1)}$ and $Q^{(2)}$ of this
paper ($n=1$). (See also
\eqref{eq:as2:high-temperature-mean-field-eqns} below.) Furthermore, a
straightforward application of the Cauchy-Schwarz inequality shows that the
matrices $q$ and $\rho$ actually satisfy 
Assumption~\ref{asm:multidimensional-sk:hadamard-squares}. We also note that
in the simultaneous diagonalisation scenario in which the matrices in
\eqref{eq:multidimesnional-sk:monotonicity-of-q-matrix-sequence} are
diagonalisable in the same orthogonal basis (see
Sections~\ref{sec:remainder:simultaneous-diagonalisation-scenario} and
\ref{sec:gaussian:simultaneous-diagonalisation-scenario}) this assumption is also satisfied.

The first main result of the present paper uses the $\text{AS}^2$
scheme to establish the upper bound on the limiting free energy $p(\beta)$ in
terms of the saddle point problem for the local 
Parisi functional \eqref{eq:as2:local-parisi-functional}. 
\begin{theorem}
\label{thm:as2:pressure-upper-bound}
For any closed set $\mathcal{V} \subset \symmetric(d)$, we have 
\begin{align}
\label{eq:as2:pressure-upper-bound}
p(\mathcal{V})
\leq
\sup_{
U
\in
\mathcal{V}
\cap
\mathcal{U}
}
\inf_{
(
x
,
\mathcal{Q}
,
\Lambda
)
}
f(x, \mathcal{Q}, \Lambda, U)
,
\end{align}
where the infimum runs over all $x$ satisfying
\eqref{eq:partition-of-the-unit-interval}, all $\mathcal{Q}$ satisfying both
\eqref{eq:multidimesnional-sk:monotonicity-of-q-matrix-sequence} and
Assumption~\ref{asm:multidimensional-sk:hadamard-squares}, and all $\Lambda \in
\symmetric(d)$.
\end{theorem}
We were not able to prove in general that the r.h.s. of
\eqref{eq:as2:pressure-upper-bound} gives also the lower bound to the
thermodynamic free energy. See, however,
Theorem~\ref{thm:gaussian:the-local-low-temperature-parisi-formula} for a positive example.

To formulate the lower bound on \eqref{eq:as2:thermodynamic-pressure} 
we need some additional definitions.

Let the \emph{comparison index space} be
$
\mathcal{A} \equiv \N^n
$.
Given 
$
\alpha^{(1)},\alpha^{(2)} \in \mathcal{A}
$, 
define 
\begin{align}
\label{eq:as2:matrix-lexicographic-overlap}
Q(\alpha^{(1)},\alpha^{(2)})
\equiv
Q^{(q_\text{L}(\alpha^{(1)},\alpha^{(2)}))}
,
\end{align}
where 
$
q_\text{L}(\alpha^{(1)},\alpha^{(2)})
$
is the (normalised) \emph{lexicographic overlap} defined as follows 
\begin{align}
\label{eq:introduction:lexicographic-overlap}
q_\text{L}(\alpha^{(1)},\alpha^{(2)})
\equiv
1
+
\begin{cases}
0
,
&
\alpha^{(1)}_1
\neq
\alpha^{(2)}_1
\\
\max
\left\{
k
\in
[1;n] \cap \N
:
[\alpha^{(1)}]_k
=
[\alpha^{(2)}]_k
\right\}
,
&
\text{otherwise.}
\end{cases}
\end{align}
Given a $d\times d$-matrix
$
M
$ and $p \in \R$,  we  denote by 
$
M^{\odot p}
$
the $d\times d$-matrix 
with entries
\begin{align*}
\left(
M^{\odot p}\right)_{u,v}
\equiv
\left(
M_{u,v}
\right)^p
.
\end{align*}
The matrix valued lexicographic overlap
\eqref{eq:as2:matrix-lexicographic-overlap} can be used to construct the multidimensional ($d\geq1$) versions of the
GREM (see, e.g., \cite{BovierKurkova2007}
and references therein for a review of the results on the one-dimensional
case of the model). Here we shall need the following two GREM-inspired
real-valued Gaussian processes:
$ 
A
\equiv
\{
A(\sigma,\alpha)
\}_{
\sigma
\in
\Sigma_N
,
\alpha
\in
\mathcal{A}
}
$ 
and
$
B
\equiv
\{
B(\alpha)
\}_{
\alpha
\in 
\mathcal{A}
}
$
with covariance structures
\begin{align*}
\E
\left[
A(\sigma^{(1)},\alpha^{(1)})
A(\sigma^{(2)},\alpha^{(2)})
\right]
=
2
\langle
R(\sigma^{(1)},\sigma^{(2)})
,
Q(\alpha^{(1)},\alpha^{(2)})
\rangle
,
\\
\E
\left[
B(\alpha^{(1)})
B(\alpha^{(2)})
\right]
=
\Vert
Q(\alpha^{(1)},\alpha^{(2)})
\Vert_\text{F}^2
.
\end{align*}
Note that the process $A$ can be represented in the following form:
\begin{align}
\label{eq:as2:a-process-definition}
A(\sigma,\alpha)
=
\left(
\frac{
2
}{
N
}
\right)^{
1/2
}
\sum_{i=1}^N
\langle
A_i(\alpha)
,
\sigma_i
\rangle
,
\end{align}
where 
$
\{
A_i
\equiv
\{
A_i(\alpha)
\}_{
\alpha
\in
\mathcal{A}
}
\}_{i=1}^N
$ are the i.i.d. (for different indices $i$) Gaussian $\R^d$-valued processes
with the following covariance structure: for  
$
i \in [1;N] \cap \N
$, 
for all 
$
\alpha^{(1)}, \alpha^{(2)} \in \mathcal{A}
$ 
and all
$
u,v \in [1;d] \cap \N
$
assume that the following holds
\begin{align*}
\E
\left[
A_i(\alpha^{(1)})_u
A_i(\alpha^{(2)})_v
\right]
=
Q(\alpha^{(1)},\alpha^{(2)})_{
u,v
}
.
\end{align*}
Given 
$
t \in [0;1]
$, 
we define the \emph{interpolating $\text{AS}^2$ Hamiltonian}
\begin{align}
\label{eq:as2:interpolating-as2-hamiltonian}
H_t(\sigma,\alpha)
\equiv
\sqrt{t}
\left(
X(\sigma)
+
B(\alpha)
\right)
+
\sqrt{1-t}
A(\sigma,\alpha)
.
\end{align}
Next, we define the random probability measure
$\pi_N\in\mathcal{M}_1(\Sigma_N\times\mathcal{A})$ through
\begin{align*}
\pi_N 
\equiv
\mu^{
\otimes N
}
\otimes
\xi
,
\end{align*}
where $\xi = \xi(x)$ is the RPC
\cite{Ruelle1987}. We denote by  
$\{
\xi(\alpha)
\}_{
\alpha
\in
\mathcal{A}
}
$
the enumeration of the atom locations of the RPC and consider
the enumeration as a random measure on $\mathcal{A}$ (independent of all other
random variables around). Define  the \emph{local 
$
\text{AS}^2
$ 
Gibbs measure} 
$
\mathcal{G}_N(t,x,\mathcal{Q},U,\mathcal{V})
$
by
\begin{align}\label{immer.1}
\mathcal{G}_N(t,x,\mathcal{Q},U,\mathcal{V})
\left[
f
\right]
\equiv
\frac
{1}
{Z_N(t,\mathcal{V})}
\int_{
\Sigma_N(\mathcal{V})
\times
\mathcal{A}
}
f(\sigma,\alpha)
\ee^{
\sqrt{N}\beta
H_t(\sigma,\alpha)
}
\dd
\pi_N
(\sigma,\alpha)
,
\end{align}
where 
$
f:
\Sigma_N\times\mathcal{A}
\to
\R
$
is an arbitrary  measurable function for which the right-hand side of
(\ref{immer.1}) is finite.
For $\mathcal{V} \subset \mathcal{U}$, define the \emph{$\text{AS}^2$
remainder term} as
\begin{align}
\label{ab.200}
&
\mathcal{R}_N(x, \mathcal{Q}, U, \mathcal{V})
\nonumber
\\
&\equiv
-
\frac{1}{2}
\int
_0^{1}
\E
\left[
\mathcal{G}_N(t,x,\mathcal{Q},U,\mathcal{V})
\otimes 
\mathcal{G}_N(t,x,\mathcal{Q},U,\mathcal{V})
\left[
\Vert
R_N(\sigma^{(1)},\sigma^{(2)})
-
Q(\alpha^{(1)},\alpha^{(2)})
\Vert
_\text{F}^2
\right]
\right]
\dd t
.
\end{align}
We define  also the \emph{limiting $\text{AS}^2$
remainder term} 
\begin{align}
\label{eq:as2:reaminder-term}
\mathcal{R}(x, \mathcal{Q}, U)
\equiv
\lim_{
\eps
\downarrow
+0
}
\lim_{
N
\uparrow
\infty
}
\mathcal{R}_N(x, \mathcal{Q}, B(U,\eps))
\leq
0
,
\end{align}
where 
$
B(U,\eps)
$
is the ball with centre $U$ and radius $\eps$. (The existence of the limiting
remainder term is proved in Theorem~\ref{thm:as2:pressure-lower-bound}.)

The second main result of this paper uses the $\text{AS}^2$ scheme to
establish a lower bound on
\eqref{eq:as2:thermodynamic-pressure} 
in terms of the same saddle point Parisi-type functional as in the upper bound
which includes, however, the non-positive remainder term
\eqref{eq:as2:reaminder-term}. 
In one-dimensional situations Talagrand \cite{TalagrandParisiFormula2006} and
Panchenko \cite{panchenko-free-energy-generalized-sk-2005}, respectively, have shown that
the corresponding error term vanishes on the optimiser of the Parisi functional. 
\begin{theorem}
\label{thm:as2:pressure-lower-bound}
For any open set $\mathcal{V} \subset \symmetric(d)$, we have
\begin{align}
\label{vektor-sk-lower-bound-plus-remainder}
p(
\mathcal{V}
)
\geq
\sup_{
U
\in
\mathcal{V}
\cap
\mathcal{U}
}
\inf_{
(
x
,
\mathcal{Q}
,
\Lambda
)
}
\left[
f(x, \mathcal{Q}, \Lambda, U)
+
\mathcal{R}(x, \mathcal{Q}, U)
\right]
,
\end{align}
where the infimum runs over all $x$ satisfying
\eqref{eq:partition-of-the-unit-interval}, 
all $\Lambda \in\symmetric(d)$, and all $\mathcal{Q}$ satisfying both
\eqref{eq:multidimesnional-sk:monotonicity-of-q-matrix-sequence} and
Assumption~\ref{asm:multidimensional-sk:hadamard-squares}.
\end{theorem}
\begin{remark}
The comparison scheme of Guerra \cite{Guerra2003a} (see also more recent accounts
\cite{Talagrand2007a}, \cite{GuerraReview2005} and
\cite{AizenmanSimsStarr2006}) is also applicable to our model and is covered by our quenched LDP approach, see Theorems~\ref{thm:multidimensional-sk:free-energy-upper-bound} and
\ref{eq:multidimensional-sk:vektor-sk-lower-bound-plus-remainder} for the
formal statements. Guerra's scheme seems to be more amenable
(compared to the Aizenman-Sims-Starr one) for Talagrand's remainder estimates \cite{TalagrandParisiFormula2006}, see Section
\ref{sec:remainder:talagrands-remainder-estimates}. The scheme is based on
the following interpolation
\begin{align}
\label{eq:as2:guerras-interpolation}
\widetilde{H}_t(\sigma,\alpha)
\equiv
\sqrt{t}
X(\sigma)
+
\sqrt{1-t}
A(\sigma,\alpha)
\end{align}
which induces the corresponding local 
Gibbs measure \eqref{immer.1} and remainder
term \eqref{ab.200} by substituting
\eqref{eq:as2:interpolating-as2-hamiltonian} with
\eqref{eq:as2:guerras-interpolation}. Guerra's scheme does not include the
process $B$ and, hence, does not require
Assumption~\ref{asm:multidimensional-sk:hadamard-squares}. Recovering the terms
corresponding to $\Phi_N(x,\mathcal{U})[B]$ (see, \eqref{eq:as2:phi-of-b}) in
the Parisi functional requires then a short additional calculation (Lemma~\ref{thm:as2:guerras-interpolation}).
\end{remark}
Note that the results of Talagrand \cite[Theorems~2.13.2 and
2.13.3]{TalagrandSpinGlassesBook2003} imply that at least in the high
temperature region (i.e., for small enough $\beta$) the Parisi formula for the
SK model with multidimensional spins is valid with $n=1$
\begin{align}
\label{eq:as2:high-temperature-parisi-formula}
p(\beta)
=
f(x, \mathcal{Q}^*, 0, U^*)
=
\sup_{
U
\in
\mathcal{U}
}
\inf_{
(
\mathcal{Q}
,
\Lambda
)
}
f(x, \mathcal{Q}, \Lambda, U)
,
\end{align}
where the matrices $Q^{*(2)} = U^*$ and $Q^{*(1)}$ solve the following system
of equations:
\begin{align}
\label{eq:as2:high-temperature-mean-field-eqns}
\begin{cases}
\partial_{
Q^{(2)}_{u,v}
}
f(x, \mathcal{Q}^*, 0, U^*)
=
0
,
&
u,v \in [1;d] \cap \N
,
\\
\partial_{
Q^{(1)}_{u,v}
}
f(x, \mathcal{Q}^*, 0, U^*)
=
0
,
&
u,v \in [1;d] \cap \N
.
\end{cases}
\end{align}
Note that the system \eqref{eq:as2:high-temperature-mean-field-eqns} coincides
with the mean-field equations obtained in \cite[see (2.469) and (2.470)]{TalagrandSpinGlassesBook2003}.

Let $\Sigma \equiv \R^d$ and fix some vector $h \in \R^d$. Let $\mu
\in \mathcal{M}_{\text{f}}(\Sigma)$ be the finite measure with the following
density (with respect to the Lebesgue measure $\lambda$ on $\Sigma$) 
\begin{align}
\label{gaussian-spins:a-priori-measure}
\frac{
\dd \mu
}
{
\dd \lambda
}
(\sigma)
=
\left(
\frac{
\det C
}{
\left(
2\pi
\right)^{d}
}
\right)^{1/2}
\exp
\left(
-
\frac{1}{2}
\langle
C
\sigma
,
\sigma
\rangle
+
\langle
h
,
\sigma
\rangle
\right)
,
\end{align}
where $C \in \symmetric^+(d)$. Note that, given $m \in \R^d$ and $C \in
\symmetric^+(d)$ such that $\det C \neq 0$, the density
\eqref{gaussian-spins:a-priori-measure} with
$h \equiv C m$ coincides (up to the constant factor
$\exp\left(-\frac{1}{2}\langle C m, m \rangle\right)$) with the Gaussian
density with covariance matrix $C^{-1}$ and mean $m$.
\begin{remark}
It turns out that only matrices $C$ with sufficiently large
eigenvalues will result in finite global free energy, cf. Lemma
\ref{lem:lecture-03:gaussian-spins:saddle-point-of-the-crisanti-sommers-functional}.
The local free energy is, in contrast, always finite, see
Lemma~\ref{lem:lecture-03:optimum-of-the-crisanti-sommers-functional} and
Theorem~\ref{thm:gaussian:the-local-low-temperature-parisi-formula}.
\end{remark}
Consider the function 
$
f: (0:+\infty)^2 \to \R
$ 
given by
\begin{align}
\label{eq:lecture-03:inf-spherical-cs-two-clause-formula}
f(c,u)
=
\begin{cases}
\beta^2 u^2
+
\log c u 
-
c u
+
1
,
&
u
\in
(0
;
\frac{
\sqrt{2}
}{
2\beta
}
]
,
\\
(
2\sqrt{2} \beta 
-
c
)
u
+
\log 
\frac{c}{\beta}
-
\frac{1}{2}
\left(
1+\log 2
\right)
,
&
u
\in 
(
\frac{
\sqrt{2}
}{
2\beta
}
;
+\infty
]
.
\end{cases}
\end{align}
The following result shows that, at least, in the highly symmetric situation
\eqref{gaussian-spins:a-priori-measure} with $h=0$ the multidimensional Parisi
formula indeed holds true (see
Lemma~\ref{lem:lecture-03:optimum-of-the-crisanti-sommers-functional} for an
explanation why the result is indeed a Parisi formula).
\begin{theorem}
\label{thm:gaussian:the-local-low-temperature-parisi-formula}
Let $\mu$ satisfy \eqref{gaussian-spins:a-priori-measure} with $h = 0$.
Assume that the matrices $U$ and $C$ are
simultaneously diagonalisable in the same basis. Denote by 
$
\{ c_{v} \in \R_+ \}_{v=1}^d
$ 
and 
$
\{ u_{v} \in \R_+ \}_{v=1}^d
$
the eigenvalues of the matrices $C$ and $U$, respectively.
Moreover, assume that
$
\min_v u_v > 0
$ 
and 
$
\min_v c_v > 0
$.

Then we have
\begin{align*}
\lim_{
\eps
\downarrow +0
}
\lim_{N \uparrow +\infty}
p_N(
\Sigma_N(B(U,\eps))
)
=
\sum_{v=1}^d
f(c_v,u_v)
.
\end{align*}

\end{theorem}
\begin{remark}
Close results have previously been obtained in the case of the spherical model
in \cite{PanchenkoTalagrandMultipleSKModels2006}, from where we borrow the
general methodology of the proof of the
Theorem~\ref{thm:gaussian:the-local-low-temperature-parisi-formula}. As noted
in \cite{PanchenkoTalagrandMultipleSKModels2006}, another more straightforward way to obtain the
Theorem~\ref{thm:gaussian:the-local-low-temperature-parisi-formula} is to
diagonalise the interaction matrix $G$ and use the properties of the
corresponding random matrix ensemble.
\end{remark}

\subsection*{Organisation of the paper}
The rest of the present paper is organised as follows. In
Section~\ref{sec:some-basic-results} we record some basic properties
of the covariance structure of the process $X$ and 
establish the relevant concentration of measure results.
The section contains also the tools allowing to compare and interpolate between
the free energy-like functionals of different
Gaussian processes. In Section~\ref{sec:sk-with-multidimensional-spins:gaertner-ellis} we 
derive a quenched LDP of the Gärtner-Ellis type under measure concentration
assumptions. Section~\ref{sec:as2:free-energy-upper-and-lower-bounds} contains the
derivation (based on the $\text{AS}^2$
scheme) of the upper and lower bounds on
the free energy of the SK model with multidimensional spins in terms of the
saddle point of the Parisi-like functional. In Section~\ref{sec:as2:guerras-scheme} we employ the ideas of Guerra's
comparison scheme in order to obtain the upper and lower bounds on the free
energy and
we also get a useful analytic representation of the remainder term. In
Section~\ref{sec:remainder:the-parisi-functional-in-terms-of-differential-equations} we study the properties of the multidimensional Parisi functional. 
Section~\ref{sec:remainder:talagrands-remainder-estimates} contains the
partial extension of Talagrand's remainder term estimates to the case of multidimensional spins. 
In Section~\ref{sec:gaussian-spins:gaussian-spins} a case of
Gaussian a priori distribution of spins is considered and the corresponding
local Parisi formula is proved. In the appendix we prove the almost
super-additivity of the local free energy, as an application of the Gaussian
comparison results of
Subsection~\ref{sec:gaussian-comparison-inequalities-for-free-energy-like-functionals}.

\section{Some preliminary results}
\label{sec:some-basic-results}

\subsection{Covariance structure}

Our definition of the overlap matrix in
\eqref{eq:as2:overlap-matrix-entries} 
is  motivated by the
fact that, as can be seen from a  straightforward computation
\begin{align}
\label{eq:multidimensional-sk:covariance}
\E 
\left[
X_N(\sigma^{(1)})X_N(\sigma^{(2)}) 
\right]
= 
\sum_{u,v=1}^d 
\left(
R_N(\sigma^{(1)},\sigma^{(2)})_{
u,v
}
\right)^2
=
\Vert 
R_N(\sigma^{(1)},\sigma^{(2)}) 
\Vert_{2}^2
,
\end{align}
that is, the
the covariance structure of the process $X_N(\sigma)$
is given by the square of the Frobenius (Hilbert-Schmidt) norm of the 
matrix $R_N(\sigma^{(1)},\sigma^{(2)})$.
The basic properties of the overlap matrix
are summarised in the following proposition.
\begin{proposition}
\label{thm:lecture-01-matrix-overlap-properties}
We have, for all $\sigma^{(1)},\sigma^{(2)},\sigma \in \Sigma_N$,
\begin{enumerate}
\item
\emph{Matrix representation.}
$
R_N(\sigma^{(1)},\sigma^{(2)}) 
= 
\frac{1}{N}
\left(
\sigma^{(1)}
\right)^*
\sigma^{(2)}
$.
\item
\emph{Symmetry \#1.}
$
R^{u,v}_N(\sigma^{(1)},\sigma^{(2)})
=
R^{v,u}_N(\sigma^{(2)},\sigma^{(1)})
$.
\item 
\emph{Symmetry \#2.}
$
R^{u,v}_N(\sigma,\sigma)
=
R^{v,u}_N(\sigma,\sigma)
$.
\item
\emph{Non-negative definiteness \#1.}
$
R_N(\sigma,\sigma) \succeq 0
$.
\item
\emph{Non-negative definiteness \#2.}
\begin{align*}
\begin{bmatrix}
R_N(\sigma^{(1)},\sigma^{(1)})
&
R_N(\sigma^{(1)},\sigma^{(2)})
\\
R_N(\sigma^{(1)},\sigma^{(2)})^*
&
R_N(\sigma^{(2)},\sigma^{(2)})
\end{bmatrix}
\succeq
0
.
\end{align*}
\item
Suppose 
$
U \equiv R_N(\sigma^{(1)},\sigma^{(1)}) 
=
R_N(\sigma^{(2)},\sigma^{(2)})
$, 
then
\begin{align*}
\Vert
R(\sigma^{(1)},\sigma^{(2)})
\Vert_\text{F}^2
\leq
\Vert
U
\Vert_\text{F}^2
.
\end{align*}
\end{enumerate}
\end{proposition}
\begin{proof}
The proof is straightforward.
\end{proof}
\subsection{Concentration of measure}
The following concentration of measure result for the free energy is standard.
\begin{proposition}
\label{vektor-sk-concentration-of-measure-for-free-energies}
Let $(\Sigma,\mathfrak{S})$ be a Polish space. Suppose $\mu$ is a random
finite measure on $\Sigma$. Suppose, moreover, that $X(\sigma)$, $\sigma \in
\Sigma$ is the family of Gaussian random variables independent of $\mu$ which
possesses a bounded covariance, i.e.,
\begin{align}
\label{vektor-sk-covariance-lemma-assumption}
\text{
there exists $K>0$ such that
}
\sup_{
\sigma^{(1)},\sigma^{(2)}
\in
\Sigma
}
\vert
\cov 
(X(\sigma^{(1)}),X(\sigma^{(2)}))
\vert
\leq 
K
.
\end{align}
Assume that 
\begin{align*}
f(X) 
\equiv
\log
\int_\Sigma
\ee^{
X(\sigma)
} 
\dd \mu(\sigma)
<
\infty
.
\end{align*}
Then 
\begin{align*}
\P
\left\{
\vert
f(X)
-
\E[f(X)]
\vert
\geq t
\right\}
\leq 
2\exp\left(-\frac{t^2}{4K}\right)
.
\end{align*}
\end{proposition}
\begin{remark}
An analogous result was given in a somewhat more specialised case
in \cite{panchenko-free-energy-generalized-sk-2005}. 
\end{remark}
\begin{proof}
This is an adaptation of the proof of 
\cite[Theorem~2.2.4]{TalagrandSpinGlassesBook2003}. 
We can not apply the comparison Theorem
\ref{prp:as2:gaussian-comparison-of-free-energy} directly, so we resort to
the basic interpolation argument as stated in
Proposition~\ref{basic-facts-generic-gaussian-integration-by-parts}. For
$j=1,2$, let the processes $X_{j}(\cdot)$ be the two independent copies of the
process $X(\cdot)$. For $t \in [0;1]$, let 
\begin{align*}
X_{j,t}
\equiv
\sqrt{t}
X_{j}
+
\sqrt{1-t}
X
\end{align*}
and
\begin{align*}
F_j(t)
\equiv
\log
\int_{\Omega}
\exp
\left(
X_{j,t}(\sigma)
\right)
\dd \mu(\sigma)
.
\end{align*}
For $s \in \R$, let
\begin{align*}
\varphi_s(t)
\equiv
\E
\left[
\exp
\left(
s(F_1-F_2)
\right)
\right]
.
\end{align*}
Hence, differentiation gives
\begin{align}
\label{vektor-sk-dot-varphi-s-t-equals}
\dot{\varphi}_s(t)
=
s
\E
\left[
\exp
\left(
s(F_1-F_2)
\right)
(\dot{F}_1-\dot{F}_2)
\right]
\end{align}
(the dots indicate the derivatives with respect to $t$) and
also
\begin{align}
\label{vektor-sk-dot-f-j-t-equals}
\dot{F}_j(t)
=
&
\frac{1}{2}
\left(
\vphantom{\sum}
\smash{
\int_\Sigma
}
\exp
\left(
X_{j,t}(\sigma)
\right)
\dd \mu(\sigma)
\right)^{-1}
\nonumber
\\
&
\times
\int_\Sigma
\left(
t^{-1/2}
X_{j}(\sigma)
-
(1-t)^{-1/2}
X(\sigma)
\right)
\exp
\left(
X_{j,t}(\sigma)
\right)
\dd \mu(\sigma)
.
\end{align}
Now, we substitute \eqref{vektor-sk-dot-f-j-t-equals} back to
\eqref{vektor-sk-dot-varphi-s-t-equals} and apply Corollary \ref{basic-facts-generic-gaussian-integration-by-parts} to
the result. After some tedious but elementary calculations we get 
\begin{align*}
\dot{\varphi}_s(t)
=
&
s^2
\E
\left[
\exp
\left(
s(F_1-F_2)
\right)
\left(
\int_\Sigma
\exp{
X_{1,t}(\sigma)
}
\dd \mu(\sigma)
\int_\Sigma
\exp{
X_{2,t}(\sigma)
}
\dd \mu(\sigma)
\right)^{-1}
\right.
\\
&
\left.
\int_\Sigma
\cov(X(\sigma^{(1)}),X(\sigma^{(2)}))
\exp{
\left(
X_{1,t}(\sigma^{(1)})
+
X_{2,t}(\sigma^{(2)})
\right)
}
\dd \mu(\sigma^{(1)})
\dd \mu(\sigma^{(2)})
\right]
.
\end{align*}
Thus, thanks to \eqref{vektor-sk-covariance-lemma-assumption}, we obtain
\begin{align*}
\dot{\varphi}_s(t)
\leq
K
s^2
\varphi_s(t)
.
\end{align*}
The conclusion of the theorem follows now exactly as in the proof of 
\cite[Theorem~2.2.4]{TalagrandSpinGlassesBook2003}.

\end{proof}
We now apply this general result to the our model and also to the free
energy-like functional of the GREM-inspired process $A$.
\begin{proposition}
\label{vektor-sk-concentration-sk-free-energy-grem-free-energy}
Suppose $\Sigma \subset B(0,r)$, for $r>0$. For $\Omega \subset \Sigma_N$,
denote
\begin{align*}
P^\text{SK}_{N}(\beta, \Omega)
\equiv
\log
\int
_\Omega
\exp
\left(
\sqrt{N}\beta
X_N(\sigma)
\right)
\dd
\mu^{\otimes N}(\sigma)
,
\end{align*}
and
\begin{align*}
P^\text{GREM}_{N}(\beta, \Omega)
\equiv
\log
\int_{
\Omega
\times
\mathcal{A}
}
\exp
\left(
\beta
\sqrt{2}
\vphantom{\sum}
\smash{
\sum_{
i=1
}^{
N
}
}
\langle
A_i(\alpha)
,
\sigma_i
\rangle
\right)
\dd
\pi_N(\sigma,\alpha)
.
\end{align*}
Then, for all $\Omega \subset \Sigma_N$, we have
\begin{enumerate}
\item 
For any $t > 0$,
\begin{align}
\label{vektor-sk-concentration-for-sk-free-energy}
\P
\left\{
\left\vert
P^\text{SK}_{N}(\beta, \Omega)
-
\E
\left[
P^\text{SK}_{N}(\beta, \Omega)
\right]
\right\vert
>
t
\right\}
\leq
2
\exp
\left(
-
\frac{t^2}{
4
\beta^2 
r^4
N
}
\right)
.
\end{align}
\item
For any $t > 0$,
\begin{align}
\label{vektor-sk-concentration-for-grem-free-energy}
\P
\left\{
\left\vert
P^\text{GREM}_{N}(\beta, \Omega)
-
\E
\left[
P^\text{GREM}_{N}(\beta, \Omega)
\right]
\right\vert
>
t
\right\}
\leq
2
\exp
\left(
-
\frac{t^2}{
8 \beta^2 r^4
N
}
\right)
.
\end{align}
\end{enumerate}
\end{proposition}
\begin{proof}
\begin{enumerate}
\item 
We would like to use Proposition
\ref{vektor-sk-concentration-of-measure-for-free-energies}. By
\eqref{eq:multidimensional-sk:covariance} and the Cauchy-Bouniakovsky-Schwarz
inequality, we have, for all $N \in \N$, $\sigma^{(1)}, \sigma^{(2)} \in
\Sigma_N$, that
\begin{align}
\label{vektor-sk-sk-covariance-norm-bound}
\cov(
X_N(\sigma^{(1)}, \sigma^{(2)})
)
=
\Vert
R_N(\sigma^{(1)}, \sigma^{(2)})
\Vert_\text{F}^2
=
\frac{1}{N^2}
\sum_{
i,j=1
}^{N}
\langle
\sigma^{(1)}_i
,
\sigma^{(1)}_j
\rangle
\langle
\sigma^{(2)}_i
,
\sigma^{(2)}_j
\rangle
\leq
r^4
.
\end{align}
Hence, for all $N \in \N$ and all subsets $\Omega$ of $\Sigma_N$, we obtain
\begin{align*}
\sup_{
\sigma^{(1)},\sigma^{(2)}
\in
\Sigma
}
\vert
\cov 
(X(\sigma^{(1)}),X(\sigma^{(2)}))
\vert
\leq
r^4
.
\end{align*}
Thus \eqref{vektor-sk-concentration-for-sk-free-energy} is proved.
\item
We fix an
arbitrary $N \in \N$,
$\sigma^{(1)},\sigma^{(2)} \in \Sigma_N$, $\alpha^{(1)}, \alpha^{(2)} \in \mathcal{A}$.
We have
\begin{align*}
\cov(
A(\sigma^{(1)},\alpha^{(1)})
,
A(\sigma^{(2)},\alpha^{(2)})
)
&
=
\E
\left[
A(\sigma^{(1)},\alpha^{(1)})
A(\sigma^{(2)},\alpha^{(2)})
\right]
\\
&
=
\sum_{
i=1
}^{N}
\langle
Q(\alpha^{(1)},\alpha^{(2)})
\sigma^{(1)}_i
,
\sigma^{(2)}_i
\rangle
.
\end{align*}
Bound \eqref{vektor-sk-sk-covariance-norm-bound} implies that, for any $U \in \mathcal{U}$, we have 
$
\left\Vert
U
\right\Vert_2
\leq
r^2
$.
Since $Q(\alpha^{(1)},\alpha^{(2)}) \in \mathcal{U}$, we obtain
\begin{align*}
\vert
\langle
Q(\alpha^{(1)},\alpha^{(2)})
\sigma^{(1)}_i
,
\sigma^{(2)}_i
\rangle
\vert
&
\leq
\Vert
Q(\alpha^{(1)},\alpha^{(2)})
\Vert_2
\Vert
\sigma^{(1)}_i
\Vert_2
\Vert
\sigma^{(2)}_i
\Vert_2
\\
&
\leq
\Vert
Q(\alpha^{(1)},\alpha^{(2)})
\Vert_2
r^2
\leq
r^4
.
\end{align*}
\end{enumerate}
Therefore, using Proposition
\ref{vektor-sk-concentration-of-measure-for-free-energies}, we obtain
\eqref{vektor-sk-concentration-for-grem-free-energy}.

\end{proof}
\subsection{Gaussian comparison inequalities for free energy-like functionals}
\label{sec:gaussian-comparison-inequalities-for-free-energy-like-functionals}
We begin by  recalling  well-known integration by parts formula
which is the source of many comparison results for  functionals of Gaussian
processes. 

 Let $F: X \to \R$ be a functional on a linear space $X$.
Given 
$
x \in X
$ and
$
e \in X
$, 
a \emph{directional (G\^ateaux) derivative} of $F$ at $x$ along the
direction $e$ is
\begin{align}
\label{ab.201}
\partial_{x \leadsto e}
F(x)
\equiv 
\partial_t
F(x + t e)
\Big\vert_{t=0}
.
\end{align}
With this notation the following lemma holds.
\begin{lemma}
\label{basic-facts-generic-gaussian-integration-by-parts}
Let 
$
\{
g(i)
\}_{
i \in I
}
$
be a real-valued Gaussian process (the set $I$ is an arbitrary index set), 
and $h$ be some Gaussian random variable.
Define the vector $e \in \R^I$ as 
$
e(i)
\equiv
\E
\left[
h g(i)
\right]
$,
$
i \in I
$.
Let  $F: \R^I \to \R$ such that, for all
$
f \in \R^I
$,
the function
\begin{align}
\label{eq:review:integration-by-parts-growth-condition}
\R
\ni
t
\mapsto
F(f + t e)
\in
\R
\end{align}
is either locally absolute continuous or everywhere differentiable on $\R$.
Moreover, assume that the random variables $
h F(g)
$ 
and 
$
\partial_{g \leadsto e}
F(g)
$ are in 
$
L^1
$.

Then
\begin{align}
\label{gauss-integr-by-parts}
\E[h F(g)]
=
\E
\left[
\partial_{g \leadsto e}
F(g)
\right]
.
\end{align}
\end{lemma}
The previous proposition coincides with
\citep[Lemma~4]{panchenko-free-energy-generalized-sk-2005} (modulo the
differentiability condition on \eqref{eq:review:integration-by-parts-growth-condition} and the integrability assumptions which are needed, e.g., for 
\cite[Theorem~5.1.2]{Bogachev1998}).

The following proposition connects the computation of the derivative of the
free energy with respect to the parameter that linearly occurs in the Hamiltonian
with a certain Gibbs average for a replicated system.

\begin{proposition}
\label{free-energy-beta-derivative}
Consider a Polish measure space $(\Sigma,\mathfrak{S})$ and a random measure
$\mu$ on it. Let 
$
X = 
\{
X(\sigma)
\}_{
\sigma \in \Sigma
}
$
and 
$
Y
\equiv
\{
Y(\sigma)
\}_{
\sigma \in \Sigma
}
$ be two independent Gaussian real-valued processes.
For 
$
u \in \R
$, 
we define
\begin{align*}
H_u(\sigma) \equiv u X(\sigma)+Y(\sigma)
.
\end{align*}
Assume that, for all 
$
u \in [a,b] \Subset \R
$, 
we have
\begin{align*}
\int
\exp
\left(
H_u(\sigma)
\right)
\dd 
\mu(\sigma)
<
\infty
,
\int
X(\sigma)
\exp
\left(
H_u(\sigma)
\right)
\dd 
\mu(\sigma)
<
\infty
\end{align*}
almost surely, and also that
\begin{align*}
\E
\left[
\log
\int
\exp
\left(
H_u(\sigma)
\right)
\dd 
\mu(\sigma)
\right]
<
\infty
.
\end{align*}
Then we have
\begin{align*}
\frac{\dd}{\dd u}
&
\E
\left[
\log
\int
\ee^{H_u(\sigma)} 
\dd \mu(\sigma)
\right]
= 
u
\E
\left[
\mathcal{G}(u) \otimes \mathcal{G}(u)
\left[
\var X(\sigma)
-
\E
\left[
X(\sigma)
,
X(\tau)
\right]
\right]
\right]
,
\end{align*}
where 
$
\mathcal{G}(u)
$ 
is the random element of 
$
\mathcal{M}_1(\Sigma)
$
which, for any measurable
$
f:
\Sigma
\to
\R
$
,
satisfies
\begin{align*}
\mathcal{G}(u)
\left[
f
\right]
=
\frac
{
1
}
{
Z(u)
}
\int
f(\sigma)
\exp
\left(
H_u(\sigma)
\right)
\dd \mu(\sigma)
.
\end{align*}
\end{proposition}
\begin{proof}
We write
\begin{align}
\label{d-d-u-pressure}
\frac{\dd}{\dd u}
\log
\int
\ee^{H_u(\sigma)}
\dd \mu(\sigma)
=
\int
X(\sigma)
\frac
{
\ee^{H_u(\sigma)}
}
{
Z_u(\beta)
}
\dd \mu(\sigma)
,
\end{align}
where 
$
Z_u(\beta) 
\equiv 
\int
\ee^{\beta H_u(\sigma)}
\dd \mu(\sigma)
$.
The main ingredient of the proof is the Gaussian integration 
by parts formula. Denote, for 
$
\tau \in \Sigma
$,
$
e(\tau)
\equiv
\E
\left[
X(\sigma)H_u(\tau)
\right]
$. 
By \eqref{gauss-integr-by-parts}, we have
\begin{align}
\label{directional-derivative-h-u}
\E
\left[
X(\sigma)
\frac
{
\ee^{H_u(\sigma)}
}
{
Z_u(\beta)
}
\right]
=
\E
\left[
\partial_X
\left(
\frac
{
\ee^{H_u(\sigma)}
}
{
\int
\ee^{H_u(\tau)}
\dd \mu(\tau)
}
\right)
\left(
X
;
e
\right)
\right]
.
\end{align}
Due to the independence, we have
\begin{align*}
\E
\left[
X(\sigma)H_u(\tau)
\right]
=
u 
\E
\left[
X(\sigma)
,
X(\tau)
\right]
.
\end{align*}
Henceforth, the computation of the directional derivative in 
\eqref{directional-derivative-h-u} amounts to
\begin{align}\label{all.1}\nonumber
&
\frac
{
\partial
}
{
\partial
t
}
\left[
\frac
{
\ee^{H_u(\sigma)+t u \var(\sigma)}
}
{
\int
\ee^{H_u(\tau)+t u \cov(\sigma, \tau)}
\dd \mu(\tau)
}
\right]
\\
\nonumber
&
\quad\quad
=
\left(
\int
\ee^{H_u(\sigma)}
\dd \mu(\sigma)
\right)
^{-2}
\left(
u \var X(\sigma)
\ee^{H_u(\sigma)}
\int
\ee^{H_u(\tau)}
\dd \mu(\tau)
\right.
\\
&
\left.
\quad\quad\quad\quad
-
\ee^{H_u(\sigma)}
\int
u \cov \left[X(\sigma),X(\tau)\right]
\ee^{H_u(\tau)}
\dd \mu(\tau)
\right)
.
\end{align}
Substituting the r.h.s. of \eqref{all.1} into \eqref{d-d-u-pressure},
we obtain the assertion of the proposition.
\end{proof}
The following proposition gives a short differentiation formula, which is useful
in getting comparison results between the (free energy-like) functionals of Gaussian processes.
\begin{proposition}
\label{prp:as2:gaussian-comparison-of-free-energy}
Let $(X(\sigma))_{\sigma\in\Sigma}$, $(Y(\sigma))_{\sigma\in\Sigma}$
be two independent Gaussian processes as before. Set
\begin{align*}
H_t(\sigma)
\equiv
\sqrt{t}X(\sigma)+\sqrt{1-t}Y(\sigma)
.
\end{align*}
Assume that
\begin{align*}
\int \ee^{H_t(\sigma)}
\dd \mu(\sigma) 
< \infty
,
\int 
X(\sigma)
\ee^{H_t(\sigma)}
\dd \mu(\sigma) 
< \infty
,
\\
\int 
Y(\sigma)
\ee^{H_t(\sigma)}
\dd \mu(\sigma) 
< \infty
\end{align*}
almost surely, and also that, for all $t \in [0;1]$,
\begin{align*}
\E
\left[
\log
\int 
\ee^{H_t(\sigma)}
\dd \mu(\sigma) 
\right]
<
\infty
.
\end{align*}
Then we have
\begin{align}
\label{eq:as2-scheme-comparison-formula}
&
\E
\left[
\log
\int \ee^{X(\sigma)} \dd \mu(\sigma) 
\right]
=
\E
\left[
\log
\int \ee^{Y(\sigma)} \dd \mu(\sigma)  
\right]
\nonumber
\\
&
-
\frac{1}{2}
\int_0^1
\mathcal{G}(t) \otimes \mathcal{G}(t)
\left[
\left(
\var X(\sigma^{(1)}) 
- 
\var Y(\sigma^{(1)}) 
\right)
\right.
\nonumber
\\
&
\quad
-
\left.
\left(
\cov
\left[
X(\sigma^{(1)})
,
X(\sigma^{(2)})
\right]
- 
\cov
\left[
Y(\sigma^{(1)})
,
Y(\sigma^{(2)})
\right]
\right)
\right]
\dd t
,
\end{align}
where 
$
\mathcal{G}(t)
$ 
is the random element of 
$
\mathcal{M}_1(\Sigma)
$
which, for all measurable
$
f: \Sigma \to \R
$, 
satisfies
\begin{align}
\label{eq:multidimensional-sk:abstract-interpolated-gibbs-measure}
\mathcal{G}(t)
\left[
f
\right]
=
\frac{1}{
Z(t)
}
\int_{\Sigma}
f(\sigma)
\exp
\left(
H_t(\sigma)
\right)
\dd \mu(\sigma)
.
\end{align}
\end{proposition}
\begin{proof}
Let us introduce the process
\begin{align*}
W_{u,v}(\sigma)
\equiv 
u X(\sigma)
+
v Y(\sigma)
.
\end{align*}
Hence,
\begin{align}
\label{vektor-sk-h-t-sigma-as-composition}
H_t(\sigma) 
= 
W_{\sqrt{t},\sqrt{1-t}}(\sigma)
.
\end{align}
Thus
\begin{align*}
\frac{\dd}{\dd t}
\E
\left[
\log 
\int \ee^{H_t(\sigma)}\dd \mu(\sigma) 
\right]
&
=
\frac{1}{2}
\left(
\frac{1}{\sqrt{t}}
\frac{\partial}{\partial u}
\E
\left[
\log 
\int \ee^{W_{u,v}(\sigma)}\dd \mu(\sigma) 
\right]
\right.
\\
&
\quad
\quad
\left.
\left.
-
\frac{1}{\sqrt{1-t}}
\frac{\partial}{\partial v}
\E
\left[
\log 
\int \ee^{W_{u,v}(\sigma)}\dd \mu(\sigma) 
\right]
\right)
\right\vert_{u=\sqrt{t},v=\sqrt{1-t}}
.
\end{align*}
Applying Proposition \ref{free-energy-beta-derivative} and 
$
\int_0^1 \cdot \dd t
$
to the previous formula, we conclude the proof.
\end{proof}

\section{Quenched G\"artner-Ellis type LDP}
\label{sec:sk-with-multidimensional-spins:gaertner-ellis}

In this section, we 
derive a quenched LDP under measure concentration assumptions. Theorems~\ref{thm:as2:generic-quenched-ldp-upper-bound} and 
\ref{vektor-sk-generic-quenched-ldp-lower-bound} give the corresponding LDP
upper and lower bounds, respectively. 
The proofs
of the LDP bounds will be adapted to get the proofs of the upper
and lower bounds on the free energy of the SK model with multidimensional spins.
However, they may be of independent interest. 

Note that the existing ``level-2"
quenched large deviation results of Comets \cite{Comets1989} are applicable only to a
certain class of mean-field random Hamiltonians which are required to be
``macroscopic'' functionals of the joint empirical distribution of the random
variables representing the disorder and the independent spin variables. The SK
Hamiltonian can not be represented in such form, since the interaction matrix
consists of i.i.d. random variables. Moreover, it is assumed in \cite{Comets1989} 
that the Hamiltonian has the form 
$
H_N(\sigma) = N
V(\sigma)
$,
where $\{ V(\sigma) \}_{\sigma \in \Sigma_N}$
is a random process taking values in some fixed bounded subset of $\R$. 
Since the Hamiltonian of our model is a Gaussian process, this assumption is
also not satisfied, due to the unboundedness of the Gaussian distribution.

\subsection{Quenched LDP upper bound}
The following assumption will be satisfied for the applications we have in
mind. As is clear from what follows, much weaker concentration functions are
also allowed.
\begin{assumption}
\label{asm:as2:measure-concentration}
Suppose $\{Q_N\}_{N=1}^\infty$ is a sequence of random
measures on a Polish space $(\mathcal{X}, \mathfrak{X})$.
Assume that there exists some $L>0$ such that for any
$Q_N$-measurable set 
$
A 
\subset
\mathcal{X}
$
we have
\begin{align}
\label{eq:as2:subgaussian-concentration}
\P
\left\{
\left\vert
\log
Q_N(A)
-
\E
\left[
\log
Q_N(A)
\right]
\right\vert
>
t
\right\}
\leq
\exp
\left(
-\frac{
t^2
}{
L
N
}
\right)
.
\end{align}
\end{assumption}
Note that Assumption~\ref{asm:as2:measure-concentration} will hold in the cases
we are interested in due to Proposition \ref{vektor-sk-concentration-of-measure-for-free-energies}.

\begin{lemma}
\label{vektor-sk-lemma-the-largest-exponent-wins-quenched}
Suppose $\{Q_N\}_{N=1}^\infty$ is a sequence of random
measures on a Polish space $(\mathcal{X},\mathfrak{X})$ and for
$
\{
A_r
\subset
\mathcal{X}
:
r 
\in 
\{
1,\dots,p
\}
\}
$ 
is a sequence of $Q_N$-measurable sets such that, for some absolute constant
$L > 0$ and some concentration function $\eta_N(t): \R_+ \to \R_+$ with the
property
\begin{align}
\label{eq:as2:concentration-function-decay}
\int_0^{+\infty} \eta_N(t N) \dd t
\xrightarrow[N \uparrow +\infty]{}
0
,
\end{align}
we have
\begin{align}
\label{vektor-sk-exponential-concentration-for-q-n-assumption}
\P
\left\{
\left\vert
\log
Q_N(A_r)
-
\E
\left[
\log
Q_N(A_r)
\right]
\right\vert
>
t
\right\}
\leq
\eta_N(t)
.
\end{align}
Then we have
\begin{align}
\label{vektor-sk-the-largest-exp-wins-ineq}
\lim_{
N
\uparrow
+\infty
}
\frac{1}{N}
\E
\left[
\left\vert
\log
Q_N
\left(
\smash{
\bigcup_{
r=1
}^{
p
}
}
\vphantom{
\bigcup
}
A_r
\right)
-
\max_{
r
\in
\{
1,\dots,p 
\}
}
\E
\left[
\log
Q_N(A_r)
\right]
\right\vert
\right]
=
0
.
\end{align}
\end{lemma}
\begin{remark}
As is easy to extract from Assumption~\ref{asm:as2:measure-concentration}, we
will apply this result in the very pleasant situation, where
\begin{align*}
\gamma_N(t)
=
\exp
\left(
-\frac{
t^2
}{
L
N
}
\right)
.
\end{align*}
However, our subsequent results hold for substantially worse concentration
functions satisfying \eqref{eq:as2:concentration-function-decay}.
\end{remark}

\begin{proof}[Proof of
Lemma~\ref{vektor-sk-lemma-the-largest-exponent-wins-quenched}] First, \eqref{vektor-sk-exponential-concentration-for-q-n-assumption} gives 
\begin{align*}
\P
\left\{
\max_{
r 
\in 
\{ 1,\dots,p \}
}
\left\vert
\log
Q_N(A_r)
-
\E
\left[
\log
Q_N(A_r)
\right] 
\right\vert
\geq
t
\right\}
\leq
p
\eta_N(t)
.
\end{align*}
Since, for 
$
a, b \in \R^p
$, 
the following elementary inequality holds
\begin{align*}
\left\vert
\max_{
r
} 
a_r
-
\max_{
r
} 
b_r
\right\vert
\leq
\max_{
r
}
\vert
a_r
-
b_r
\vert
,
\end{align*}
we get
\begin{align*}
\P
\left\{
\left\vert
\max_{
r 
\in 
\{ 1,\dots,p \}
}
\log
Q_N(A_r)
-
\max_{
r 
\in 
\{ 1,\dots,p \}
}
\E
\left[
\log
Q_N(A_r)
\right] 
\right\vert
\geq
t
\right\}
\leq
p
\eta_N(t N)
.
\end{align*}
The last equation in turn implies that
\begin{align}
\label{as2:upper-bound:max-expect-are-almost-commutative}
\frac{1}{N}
\E
\left[
\left\vert
\max_{
r 
\in 
\{ 1,\dots,p \}
}
\log
Q_N(A_r)
-
\max_{
r 
\in 
\{ 1,\dots,p \}
}
\E
\left[
\log
Q_N(A_r)
\right]
\right\vert
\right]
&
\leq
p
\int_{0}^{
+\infty
}
\eta_N(t N)
\dd t
,
\end{align}
and the r.h.s. of the previous formula vanishes as $N \uparrow +\infty$ due to
\eqref{eq:as2:concentration-function-decay}.
\end{proof}
Let
$
Q_{N}
\in
\mathcal{M}(
\mathcal{X}
)
$,
$
N \in \N
$
be a family of random measures on $(\mathcal{X},\mathfrak{X})$.
Define the Laplace transform 
\begin{align*}
L_{N}(\Lambda)
\equiv
\int_{
\mathcal{X}
}
\ee^{
N
\langle
x,
\Lambda
\rangle
}
\dd
Q_{N}(x)
.
\end{align*}
Suppose that, for all $\Lambda \in \R^d$, we have
\begin{align}
\label{vektor-sk-i-lambda-defined}
I(\Lambda)
\equiv
\lim_{
N \uparrow \infty
}
\frac{1}{N}
\E
\left[
\log
L_N(\Lambda)
\right]
\in
\overline{\R}
=
\R
\cup
\{
-\infty
, 
+\infty
\}
.
\end{align}
Define the Legendre transform
\begin{align}
\label{vector-sk-chap-i-star-definition}
I^{*}
(x)
\equiv
\inf_{
\Lambda
}
\left[
-
\langle
x
,
\Lambda
\rangle
+
I(\Lambda)
\right]
.
\end{align}
Define, for $\delta > 0$,
\begin{align}
\label{eq:as2:i-star-delta}
I^{*}_{
\delta
}
(x)
\equiv
\max
\left\{
I^{*}
(x)
+
\delta
,
-
\frac{1}{\delta}
\right\}
.
\end{align}
\begin{lemma}
\label{vektor-sk-i-star-is-a-good-rate-function}
Suppose
\begin{align}
\label{vektor-sk-zero-in-interior-of-domain-of-I}
0 
\in
\interior
\mathcal{D}(I)
\equiv
\interior
\{
\Lambda
:
I(\Lambda) 
<
+\infty
\}
.
\end{align}
Then
\begin{enumerate}
\item 
The mapping $I^*(\cdot): \mathcal{X} \to \R$ is upper semi-continuous and
concave.
\item
For all $M>0$, 
\begin{align*}
\text{
$
\{
x
\in
\mathcal{X}
:
I^*(x)
\leq
M
\}
$
is a compact.
}
\end{align*}
\end{enumerate}
\end{lemma}
\begin{proof}
\begin{enumerate}
\item 
Since, for all $\Lambda \in \mathcal{D}(I)$, the linear
mappings 
\begin{align*}
x
\mapsto
-
\langle
\Lambda
,
x
\rangle
+
I(\Lambda)
\end{align*}
are obviously concave, the infimum of this family is upper
semi-continuous and concave.
\item
See, e.g., \cite{denHollanderLDPBook2000} for the proof.
\end{enumerate}

\end{proof}
\begin{theorem}
\label{thm:as2:generic-quenched-ldp-upper-bound}
Suppose that
\begin{enumerate}
\item 
The family $\{Q_N\}$ satisfies condition
\eqref{vektor-sk-the-largest-exp-wins-ineq}.
\item
Condition \eqref{vektor-sk-i-lambda-defined} is satisfied.
\item
Condition \eqref{vektor-sk-zero-in-interior-of-domain-of-I} is satisfied.
\end{enumerate}
Then, for any closed set 
$
\mathcal{V} \subset \R^d
$, 
we have
\begin{align}
\label{vektor-sk-varlimsup-expect-log-q-n-v-leq-sup-i-star}
\varlimsup_{
N
\uparrow
\infty
}
\frac{1}{N}
\E
\left[
\log
Q_{N}
(\mathcal{V})
\right]
\leq
\sup_{
x
\in
\mathcal{V}
}
I^*
(x)
.
\end{align}
\end{theorem}
\begin{proof}
\begin{enumerate}
\item 

Suppose at first that $\mathcal{V}$ is a compact. 

Thanks to \eqref{vector-sk-chap-i-star-definition}, for any $x \in \mathcal{X}$,
there exists $\Lambda(x) \in \mathcal{X}$ such that
\begin{align}
\label{eq:as2:supremum-approximant}
-
\langle
x
,
\Lambda(x)
\rangle
+
I(\Lambda(x))
\leq
I^{*}_{
\delta
}
(x)
.
\end{align}
For any $x \in \mathcal{X}$, 
there exists a neighbourhood $A(x)\subset\mathcal{X}$ of $x$ such that 
\begin{align*}
\sup_{
y
\in
A(x)
}
\langle
y-x
,
\Lambda(x)
\rangle
\leq
\delta
.
\end{align*}
By compactness, the covering 
$
\bigcup_{
x \in \mathcal{Y}
}
A(x)
\supset
\mathcal{V}
$
has the finite subcovering, say
$
\bigcup_{
r=1
}^{
p
}
A(x_r)
\supset
\mathcal{V}
$.
Hence,
\begin{align}
\label{eq:as2:log-p-n-u-upper-bound}
\frac{1}{N}
\log
Q_N(\mathcal{V})
&
\leq
\frac{1}{N}
\log
\left(
\bigcup_{
r=1
}^{
p
}
Q_N(A(x_r))
\right)
.
\end{align}
Applying condition \eqref{vektor-sk-the-largest-exp-wins-ineq}, we get
\begin{align}
\label{as2:upper-bound:proof:max-expect-are-almost-commutative}
\varlimsup_{
N
\uparrow
\infty
}
\frac{1}{N}
\E
\left[
\max_{
r 
\in 
\{ 1,\dots,p \}
}
\log
Q_N(A(x_r))
-
\max_{
r 
\in 
\{ 1,\dots,p \}
}
\E
\left[
\frac{1}{N}
\log
Q_N(A(x_r))
\right]
\right]
\leq
0
.
\end{align}
By the Chebyshev inequality,
\begin{align}
\label{eq:as2:q-n-a-x}
Q_N(A(x)) 
&
\leq
Q_N
\left\{
y
\in
\mathcal{X}
:
\langle
y-x
,
\Lambda(x)
\rangle
\leq
\delta
\right\}
\nonumber
\\
&
\leq
\ee^{
-
\delta
N
}
\int_{
\mathcal{X}
}
\ee^{
N
\langle
y-x
,
\Lambda(x)
\rangle
}
\dd
Q_N(y)
\nonumber
\\
&
=
\ee^{
-
\delta
N
}
\ee^{
-
N
\langle
x
,
\Lambda(x)
\rangle
}
L_N(\Lambda(x))
.
\end{align}
Hence, \eqref{eq:as2:q-n-a-x} together with \eqref{eq:as2:supremum-approximant}
yields
\begin{align}
\label{eq:as2:abstract:log-chebyshev}
\varlimsup_{
N
\uparrow
+\infty
}
\frac{1}{N}
\E
\left[
\log
Q_N(A(x_r))
\right]
&
\leq
\lim_{
N
\uparrow
+\infty
}
\left[
-
\langle
x_r
,
\Lambda(x_r)
\rangle
+
\frac{1}{N}
\log
L_N(\Lambda(x_r))
\right]
-
\delta
\nonumber
\\
&
=
-
\langle
x_r
,
\Lambda(x_r)
\rangle
+
I(\Lambda(x_r))
-
\delta
\nonumber
\\
&
\leq
I
^{*}_{
\delta
}
(x_r)
-
\delta
.
\end{align}
Combining \eqref{eq:as2:log-p-n-u-upper-bound},
\eqref{as2:upper-bound:proof:max-expect-are-almost-commutative}, 
\eqref{eq:as2:abstract:log-chebyshev}, we obtain 
\begin{align*}
\varlimsup_{
N
\uparrow
+\infty
}
\frac{1}{N}
\E
\left[
\log
Q_N(\mathcal{V})
\right]
&
\leq
\max_{
r 
\in 
\{ 1,\dots,p \}
}
I
^{*}_{
\delta
}
(x_r)
-
\delta
\\
&
\leq
\sup_{
x \in \mathcal{V}
}
I
^{*}_{
\delta
}
(x)
-
\delta
.
\end{align*}
Taking $\delta \downarrow +0$ limit, we get the assertion of the theorem.
\item
Let us allow now the set $\mathcal{V}$ to be unbounded.
We first prove that the family 
$
Q_{
N
}
$ is quenched
 exponentially tight. 
For that purpose, let
\begin{align*}
R_N(M)
\equiv
\frac
{1}{
N
}
\E
\left[
\log
Q_{
N
}
(
\mathcal{X}
\setminus
[-M;M]^{
d}
)
\right]
,
\end{align*}
and denote
\begin{align*}
R(M)
\equiv
\varlimsup_{
N
\uparrow
+\infty
}
R_N(M)
.
\end{align*}
We want to prove that
\begin{align}
\label{vektor-sk-lim-r-kappa-m-equals-minus-infty}
\lim_{
M
\uparrow
+\infty
}
R(M)
=
-\infty
.
\end{align}
Fix some 
$
u
\in 
\{1,\dots,d\}
$.
Suppose $\delta_{u,p} \in \{0,1\}$ is the standard Kronecker symbol. Let
$e_{u} \in \R^d$ be an element of the standard basis of
$\R^d$, i.e., for all $p \in \{1,\dots,d\}$, we have
\begin{align*}
(e_{u})
_{p}
\equiv
\delta_{u,p}
.
\end{align*}
Thanks to the Chebyshev inequality, we have
\begin{align}
\label{eq:as2:q-n-x-u-leq--m}
Q_{
N
}
\{
x_{u}
\leq
-M
\}
\leq
\ee^{
-
NM
}
\int_{
\R^d
}
\ee^{
-
N
\langle
x
,
e_{u}
\rangle
}
\dd 
Q_{
N
}
(x)
,
\text{
a.s.
}
\end{align}
Now, we get
\begin{align}
\label{eq:as2:int-e-n-x}
\int_{
\R^d
}
\ee^{
-
N
\langle
x
,
e_{u}
\rangle
}
\dd 
Q_{
N
}
(x)
&
=
\frac{1}{
L_{
N
}
(\Lambda_e)
}
\int_{
\R^d
}
\ee^{
N
\langle
x
,
\Lambda_e
-
e_{u}
\rangle
}
\dd 
Q_{
N
}
(x)
\nonumber
\\
&
=
\frac{
L_N(\Lambda_e-e_{u})
}{
L_N(\Lambda_e)
}
,
\text{
a.s.
}
\end{align}
Hence, combining \eqref{eq:as2:q-n-x-u-leq--m} and \eqref{eq:as2:int-e-n-x}, we
obtain
\begin{align}
\label{vektor-sk-hat-p-lower-tail}
\frac{1}{N}
\E
\left[
\log
Q_{
N
}
\{
x_{u}
\leq
-M
\}
\right]
\leq
-M
+
I_N(\Lambda_e-e_{u})
-
I_N(\Lambda_e)
.
\end{align}
Using the same argument, we also get
\begin{align}
\label{vektor-sk-hat-p-upper-tail}
\frac{1}{N}
\E
\left[
\log
Q_{
N
}
\{
x_{u}
\geq
M
\}
\right]
\leq
-M
+
I_N(\Lambda_e+e_{u})
-
I_N(\Lambda_e)
.
\end{align}
We obviously have
\begin{align}
\label{vektor-sk-r-n-m-leq}
R_N(M)
&
\leq
\frac
{1}{
N
}
\E
\left[
\log
Q_{
N
}
\left(
\smash{
\bigcup_{
u=1
}^{
d
}
}
\vphantom{
\bigcup
}
\left(
\left\{
x_{u}
\leq
-M
\right\}
\cup
\left\{
x_{u}
\geq
M
\right\}
\right)
\right)
\right]
.
\end{align}
Applying condition \eqref{vektor-sk-the-largest-exp-wins-ineq} to
\eqref{vektor-sk-r-n-m-leq}, we get 
\begin{align}
\label{eq:as2:varlimsup-n-to-inf-e-log-q-n}
\varlimsup_{
N
\uparrow
+\infty
}
\frac{1}{N}
\E
\Bigl[
&
\log
Q_N
\Bigl(
\bigcup_{
u=1
}^{
d
}
\Bigl(
\left\{
x_{u}
\leq
-M
\right\}
\cup
\left\{
x_{u}
\geq
M
\right\}
\Bigr)
\Bigr)
\nonumber
\\
&
-
\max_{
u
\in
\left\{
1,\dots,d
\right\}
}
\max
\Bigl\{
\E
[
\log
Q_{
N
}
\left(
\left\{
x_{u}
\leq
-M
\right\}
\right)
]
,
\E
[
\log
Q_{N}
\left(
\left\{
x_{u}
\geq
M
\right\}
\right)
]
\Bigr\}
\Bigr]
\leq
0
.
\end{align}
Applying \eqref{vektor-sk-hat-p-lower-tail} and
\eqref{vektor-sk-hat-p-upper-tail} in
\eqref{eq:as2:varlimsup-n-to-inf-e-log-q-n}, we get
\begin{align}
\label{eq:as2:varlimsup-n-to-infty-e-log-q-n-2}
&
\varlimsup_{
N
\uparrow
+\infty
}
\frac{1}{N}
\E
\left[
\log
Q_N
\left(
\smash{
\bigcup_{
u=1
}^{
d
}
}
\vphantom{
\bigcup
}
\left(
\left\{
x_{u}
\leq
-M
\right\}
\cup
\left\{
x_{u}
\geq
M
\right\}
\right)
\right)
\right]
\nonumber
\\
&
\leq
-
M
-
I(\Lambda_e)
+
\max_{
u 
\in
\left\{
1,\dots,d
\right\}
}
\max
\left\{
I(\Lambda_e-e_{u})
,
I(\Lambda_e+e_{u})
\right\}
.
\end{align}
The bound \eqref{eq:as2:varlimsup-n-to-infty-e-log-q-n-2}
assures \eqref{vektor-sk-lim-r-kappa-m-equals-minus-infty}.
Now, since we have (with the help of \eqref{vektor-sk-the-largest-exp-wins-ineq}
and \eqref{vektor-sk-varlimsup-expect-log-q-n-v-leq-sup-i-star})
\begin{align}
\label{eq:as2:varlimsup-n-to-infty-e-log-q-n-3}
\varlimsup_{
N
\uparrow
+\infty
}
\frac{1}{N}
\E
\left[
\log
Q_{
N
}
(\mathcal{V})
\right]
&
\leq
\varlimsup_{
N
\uparrow
+\infty
}
\frac{1}{N}
\E
\left[
\log
Q_{
N
}
(
(
\mathcal{V}
\cap
[-M;M]^d
)
\cup
(
\mathcal{X}
\setminus
[-M;M]^d
)
)
\right]
\nonumber
\\
&
\leq
\max
\left\{
\smash{
\sup_{
x
\in
(
\mathcal{V}
\cap
[-M;M]^{d}
)
}
}
\vphantom{\sup_A}
I^*
(x)
,
R(M)
\right\}
,
\end{align}
the assertion of the theorem follows from
\eqref{vektor-sk-lim-r-kappa-m-equals-minus-infty} by taking the 
$
\varlimsup_{
M
\uparrow
+\infty
}
$
in the bound \eqref{eq:as2:varlimsup-n-to-infty-e-log-q-n-3}.
\end{enumerate}

\end{proof}

\subsection{Quenched LDP lower bound}
Suppose that, for some  $\Lambda \in \R^d$ and all $N \in \N$, we have 
\begin{align*}
\int_{
\mathcal{X}
}
\ee^{
N
\langle
y
,
\Lambda
\rangle
}
\dd
Q_N(y)
<
+\infty
.
\end{align*}
Let 
$
\widetilde{Q}_{
N
,
\Lambda
} 
\in
\mathcal{M}
(\mathcal{X})
$
be the random measure defined by
\begin{align}
\label{vektor-sk-tilde-q-n-lambda-measure-definition}
\widetilde{Q}_{
N
,
\Lambda
}
(A) 
=
\int_{
A
}
\ee^{
N
\langle
y
,
\Lambda
\rangle
}
\dd
Q_N(y)
,
\end{align}
for any $Q_N$ measurable $A \subset \mathcal{X}$. 
\begin{lemma}
\label{vektor-sk-log-laplace-trasform-concentration}
Suppose the family of random measures $Q_N$ satisfies the following
assumptions.
\begin{enumerate}
\item
\label{itm:as2:measure-concentration}
Measure concentration. For all $N \in \N$, there exists some $L>0$ and $\eta_N:
\R_+ \to \R_+$ such that, for any $Q_N$-measurable set $
A 
\subset
\mathcal{X}
$,
we have
\begin{align*}
\P
\left\{
\left\vert
\log
Q_N(A)
-
\E
\left[
\log
Q_N(A)
\right]
\right\vert
>
t
\right\}
\leq
\eta_N(t)
.
\end{align*}
Assume, in addition, that, for some $p > 0$, the concentration function
satisfies
\begin{align}
\label{eq:as2:concentration-function-decay-lower-bound}
N^p
\int_0^{+\infty}
\eta_N(N t)
\dd t
\xrightarrow[N \uparrow +\infty]{}
0
.
\end{align}
\item
\label{itm:as2:tails-decay-condition}
Tails decay condition. Let
\begin{align*}
C(M) 
\equiv
\{
x
\in
\mathcal{X}
:
\Vert
x
\Vert
<
M
\}
.
\end{align*}
There exists $p \in \N$ such that
\begin{align}
\label{eq:as2:lower-bound-quantitative-tightness-like-assumption}
\lim_{
K
\uparrow
+\infty
}
\varlimsup_{
N
\uparrow
\infty
}
\int_{
0
}^{
+\infty
}
\P
\left\{
\frac{1}{N}
\log
\widetilde{Q}_{
N
,
\Lambda
}(
\mathcal{X}
\setminus
C(N^p)
) 
>
-K
+
t
\right\}
\dd t
=
0
.
\end{align}
\item
Non-degeneracy. The family of the sets
$
\left\{
B_j
\subset
\mathcal{X}
:
j 
\in
\{
1,\dots,q
\}
\right\}
$
satisfies the following condition
\begin{align}
\label{vektor-sk-varliminf-log-tilde-q-b-j-geq-minus-infty}
\text{
there exists some 
$
j_0 
\in
\{
1,\dots,q
\}$
such that
}
\varliminf_{
N
\uparrow
\infty
}
\frac{1}{N}
\E
\left[
\log
\widetilde{Q}_{
N,\Lambda
} 
(
B_{j_0}
)
\right]
>
-\infty
.
\end{align}
\end{enumerate}
Then, for any $\Lambda \in \R^d$, 
we have
\begin{align}
\label{vektor-sk-varlimsup-expect-max-log-leq-varlimsup-max-expect-log}
\varlimsup_{
N
\uparrow
\infty
}
\frac{1}{N}
\E
\left[
\log
\widetilde{Q}_{
N
,
\Lambda
}
\left(
\vphantom{
\bigcup
}
\smash{
\bigcup_{
j=1
}^{q}
}
B_j
\right) 
-
\max_{
j 
\in
\{
1,\dots,q
\}
}
\E
\left[
\log
\widetilde{Q}_{
N
,
\Lambda
}
(B_j)
\right]
\right]
\leq
0
.
\end{align}
\end{lemma}
\begin{remark}
The polynomial growth choice of $M = M_N \equiv N^p$ made in assumptions
\eqref{eq:as2:lower-bound-quantitative-tightness-like-assumption} and
\eqref{eq:as2:concentration-function-decay-lower-bound} is made for
specificity. Inspecting the following proof, one can easily restate the
conditions \eqref{eq:as2:lower-bound-quantitative-tightness-like-assumption} and
\eqref{eq:as2:concentration-function-decay-lower-bound} for general $M_N$
dependencies. Effectively, the growth rate of $M_N$ is related to the covering
dimension of the Polish space $(\mathcal{X}, \mathfrak{X})$.
\end{remark} 
\begin{proof}[Proof of Lemma~\ref{vektor-sk-log-laplace-trasform-concentration}]
We fix some 
$
j
\in
\{
1,\dots,q
\}
$.
Take an arbitrary $\eps > 0$, $M>0$ and denote
$
J_{
M,\eps
}
\equiv
\Z
\cap
[
-
\Vert
\Lambda
\Vert
M
/
\eps
;
\Vert
\Lambda
\Vert
M
/
\eps
]
$. 
Consider, for 
$
i
\in
J_{M,\eps} 
$, 
the following closed sets 
\begin{align*}
A_{
i,j
}
\equiv
\{
x 
\in
B_j
:
(j-1)
\eps
\leq
\langle
\Lambda
,
x
\rangle
\leq
j
\eps
\}
.
\end{align*}
We get
\begin{align}
\label{vektor-sk-log-tilde-q-cup-b-j-leq}
\frac{1}{N}
\log
\widetilde{Q}_{
N
,
\Lambda
}
\left(
\vphantom{
\bigcup
}
\smash{
\bigcup_{
j=1
}^{q}
}
B_j
\right)
&
\leq
\frac{1}{N}
\log
\widetilde{Q}_{
N
,
\Lambda
}
\left(
\left(
\vphantom{
\bigcup
}
\smash{
\bigcup_{
j=1
}^{q}
}
B_j
\cap
C(M)
\right)
\cup
\left(
\mathcal{X}
\setminus
C(M)
\right)
\right)
\nonumber
\\
&
\leq
\frac{1}{N}
\max
\left\{
\smash{
\max_{
j
\in
\{1,\dots,q\}
}
}
\log
\widetilde{Q}_{
N
,
\Lambda
}
(
B_j
\cap
C(M)
)
,
\right.
\nonumber
\\
&
\quad
\left.
\log
\widetilde{Q}_{
N
,
\Lambda
}
(
\mathcal{X}
\setminus
C(M)
)
\right\}
+
\frac{
\log(q+1)
}{
N
}
.
\end{align}
We have
\begin{align}
\label{vektor-sk-log-tilde-q-b-j-cap-c-m-leq}
\frac{1}{N}
&
\log
\widetilde{Q}_{
N
,
\Lambda
}
(
B_j
\cap
C(M)
)
\leq
\frac{1}{N}
\log
\vphantom{\sum_A}
\left(
\vphantom{\sum}
\smash{
\sum_{
i
\in
J_{M,\eps}
}
}
\ee^{
N
i
\eps
}
Q_N(A_{
i,j
}
)
\right)
\nonumber
\\
&
\leq
\max_{
i
\in
\left\{
1,\dots,p
\right\}
}
\left[
i
\eps
+
\frac{1}{N}
\log
Q_N(A_{
i
,
j
}
)
\right]
+
\frac{
\log
(
\card 
J_{
M,\eps
}
)
}{
N
}
.
\end{align}
Denote
\begin{align*}
\alpha_{
N
}
(\eps)
\equiv
\max_{
j 
\in
\{
1,\dots,q
\}
}
\max_{
i
\in
J_{M,\eps}
}
\left(
i
\eps
+
\frac{1}{N}
\log
Q_N(A_{
i
,
j
}
)
\right)
,
\end{align*}
and
\begin{align*}
\beta_{
N
}
\equiv
\max_{
j 
\in
\{
1,\dots,q
\}
}
\E
\left[
\log
\widetilde{Q}_{
N
,
\Lambda
}
(B_j)
\right]
,
\end{align*}
\begin{align*}
\widetilde{\beta}_{
N
}
(\eps)
\equiv
\max_{
j 
\in
\{
1,\dots,q
\}
}
\E
\left[
\max_{
i
\in
J_{M,\eps}
}
\left(
i
\eps
+
\frac{1}{N}
\log
Q_N(A_{
i
,
j
}
)
\right)
\right]
,
\end{align*}
\begin{align*}
\gamma_{
N
}
(M)
\equiv
\frac{1}{N}
\log
\widetilde{Q}_{
N
,
\Lambda
}
(
\mathcal{X}
\setminus
C(M)
)
.
\end{align*}
We also have
\begin{align}
\label{vektor-sk-log-tilde-q-b-j-geq}
\frac{1}{N}
\log
\widetilde{Q}_{
N
,
\Lambda
}
(
B_j
)
&
\geq
\frac{1}{N}
\log
\widetilde{Q}_{
N
,
\Lambda
}
(
B_j
\cap
C(M)
)
\nonumber
\\
&
\geq
\max_{
i
\in
J_{M,\eps}
}
\left[
(i-1)
\eps
+
\frac{1}{N}
\log
Q_N(A_{
i
,
j
}
)
\right]
\nonumber
\\
&
=
\max_{
i
\in
J_{M,\eps}
}
\left[
i
\eps
+
\frac{1}{N}
\log
Q_N(A_{
i
,
j
}
)
\right]
-
\eps
.
\end{align}
Due to condition \eqref{itm:as2:measure-concentration}, we have
\begin{align}
\label{vektor-sk-prob-abs-alpha-minus-tilde-beta-abs-is-small}
\P
\left\{
\left\vert
\alpha_{
N
}
(\eps)
-
\widetilde{\beta}_{
N
}
(\eps)
\right\vert
>
t
\right\}
\leq
\eta_N(t N)
q
\card 
J_{
M,\eps
}
.
\end{align}
We put 
$
M
\equiv
M_N
\equiv
N^p
$,
and we get
\begin{align}
\label{vektor-sk-card-j-m-eps-leq}
\card 
J_{
M,\eps
}
&
\leq
2
\Vert
\Lambda
\Vert
M
/
\eps
+
1
\nonumber
\\
&
\leq
2
\Vert
\Lambda
\Vert
N^p
/
\eps
+
1
.
\end{align}
Let
\begin{align*}
X_N(M,\eps)
\equiv
\max
\{
\gamma_N
(M)
,
\alpha_N(\eps)
\}
-
\beta_N
,
\end{align*}
then we have
\begin{align}
\label{vektor-sk-prob-x-n-k-eps-greater-t}
\P
\{
X_N(K,\eps) 
> 
t
\}
&
\leq
\P
\{
\gamma_N
(M)
>
\beta_N
+
t
\}
+
\P
\{
\alpha_N(\eps)
>
\beta_N
+
t
\}
.
\end{align}
Due to property \eqref{vektor-sk-varliminf-log-tilde-q-b-j-geq-minus-infty},
there exists $K>0$ such that we have 
\begin{align}
\label{vektor-sk-prob-gamma-bigger-beta-plus-t}
\P
\{
\gamma_N
(M)
>
\beta_N
+
t
\}
\leq
\P
\{
\gamma_N
(M)
>
-K
+
t
\}
.
\end{align}
Thanks to \eqref{vektor-sk-log-tilde-q-b-j-geq}, we have
\begin{align}
\label{eq:as2:p-alpha-n-eps-greater-beta-n-plus-t}
\P
\{
\alpha_N(\eps)
>
\beta_N
+
t
\}
\leq
\P
\{
\alpha_N(\eps)
>
\widetilde{\beta}_N
(\eps)
+
t
-\eps
\}
.
\end{align}
For $t > \eps$, we apply
\eqref{vektor-sk-prob-abs-alpha-minus-tilde-beta-abs-is-small} and
\eqref{vektor-sk-card-j-m-eps-leq} to
\eqref{eq:as2:p-alpha-n-eps-greater-beta-n-plus-t} to obtain
\begin{align}
\label{vektor-sk-prob-alpha-n-eps-bigger-beta-n-plus-t-leq}
\P
\{
\alpha_N(\eps)
>
\beta_N
+
t
\}
\leq
\left(
2
\Vert
\Lambda
\Vert
N^p
/
\eps
+
1
\right)
q
\eta_N(t N)
.
\end{align}
Combining \eqref{vektor-sk-log-tilde-q-cup-b-j-leq} and
\eqref{vektor-sk-log-tilde-q-b-j-cap-c-m-leq}, we get 
\begin{align}
\label{vektor-sk-expect-log-tilde-q-union-b-j-max-epect-log-tilde-q-b-j-leq}
\E
&
\left[
\log
\widetilde{Q}_{
N
,
\Lambda
}
\left(
\vphantom{
\bigcup
}
\smash{
\bigcup_{
j=1
}^{q}
}
B_j
\right) 
-
\max_{
j 
\in
\{
1,\dots,q
\}
}
\E
\left[
\log
\widetilde{Q}_{
N
,
\Lambda
}
(B_j)
\right]
\right]
\nonumber
\\
&
\leq
\E
\left[
X_N(M,\eps)
\right]
+
\frac{
\log(q+1)
}{
N
}
+
\frac{
\log
\left(
2
\Vert
\Lambda
\Vert
N^p
/
\eps
+
1
\right)
}{
N
}
.
\end{align}
Now, \eqref{vektor-sk-prob-x-n-k-eps-greater-t},
\eqref{vektor-sk-prob-gamma-bigger-beta-plus-t} and
\eqref{vektor-sk-prob-alpha-n-eps-bigger-beta-n-plus-t-leq} 
imply
\begin{align}
\label{eq:as2:e-x-n-m-eps}
\E
\left[
X_N(M,\eps)
\right]
&
\leq
\int_{0}^{
+\infty
}
\P
\left\{
X_N(M,\eps)
>
t
\right\}
\dd t
\nonumber
\\
&
\leq
\int_{\eps}^{
+\infty
}
\P
\left\{
X_N(M,\eps)
>
t
\right\}
\dd t
+
\eps
\nonumber
\\
&
\leq
\int_{\eps}^{
+\infty
}
\P
\{
\gamma_N
(M)
>
-K
+
t
\}
\dd t
\nonumber
\\
&
\quad
+
\left(
2
\Vert
\Lambda
\Vert
N^p
/
\eps
+
1
\right)
q
\int_{\eps}^{
+\infty
}
\eta_N(t N)
\dd t
+
\eps
.
\end{align}
Therefore, taking sequentially 
$ 
\varlimsup_{
N 
\uparrow
+\infty
}
$,
$
\lim_{
K 
\uparrow
+\infty
}
$
and
$
\lim_{
\eps
\uparrow
+
0
}
$
in \eqref{eq:as2:e-x-n-m-eps}, and using
\eqref{eq:as2:concentration-function-decay-lower-bound}, we arrive at
\begin{align}
\label{eq:as2:varlimsup-n-to-infty-e-x-n-n-eps}
\varlimsup_{
N
\uparrow
\infty
}
\E
\left[
X_N(M,\eps)
\right]
\leq
0
.
\end{align}
Bound \eqref{eq:as2:varlimsup-n-to-infty-e-x-n-n-eps} together with
\eqref{vektor-sk-expect-log-tilde-q-union-b-j-max-epect-log-tilde-q-b-j-leq}
implies the assertion of the lemma.

\end{proof}
Let $
\hat{Q}_{
N
,
\Lambda
}
$
be the (random) probability measure defined by
\begin{align*}
\hat{Q}_{
N
,
\Lambda
}
\equiv
\frac{
\widetilde{Q}_N
}{
L_N
(
\Lambda
)
}
.
\end{align*}
\begin{lemma}
\label{vektor-sk-lemma-hat-q-tightness-like-stuff}
Suppose that the measure $Q_N$ satisfies the assumptions of the previous lemma.

Then \eqref{vektor-sk-varlimsup-expect-max-log-leq-varlimsup-max-expect-log}
is valid also for 
$
\hat{Q}_{
N
,
\Lambda
}
$.
\end{lemma}
\begin{proof}
Similar to the one of the previous lemma. 
\end{proof}
\begin{remark}
Recall that a point $x \in \mathcal{X}$ is called an exposed point of the 
concave mapping $I^*$ if there exists $\Lambda \in \R^d$ such that, for all 
$
y 
\in
\mathcal{X}
\setminus
\{
x
\}
$,
we have
\begin{align}
\label{vektor-sk-exposed-point-definition}
I^*(y)
-
I^*(x)
<
\langle
y-x
,
\Lambda
\rangle
.
\end{align}
\end{remark}
\begin{theorem}
\label{vektor-sk-generic-quenched-ldp-lower-bound}
Suppose 
\begin{enumerate}
\item 
The
family 
$
\left\{
Q_N
:
N \in \N
\right\}
\subset
\mathcal{M}(\R^d)
$
satisfies the assumptions of Lemma
\ref{vektor-sk-log-laplace-trasform-concentration}.
\item
$
\mathcal{G} 
\subset 
\mathcal{X}
$ is an open set.
\item
$
\emptyset
\neq
\mathcal{E}(I^*)
\subset 
\mathcal{D}(I^*)
$ 
is the set of the exposed points of the mapping $I^*$.
\item
Condition \eqref{vektor-sk-zero-in-interior-of-domain-of-I} is satisfied.
\end{enumerate}
Then
\begin{align}
\label{eq:as2:vektor-sk-generic-lower-bound}
\varliminf_{
N
\uparrow
+\infty
}
\frac{1}{N}
\E
\left[
\log
Q_N
(
\mathcal{G}
\cap
\mathcal{E}
)
\right]
\geq
\sup_{
x
\in
\mathcal{G}
}
I^*(x)
.
\end{align}
\end{theorem}
\begin{proof}
Let $B(x,\eps)$ be a ball of radius $\eps>0$ around
some arbitrary $x \in \mathcal{X}$. It suffices to prove that
\begin{align}
\label{vektor-sk-generic-lim-eps-0-N-infty}
\lim_{
\eps\downarrow+0
}
\varliminf_{
N\uparrow\infty
}
\frac{1}{N}
\E
\left[
\log
Q_N(B(x,\eps))
\right]
\geq
I^*(x)
.
\end{align}
Indeed, since we have
\begin{align}
\label{eq:as2:q-n-g}
Q_N(\mathcal{G})
\geq
Q_N(B(x,\eps))
,
\end{align}
applying $\frac{1}{N}\log(\cdot)$, taking the expectation, taking
$\varliminf_{N\uparrow+\infty}$,
$\eps\downarrow+0$ and taking the supremum over $x\in\mathcal{G}$ in
\eqref{eq:as2:q-n-g}, we get \eqref{eq:as2:vektor-sk-generic-lower-bound}.

Take any 
$
x 
\in 
\mathcal{G}
\cap
\mathcal{E}
$. 
Then we can find the corresponding vector
$
\Lambda_e 
= 
\Lambda_e(x) 
\in 
\R^d
$
orthogonal to the exposing hyperplane at the point $x$, as in
 \eqref{vektor-sk-exposed-point-definition}. Define the new (``tilted'') random
 probability measure $\hat{Q}_{N}$ on $\R^d$ by demanding that
\begin{align}
\label{vektor-sk-generic-hat-q}
\frac{
\dd
\hat{Q}_{N}
}
{
\dd
Q_{N}
}
(y)
=
\frac{1}{
L_N(\Lambda_e)
}
\ee^{
N
\langle
y,
\Lambda_e
\rangle
}
.
\end{align}
Moreover, we have 
\begin{align*}
\frac{1}{N}
\E
\left[
\log
Q_N
(
B(x,\eps)
)
\right]
&
=
\frac{1}{N}
\E
\left[
\log
\vphantom{\int}
\smash{
\int_{
B(x,\eps)
}
}
\dd 
Q_N(y)
\right]
\\
&
=
\frac{1}{N}
\E
\left[
\log
L_N(\Lambda_e)
\right]
+
\frac{1}{N}
\E
\left[
\vphantom{\int}
\smash{
\int_{
B(x,\eps)
}
}
\ee^{
-
N
\langle
y
,
\Lambda_e
\rangle
}
\dd
\hat{Q}_N(y)
\right]
\\
&
\geq
\frac{1}{N}
\E
\left[
\log
L_N(\Lambda_e)
\right]
-
\langle
x
,
\Lambda_e
\rangle
- 
\eps \Vert \Lambda_e \Vert_2
+
\frac{1}{N}
\E
\left[
\log
\hat{Q}_N(B(x,\eps))
\right]
.
\end{align*}
Hence,
\begin{align*}
\lim_{
\eps\downarrow+0
}
\varliminf_{
N\uparrow\infty
}
\frac{1}{N}
\E
\left[
\log
Q_N(B(x,\eps))
\right]
\geq
\left[
-
\langle
x
,
\Lambda_e
\rangle
+
I(\Lambda_e)
\right]
+
\lim_{
\eps\downarrow+0
}
\varliminf_{
N\uparrow\infty
}
\frac{1}{N}
\E
\left[
\log
\hat{Q}_N(B(x,\eps))
\right]
.
\end{align*}
Since we have
\begin{align*}
-
\langle
x
,
\Lambda_e
\rangle
+
I(\Lambda_e)
\geq
I^*(x)
,
\end{align*}
in order to show \eqref{vektor-sk-generic-lim-eps-0-N-infty} it
remains to prove that 
\begin{align}
\label{vektor-sk-generic-lim-log-hat-p-b-eps-equals-0}
\lim_{
\eps\downarrow+0
}
\varliminf_{
N\uparrow\infty
}
\frac{1}{N}
\E
\left[
\log
\hat{Q}_N(B(x,\eps))
\right]
=
0
.
\end{align}
The Laplace transform of $\hat{Q}_{N}$ is
\begin{align*}
\hat{L}_{N}(\Lambda)
=
\frac{
L_{
N
}
(\Lambda+\Lambda_e)
}
{
L_{
N
}
(\Lambda_e)
}
.
\end{align*}
Hence, we arrive at
\begin{align*}
\hat{I}(\Lambda)
=
I(\Lambda+\Lambda_e)
-
I(\Lambda_e)
.
\end{align*}
Moreover, we have
\begin{align}
\label{vektor-sk-generic-hat-i-star-u-equals}
\hat{I}^*(x)
&
=
I^*(x)
+
\langle
x
,
\Lambda_e
\rangle
-
I(\Lambda_e)
.
\end{align}
By the assumptions of the theorem, the family $Q_N$ satisfies the assumptions of
Lemma \ref{vektor-sk-log-laplace-trasform-concentration}. Hence, due to Lemma
\ref{vektor-sk-lemma-hat-q-tightness-like-stuff}, the family $\hat{Q}_N$
satisfies \eqref{vektor-sk-the-largest-exp-wins-ineq}. Thus we can apply
Theorem \ref{thm:as2:generic-quenched-ldp-upper-bound} to obtain 
\begin{align}
\label{vektor-sk-generic-log-p-kappa-leq}
\varlimsup_{
N\uparrow+\infty
}
\frac{1}{N}
\E
\left[
\log
\hat{Q}_{
N
}
(
\R^d
\setminus
B(U,\eps)
)
\right]
\leq
\sup_{
y
\in
\mathcal{U}
\setminus
B(x,\eps)
}
\hat{I}^*(y)
.
\end{align}
Lemma \ref{vektor-sk-i-star-is-a-good-rate-function} implies that there exists
some $x_0 \in \mathcal{X} \setminus B(x,\eps)$ (note that $x_0 \neq x$) such that 
\begin{align*}
\sup_{
y
\in
\mathcal{X}
\setminus
B(x,\eps)
}
\hat{I}_{
}^*(y)
=
\hat{I}_{
}^*(x_0)
.
\end{align*}
Since $\Lambda_e$ is an exposing hyperplane, using
\eqref{vektor-sk-generic-hat-i-star-u-equals}, we get 
\begin{align}
\label{vektor-sk-generic-i-kappa-star-u0-is-negative}
\hat{I}^*(x_0)
&
=
I^*(x_0)
+
\langle
x_0
,
\Lambda_e
\rangle
-
I(\Lambda_e)
\nonumber
\\
&
\leq
[
I^*(x_0)
+
\langle
x_0
,
\Lambda_e
\rangle
]
-
[
I^*(x)
+
\langle
x
,
\Lambda_e
\rangle
]
<
0
,
\end{align}
and hence, combining \eqref{vektor-sk-generic-log-p-kappa-leq} and
\eqref{vektor-sk-generic-i-kappa-star-u0-is-negative}, we get 
\begin{align*}
\varlimsup_{
N\uparrow+\infty
}
\frac{1}{N}
\E
\left[
\log
\hat{Q
}_{
N
}(
\R^d
\setminus 
B(x,\eps)
)
\right]
<
0
.
\end{align*}
Therefore, due to the concentration of measure, we have almost surely
\begin{align*}
\varlimsup_{
N\uparrow+\infty
}
\frac{1}{N}
\log
\hat{Q}_{
N
}
(
\R^d
\setminus 
B(x,\eps)
)
<
0
\end{align*}
which implies that, for all $\eps > 0$, we have almost surely
\begin{align*}
\lim_{
N\uparrow+\infty
}
\hat{Q}_{
N
}
(
\R^d
\setminus 
B(x,\eps)
)
=
0
,
\end{align*}
and \eqref{vektor-sk-generic-lim-log-hat-p-b-eps-equals-0} follows by yet another
application of the concentration of measure. 

\end{proof}

\begin{corollary}
\label{vektor-sk-generic-quenched-ldp-lower-bound-strict-convexity-everywhere}
Suppose that in addition to the assumptions of previous Theorem
\ref{vektor-sk-generic-quenched-ldp-lower-bound} we
have
\begin{enumerate}
\item 
\label{itm:as2:differentiability-of-the-rate-function}
$I(\cdot)$ is differentiable on $\interior \mathcal{D}(I)$.
\item
\label{itm:as2:regularity-for-the-ldp-lower-bound}
Either $\mathcal{D}(I)=\mathcal{X}$ or 
\begin{align*}
\lim_{
\Lambda
\to
\partial
\mathcal{D}(I)
}
\Vert
\nabla
I(\Lambda)
\Vert
=
+\infty
.
\end{align*}
\end{enumerate}
Then 
$
\mathcal{E}(I^*) 
= 
\R^d
$,
consequently
\begin{align*}
\varliminf_{
N
\uparrow
+\infty
}
\frac{1}{N}
\E
\left[
\log
Q_N
(
\mathcal{G}
)
\right]
\geq
\sup_{
x
\in
\mathcal{G}
}
I^*(x)
.
\end{align*}
\end{corollary}
\begin{proof}
The proof is the same as in the classical Gärtner-Ellis theorem (see, e.g.,
\cite{denHollanderLDPBook2000}).

\end{proof}

\section{The Aizenman-Sims-Starr comparison scheme}
\label{sec:as2:free-energy-upper-and-lower-bounds}

In this section, we shall extend the
$\text{AS}^2$ scheme to the case of the SK model with multidimensional
spins and prove 
Theorems~\ref{thm:as2:pressure-upper-bound} and
\ref{thm:as2:pressure-lower-bound}, as stated in the introduction. We use the
Gaussian comparison results of
Section~\ref{sec:gaussian-comparison-inequalities-for-free-energy-like-functionals}
in the spirit of $\text{AS}^2$ scheme in order to relate the free energy of
the SK model with multidimensional spins with the free energy of a certain
GREM-inspired model. Comparing to \cite{AizenmanSimsStarr2003}, due to more intricate nature of spin configuration
space, some new effects occur. In particular, the remainder term of the
Gaussian comparison non-trivially depends on the variances and covariances of the Hamiltonians
under comparison. To deal with this obstacle, we use the idea of localisation
to the configurations having a given overlap (cf.
\eqref{eq:as2:local-free-energy}). This idea is formalised by adapting the
proofs of the quenched Gärtner-Ellis type LDP obtained in Section~\ref{sec:sk-with-multidimensional-spins:gaertner-ellis}.

\subsection{Naive comparison scheme}
We start by recalling the basic principles of the $\text{AS}^2$ comparison
scheme (see, e.g., \cite[Chapter 11]{BovierBook2006}). 
It is a simple idea to get the comparison inequalities
by adding some additional structure into the model. However, the way the
additional structure is attached to the model might be suggested by the
model itself. Later on we shall encounter a real-world use of this
trick. Let $(\Sigma,\mathfrak{S})$ and $(\mathcal{A},\mathfrak{A})$ be Polish
spaces equipped with measures $\mu$ and $\xi$, respectively. Furthermore, let
\begin{align*}
X 
\equiv 
\{
X(\sigma)
\}_{
\sigma \in \Sigma
}
,
A 
\equiv 
\{
A(\sigma,\alpha)
\}_{\substack{
\sigma \in \Sigma,\\
\alpha \in \mathcal{A}
}}
,
B
\equiv 
\{
B(\sigma)
\}_{
\alpha \in \mathcal{A}
}
\end{align*}
be independent real-valued Gaussian processes. 
Define the  \emph{comparison functional}
\begin{align}
\label{eq:as2:comparison-functional}
\Phi[C]
\equiv
\E\left[
\log
\int_{
\Sigma\times\mathcal{A}
}
\ee
^{C(\sigma,\alpha)}
\dd
\left(
\mu
\otimes
\xi
\right)
(
\sigma,
\alpha
)
\right]
,
\end{align}
where 
$
C \equiv
\{
C(\sigma,\alpha)
\}_{
\substack{
\sigma \in \Sigma
\\
\alpha \in \mathcal{A}
}
}
$ is a suitable real-valued Gaussian process. Theorem~4.1 of
\cite{AizenmanSimsStarr2006} is easily understood as an example of the following observation.
Suppose $\Phi[X]$ is somehow hard to compute directly, but $\Phi[A]$ and $\Phi[B]$ are manageable.
We always have the following additivity property
\begin{align}\label{all.2}
\Phi
\left[
X
+
B
\right]
=
\Phi
\left[
X
\right]
+
\Phi
\left[
B
\right]
.
\end{align}
Assume now that
\begin{align}\label{all.3}
\Phi
\left[
X
+
B
\right]
\leq
\Phi
\left[
A
\right]
\end{align}
which we can obtain, e.g., from
Proposition~\ref{prp:as2:gaussian-comparison-of-free-energy}.
Combining \eqref{all.2} and \eqref{all.3}, we get the bound
\begin{align}\label{all.4}
\Phi
\left[
X
\right]
\leq
\Phi
\left[
A
\right]
-
\Phi
\left[
B
\right]
.
\end{align}

\subsection{Free energy upper bound}
\label{vektor-sk-free-energy-upper-bound-multidim-sk-model}

Let 
$
\mathcal{V}
\subset
\symmetric(d)
$
be an arbitrary Borell set.
\begin{remark}
Note that $\mathcal{U}$ is closed and convex.
\end{remark}

Let
\begin{align}
\label{vektor-sk-sigma-n-mathcal-v-definition}
\Sigma_N(\mathcal{V})
&
\equiv
\left\{
\sigma
\in
\Sigma_N
:
R_N(\sigma,\sigma)
\in
\mathcal{V}
\right\}
\nonumber
\\
&
=
\left\{
\sigma
\in
\Sigma_N
:
R_N(\sigma,\sigma)
\in
\mathcal{V}
\cap
\mathcal{U}
\right\}
.
\end{align}
Let us define the \emph{local comparison functional}
$
\Phi_N(x,\mathcal{V})
$
as follows (cf. \eqref{eq:as2:comparison-functional})
\begin{align}
\label{eq:remainder:comparison-functional}
\Phi_N(x,\mathcal{V})
[C]
\equiv
\frac{1}{N}
\E
\left[
\log
\pi_N
\left[
\I_{
\Sigma_N(
\mathcal{V}
)
}
\exp
\left(
\beta
\sqrt{N}
C
\right)
\right]
\right]
,
\end{align}
where 
$
C \equiv
\{
C(\sigma,\alpha)
\}_{
\substack{
\sigma \in \Sigma
\\
\alpha \in \mathcal{A}
}
}
$ 
is a suitable Gaussian process.
Let us consider the following family ($N \in \N$) of random measures on the Borell subsets of $\symmetric(d)$
generated by the SK Hamiltonian,
\begin{align*}
P_N(\mathcal{V})
\equiv
\int_{
\Sigma_N(\mathcal{V})
}
\ee
^{
\beta
\sqrt{N}
X_N(\sigma)
}
\dd
\mu^{
\otimes 
N 
}
(\sigma)
,
\end{align*}
and consider also the following family of the random measures generated by the
Hamiltonian $A(\sigma, \alpha)$ 
\begin{align}
\label{vektor-sk-tilde-p-n-mathcal-v-definition}
\widetilde{P}_{N}
(\mathcal{V})
\equiv
\widetilde{P}^{
x,\mathcal{Q},U
}_N
(\mathcal{V})
\equiv
\int_{
\Sigma_N(\mathcal{V})
\times
\mathcal{A}
}
\exp
\left(
\beta
\sqrt{N}
\vphantom{\sum}
\smash{
\sum_{
i=1
}^{
N
}
}
\langle
A_i(\alpha)
,
\sigma_i
\rangle
\right)
\dd
\pi_N(\sigma,\alpha)
,
\end{align}
where the parameters $\mathcal{Q}$ and
$U$
are taken from the definition of the process $A(\alpha)$ (cf.
\eqref{eq:multidimesnional-sk:monotonicity-of-q-matrix-sequence}). The
vector $x$ defines the random measure $\xi \in \mathcal{M}(\mathcal{A})$ (cf.
\eqref{eq:partition-of-the-unit-interval}), and, hence, also the measure 
$
\pi_N 
\in 
\mathcal{M}
(\Sigma \times \mathcal{A})
$.
\begin{remark}
To lighten the notation, most of the time we shall not indicate explicitly the
dependence of the following quantities on the parameters $x$, $\mathcal{Q}$,
$U$.
\end{remark}
Consider (if it exists) the Laplace
transform of the measure
\eqref{vektor-sk-tilde-p-n-mathcal-v-definition}
\begin{align}
\label{vektor-sk-l-n-lambda-grem-definition}
\widetilde{L}_N(\Lambda)
\equiv
\int_{
\mathcal{U}
}
\ee^{
N
\langle
U,
\Lambda
\rangle
}
\dd
\widetilde{P}_N(U)
.
\end{align}
Let (if it exists)
\begin{align}
\label{vektor-sk-i-lambda-for-grem-definition}
\widetilde{I}(\Lambda)
\equiv
\lim_{
N \uparrow \infty
}
\frac{1}{N}
\E
\left[
\log
\widetilde{L}_N(\Lambda)
\right]
.
\end{align}
Define the following Legendre transform
\begin{align}
\label{eq:as2:sk-upper-bound:legendre-transform:i-star-definition}
\widetilde{I}^{*}
(U)
\equiv
\inf_{
\substack{
x
\in \mathcal{Q}^\prime(1,1)
,
\\
\mathcal{Q}
\in \mathcal{Q}^\prime(U,d)
,
\\
\Lambda 
\in
\symmetric(d)
}
}
\left[
-
\langle
U
,
\Lambda
\rangle
-
\Phi_N(x,\mathcal{V})
[B]
+
\widetilde{I}(\Lambda)
\right]
.
\end{align}
Denote, for $\delta > 0$,
\begin{align*}
\widetilde{I}^{*}_{
\delta
}
(U)
\equiv
\max
\left\{
\widetilde{I}^{*}
(U)
+
\delta
,
-
\frac{1}{\delta}
\right\}
.
\end{align*}
Let
\begin{align}
\label{vektor-sk-sk-free-energy-definition}
p(\mathcal{V})
\equiv
\lim_{
N\uparrow+\infty
}
\frac{1}{N}
\E
\left[
\log
P_N(\mathcal{V})
\right]
.
\end{align}
\begin{remark}
Note that the result of \cite{Guerra-Toninelli-Generalized-SK-2003} assures the
existence of the limit in the previous formula.
\end{remark}
\begin{lemma}
\label{vektor-sk-existence-of-laplace-trasform}
We have
\begin{enumerate}
\item 
The Laplace transform \eqref{vektor-sk-l-n-lambda-grem-definition} exists.
Moreover, for any $\Lambda \in \symmetric(d)$, we have
\begin{align}
\label{eq:as2:the-laplace-transform-for-p}
&
\int_{
\mathcal{V}
}
\ee^{
N
\langle
U,
\Lambda
\rangle
}
\dd
P_N(U)
\nonumber
\\
&
=
\int_{
\Sigma_N(\mathcal{V})
}
\exp
\left(
N
\langle
\Lambda
,
R_N(\sigma,\sigma)
\rangle
+
\beta
\sqrt{N}
X(\sigma)
\right)
\dd
\mu^{
\otimes N
}
(\sigma)
,
\\
&
\int_{
\mathcal{V}
}
\ee^{
N
\langle
U,
\Lambda
\rangle
}
\dd
\widetilde{P}_N(U)
\nonumber
\\
\label{vektor-sk-l-n-lambda-grem-representation}
&
=
\int_{
\Sigma_N(\mathcal{V})
\times
\mathcal{A}
}
\exp
\left(
N
\langle
\Lambda
,
R_N(\sigma,\sigma)
\rangle
+
\beta
\sqrt{N}
\vphantom{\sum}
\smash{
\sum_{
i=1
}^{
N
}
}
\langle
A_i(\alpha)
,
\sigma_i
\rangle
\right)
\dd \pi_N(\sigma,\alpha)
.
\end{align}
\item
The quenched cumulant generating function
\eqref{vektor-sk-i-lambda-for-grem-definition} exists in the $N\uparrow\infty$ limit, for any $\Lambda \in \symmetric(d)$. Moreover,
for all $N \in \N$, we have
\begin{align}
\label{vektor-sk-i-n-grem-equals}
I_N(\Lambda)
\equiv
\frac{1}{N}
\E
\left[
\log
L_N(\Lambda)
\right]
=
X_0(x, \mathcal{Q}, \Lambda, U)
,
\end{align}
that is $I_N(\cdot)$ in fact does not depend on $N$.
\end{enumerate}
\end{lemma}
\begin{proof}
\begin{enumerate}
\item 
We prove \eqref{vektor-sk-l-n-lambda-grem-representation}, the proof of
\eqref{eq:as2:the-laplace-transform-for-p} is similar. Since $\mathcal{U}$ is a compact, it follows that, for arbitrary
$\eps > 0$, there exists the following $\eps$-partition of $\mathcal{U}$
\begin{align*}
\mathcal{N}(\eps) 
=
\left\{
\mathcal{V}_r
\subset
\mathcal{U}
:
r
\in
\left\{
1
,
\dots
,
K
\right\}
\right\}
\end{align*}
such that 
$
\bigcup_{
r
}
\mathcal{V}_r
=
\mathcal{U}
$
,
$
\mathcal{V}_r 
\cap
\mathcal{V}_s
=
\emptyset
$,
$
\diam 
\mathcal{V}_r 
\leq 
\eps
$
and pick some $V_r \in \interior \mathcal{V}_r $,
for all $r \neq s$.

We denote
\begin{align*}
\widetilde{L}_N(\Lambda,\eps)
\equiv
\sum_{
r=1
}
^{
K
}
\ee^{
N
\langle
\Lambda
,
V_r
\rangle
}
\int_{
\Sigma_N(\mathcal{V}_r)
\times
\mathcal{A}
}
\exp
\left(
\beta
\sqrt{N}
\vphantom{\sum}
\smash{
\sum_{
i=1
}^{
N
}
}
\langle
A_i(\alpha)
,
\sigma_i
\rangle
\right)
\dd \pi_N(\sigma,\alpha)
.
\end{align*}
For small enough $\eps$, we have
\begin{align*}
\left(
1 - 2 N \Vert \Lambda \Vert \eps
\right)
\ee^{
N
\langle
\Lambda
,
R_N(\sigma,\sigma)
\rangle
}
\leq
\ee^{
N
\langle
\Lambda
,
U
\rangle
}
\leq
\ee^{
N
\langle
\Lambda
,
R_N(\sigma,\sigma)
\rangle
}
\left(
1 + 2 N \Vert \Lambda \Vert \eps
\right)
.
\end{align*}
Therefore, if we denote 
\begin{align*}
\widehat{L}_{N}(\mathcal{V},\Lambda)
\equiv
\int_{
\Sigma_N(\mathcal{V})
\times
\mathcal{A}
}
\exp
\left(
N
\langle
\Lambda
,
R_N(\sigma,\sigma)
\rangle
+
\beta
\sqrt{N}
\vphantom{\sum}
\smash{
\sum_{
i=1
}^{
N
}
}
\langle
A_i(\alpha)
,
\sigma_i
\rangle
\right)
\dd \pi_N(\sigma,\alpha)
,
\end{align*}
we get
\begin{align*}
\left(
1 - 2 N \Vert \Lambda \Vert \eps
\right)
\sum_{
r=1}^{
K
}
\widehat{L}_{N}(\mathcal{V}_r,\Lambda)
\leq
\widetilde{L}_N(\Lambda, \eps)
\leq
\left(
1 + 2 N \Vert \Lambda \Vert \eps
\right)
\sum_{
r=1}^{
K
}
\widehat{L}_{N}(\mathcal{V}_r,\Lambda)
.
\end{align*}
Hence,
\begin{align}
\label{eq:as2:three-milizman}
\left(
1 - 2 N \Vert \Lambda \Vert \eps
\right)
\widehat{L}_{N}(\mathcal{U},\Lambda)
\leq
\widetilde{L}_N(\Lambda, \eps)
\leq
\left(
1 + 2 N \Vert \Lambda \Vert \eps
\right)
\widehat{L}_{N}(\mathcal{U},\Lambda)
.
\end{align}
Let $\eps \downarrow +0$ in \eqref{eq:as2:three-milizman} and we arrive at
\begin{align*}
\widetilde{L}_N(\Lambda) 
=
\widehat{L}_N(\mathcal{U},\Lambda)
.
\end{align*}
That is, the existence of $L_N(\Lambda)$ and the representation
\eqref{vektor-sk-l-n-lambda-grem-representation} are
proved. 
\item
For all $N \in \N$, we have, by the RPC averaging property (see, e.g.,
\cite[Theorem~5.4]{AizenmanSimsStarr2006} or
Theorem~\ref{thm:remainder:some-limiting-grem-properties}, property
\eqref{the-averaging-property} below), that
\begin{align*}
\frac{1}{N}
\E
\left[
\log
\widetilde{L}_N(\mathcal{U},\Lambda)
\right]
=
\Phi_N(x,\mathcal{U})
\left[
A
+
N
\langle
\Lambda
,
R_N(\sigma,\sigma)
\rangle
\right]
=
X_0(x, \mathcal{Q}, \Lambda, U)
.
\end{align*}

\end{enumerate}
\end{proof}
\begin{proof}[Proof of Theorem~\ref{thm:as2:pressure-upper-bound}]
In essence, the proof follows almost literally the proof of Theorem
\ref{thm:as2:generic-quenched-ldp-upper-bound}. The notable difference is that
we apply the Gaussian comparison inequality
(Proposition~\ref{prp:as2:gaussian-comparison-of-free-energy}) in order to
``compute'' the rate function in a somewhat more explicit way.

Due to \eqref{vektor-sk-sigma-n-mathcal-v-definition}, we can without loss of
generality suppose that $\mathcal{V}$ is compact. For any 
$\delta > 0$
and
$
U
\in 
\mathcal{V}
$, 
by
\eqref{eq:as2:sk-upper-bound:legendre-transform:i-star-definition},
there exists 
$
\Lambda(U,\delta) \in \symmetric(d)
$, 
$
x(U,\delta) \in
\mathcal{Q}^\prime(1,1)
$ 
and 
$
Q(U,\delta)  \in
\mathcal{Q}^\prime(U,d)
$ 
such that
\begin{align}
\label{eq:as2:sk-upper-bound:supremum-approximant}
-
\langle
U
,
\Lambda(U)
\rangle
+
\widetilde{I}(\Lambda(U))
\leq
\widetilde{I}^{*}_{
\delta
}
(U)
.
\end{align}
For any $U \in \mathcal{V}$,
there exists an open neighbourhood  
$
\mathcal{V}(U)
\subset
\symmetric(d)
$
of $U$ such that 
\begin{align*}
\sup_{
V
\in
\mathcal{V}(U)
}
\langle
V-U
,
\Lambda(U)
\rangle
\leq
\delta
.
\end{align*}
Fix some $\eps > 0$.
Without loss of generality, we can suppose that all the neighbourhoods satisfy
additionally the condition $\diam \mathcal{V}(U) \leq \eps$. 
By compactness, the covering
$
\bigcup_{
U \in \mathcal{V}
}
\mathcal{V}(U)
\supset
\mathcal{V}
$
has a finite subcovering, say
$
\bigcup_{
r=1
}^{
p
}
\mathcal{V}(U^{(r)})
\supset
\mathcal{V}
$.
We denote the corresponding to this covering approximants in
\eqref{eq:as2:sk-upper-bound:supremum-approximant} by 
$
\{
x^{(r)} \in \mathcal{Q}^\prime(1,1) 
\}_{r=1}^p
$
and
$
\{
\mathcal{Q}^{(r)} \in \mathcal{Q}^\prime(U^{(r)},d) 
\}_{r=1}^p
$. 
We have
\begin{align}
\label{eq:as2:sk-case:log-p-n-u-upper-bound}
\frac{1}{N}
\log
P_N(\mathcal{V})
&
\leq
\frac{1}{N}
\log
\left(
\vphantom{
\bigcup
}
\smash{
\bigcup_{
r=1
}^{
p
}
}
P_N(\mathcal{V}(U^{(r)}))
\right)
.
\end{align}
Due to the concentration of measure Proposition
\ref{vektor-sk-concentration-sk-free-energy-grem-free-energy}, we can apply Lemma
\ref{vektor-sk-lemma-the-largest-exponent-wins-quenched} and get 
\begin{align}
\label{eq:as2:lim-n-to-infty-e-log-p-n-somethin}
\lim_{
N
\uparrow
+\infty
}
\frac{1}{N}
\E
\left[
\left\vert
\log
P_N
\left(
\smash{
\bigcup_{
r=1
}^{
p
}
}
\vphantom{
\bigcup
}
\mathcal{V}(U^{(r)})
\right)
-
\max_{
r
\in
\{
1,\dots,p 
\}
}
\E
\left[
\log
P_N(\mathcal{V}(U^{(r)}))
\right]
\right\vert
\right]
=
0
.
\end{align}
In fact, since we know that \eqref{vektor-sk-sk-free-energy-definition}
exists, \eqref{eq:as2:lim-n-to-infty-e-log-p-n-somethin} implies that
\begin{align}
\label{vektor-sk-}
\lim_{
N
\uparrow
+\infty
}
\frac{1}{N}
\E
\left[
\log
P_N
\left(
\smash{
\bigcup_{
r=1
}^{
p
}
}
\vphantom{
\bigcup
}
\mathcal{V}(U^{(r)})
\right)
\right]
=
\max_{
r
\in
\{
1,\dots,p 
\}
}
\lim_{
N
\uparrow
+\infty
}
\frac{1}{N}
\E
\left[
\log
P_N(\mathcal{V}(U^{(r)}))
\right]
.
\end{align}
For $U^{(r)}$, $x = x^{(r)}$, $\mathcal{Q} = \mathcal{Q}^{(r)}$, 
Proposition~\ref{prp:as2:gaussian-comparison-of-free-energy} gives
\begin{align}
\label{vektor-sk-comaprison-applied-to-the-restricted-sk}
\frac{1}{N}
\E
\left[
\log
P_{
N
}
(\mathcal{V}(U^{(r)}))
\right]
&
=
\frac{1}{N}
\E
\left[
\log
\widetilde{P}_N(\mathcal{V}(U^{(r)}))
\right]
-
\Phi_N(x,\mathcal{U})[B]
\nonumber
\\
&
\quad
+
\mathcal{R}_N(x^{(r)}, \mathcal{Q}^{(r)}, U^{(r)}, \mathcal{V}(U^{(r)}))
+
\mathcal{O}(\eps)
\nonumber
\\
&
\leq
\frac{1}{N}
\E
\left[
\log
\widetilde{P}_N(\mathcal{V}(U^{(r)}))
\right]
-
\Phi_N(x,\mathcal{U})[B]
+
K \eps
,
\end{align}
where $K>0$ is an absolute constant.

By the Chebyshev inequality and
Lemma~\ref{vektor-sk-existence-of-laplace-trasform}, we have
\begin{align*}
\widetilde{P}_N(\mathcal{V}(U)) 
&
\leq
\widetilde{P}_N
\left\{
V
\in
\mathcal{U}
:
\langle
V-U
,
\Lambda(U)
\rangle
\leq
\delta
\right\}
\nonumber
\\
&
\leq
\ee^{
-
\delta
N
}
\int_{
\mathcal{U}
}
\ee^{
N
\langle
V-U
,
\Lambda(U)
\rangle
}
\dd
\widetilde{P}_N(V)
\nonumber
\\
&
=
\ee^{
-
\delta
N
}
\ee^{
-
N
\langle
U
,
\Lambda(U)
\rangle
}
\widetilde{L}_N(\Lambda(U))
.
\end{align*}
Thus, using \eqref{vektor-sk-comaprison-applied-to-the-restricted-sk} and
\eqref{eq:as2:sk-upper-bound:supremum-approximant}, we get
\begin{align}
\label{vektor-sk-restricted-upper-bound}
\lim_{
N
\uparrow
+\infty
}
\frac{1}{N}
\E
\left[
\log
P_N(\mathcal{V}(U^{(r)}))
\right]
&
\leq
\lim_{
N
\uparrow
+\infty
}
\left[
-
\langle
U^{(r)}
,
\Lambda(U^{(r)})
\rangle
-
\Phi[B]
+
\frac{1}{N}
\log
\widetilde{L}_N(\Lambda(U^{(r)}))
\right]
-
\delta
+
K \eps
\nonumber
\\
&
=
-
\langle
U_r
,
\Lambda(U_r)
\rangle
-
\Phi[B]
+
\widetilde{I}(\Lambda(U_r))
-
\delta
+
K \eps
\nonumber
\\
&
\leq
\widetilde{I}^{*}_{
\delta
}
(U_r)
-
\delta
+
K \eps
.
\end{align}
Combining \eqref{eq:as2:sk-case:log-p-n-u-upper-bound},
\eqref{as2:upper-bound:proof:max-expect-are-almost-commutative}, 
\eqref{vektor-sk-restricted-upper-bound}, we obtain 
\begin{align*}
p(\mathcal{V})
=
\lim_{
N
\uparrow
+\infty
}
\frac{1}{N}
\E
\left[
\log
P_N(\mathcal{V})
\right]
&
\leq
\max_{
r 
\in 
\{ 1,\dots,p \}
}
\lim_{
N
\uparrow
+\infty
}
\frac{1}{N}
\E
\left[
\log
P_N(\mathcal{V}(U^{(r)}))
\right]
\\
&
\leq
\max_{
r 
\in 
\{ 1,\dots,p \}
}
\widetilde{I}^{*}_{
\delta
}
(U^{(r)})
+
K \eps
-
\delta
\\
&
\leq
\sup_{
U \in \mathcal{V}
}
\widetilde{I}^{*}_{
\delta
}
(V)
+
K \eps
-
\delta
.
\end{align*}
Taking $\delta \downarrow +0$ and $\eps \downarrow +0$ limits, we get
\begin{align}
\label{vektor-sk-p-upper-bound-}
p(\mathcal{V})
\leq
\sup_{
V \in \mathcal{V}
}
\widetilde{I}^{*}
(U)
.
\end{align}
The averaging property of the RPC (see, e.g.,
\cite[Theorem~5.4]{AizenmanSimsStarr2006} or property
\eqref{the-averaging-property} of
Theorem~\ref{thm:remainder:some-limiting-grem-properties}) gives
\begin{align}
\label{eq:as2:phi-of-b}
\Phi_N(x,\mathcal{U})[B]
=
\frac{\beta^2}{2}
\sum_{k=1}^{n}
x_k
\left(
\Vert
Q^{(k+1)}
\Vert_\text{F}^2
-
\Vert
Q^{(k)}
\Vert_\text{F}^2
\right)
.
\end{align}
To finish the proof it remains to show that, for any fixed $\Lambda \in
\symmetric(d)$, we have
\begin{align*}
\widetilde{I}(\Lambda)
=
X_0(x, \mathcal{Q}, \Lambda, U)
\end{align*}
which is assured by Lemma \ref{vektor-sk-existence-of-laplace-trasform}.
\end{proof}

\subsection{Free energy lower bound} 
In this subsection, we return to the notations of Section
\ref{vektor-sk-free-energy-upper-bound-multidim-sk-model}.
\begin{lemma}
\label{vektor-sk-restricted-free-enrgy-rough-lower-bound-lemma}
For any 
$
\mathcal{B} 
\subset 
\symmetric(d)
$ 
such that  
$
\interior 
\mathcal{B} 
\cap
\interior
\mathcal{U}
\neq 
\emptyset
$ 
there exists 
$
\Delta
\subset
\Sigma
$
with
$
\interior \Delta
\neq
\emptyset
$
such that 
\begin{multline}
\lim_{
N
\uparrow
\infty
}
\frac{1}{N}
\E
\left[
\vphantom{\int}
\smash{
\int_{
\Sigma_N(\mathcal{B})
\times
\mathcal{A}
}
}
\exp
\left(
N
\langle
\Lambda
,
R_N(\sigma,\sigma)
\rangle
+
\vphantom{\sum}
\smash{
\sum_{
i=1
}^{N}
}
\langle
A_i(\alpha)
,
\sigma_i
\rangle
\right)
\dd
\pi_N(\sigma, \alpha)
\right]
\\
\geq
\log
\int_{
\Delta
}
\exp
\left(
\langle
(\beta^2 U + \Lambda)
\sigma
,
\sigma
\rangle
\right)
\dd \mu(\sigma)
>
-\infty
.
\end{multline}
\end{lemma}
\begin{proof}
In view of \eqref{eq:vektor-sk-x-k-recursive-definition}, iterative
application of the Jensen inequality with respect to $\E_{z^{(k)}}$ leads to the
following
\begin{align*}
\E
\left[
X_{n+1}(x, \mathcal{Q}, \Lambda, U)
\right]
\leq
X_0(x, \mathcal{Q}, \Lambda, U)
.
\end{align*}
Performing the Gaussian integration, we get
\begin{align*}
\E
\left[
X_{n+1}(x, \mathcal{Q}, \Lambda, U)
\right]
&
\geq
\log
\int_\Delta
\exp
\left(
\langle
(\beta^2 U + \Lambda)
\sigma
,
\sigma
\rangle
\right)
\dd \mu(\sigma)
,
\end{align*}
where $\Delta \subset \Sigma$ is such that $\mu(\Delta) > 0$ and
$
\{ 
R(\sigma,\sigma)
:
\sigma \in \Delta^N
\} 
\subset 
\mathcal{B}
$.
\end{proof}
Define the following Legendre transform
\begin{align}
\label{eq:as2:sk-upper-bound:legendre-transform:i-star-hat-definition-}
\widehat{I}^{*}
(U)
\equiv
\inf_{
\substack{
x
\in \mathcal{Q}^\prime(1,1)
,
\\
\mathcal{Q}
\in \mathcal{Q}^\prime(U,d)
,
\\
\Lambda 
\in
\symmetric(d)
}
}
\left[
-
\langle
U
,
\Lambda
\rangle
-
\Phi[B]
+
\widetilde{I}(\Lambda)
+
\mathcal{R}(x, \mathcal{Q}, U)
\right]
.
\end{align}
\begin{proof}[Proof of Theorem~\ref{thm:as2:pressure-lower-bound}]
As it is the case with the proof of Theorem \ref{thm:as2:pressure-upper-bound},
this proof also follows in essence almost literally the proof of Theorem
\ref{vektor-sk-generic-quenched-ldp-lower-bound}. The notable difference is
that we apply the Gaussian comparison in order to ``compute''
the rate function in a somewhat more explicit way.  

In notations of Theorem \ref{vektor-sk-generic-quenched-ldp-lower-bound} we
are in the following situation: $\mathcal{X} \equiv \symmetric(d)$ and
$\mathfrak{X}$ is the topology induced by any norm on $\symmetric(d)$.

Let $B(U,\eps)$ be the ball (in the Hilbert-Schmidt norm) of radius
$\eps>0$ around some arbitrary $U \in \mathcal{V}$. Let us prove at first that
\begin{align}
\label{vektor-sk-ldp-lower-bound-grem-lim-eps-0-N-infty}
\lim_{
\eps\downarrow+0
}
\lim_{
N\uparrow\infty
}
\frac{1}{N}
\E
\left[
\log
P_N(B(U,\eps))
\right]
\geq
\widehat{I}^{*}(U)
.
\end{align}
Similarly to
\eqref{vektor-sk-comaprison-applied-to-the-restricted-sk}, for any
$
(x,\mathcal{Q})
$,
we have 
\begin{align}
\label{vektor-sk-expect-log-p-n-b-eps-equals-epext-log-tilde-p-b-eps-plus-remainder}
\E
&
\left[
\frac{1}{N}
\log
P_{
N
}(
B(U,\eps)
)
\right]
\nonumber
\\
&
=
\frac{1}{N}
\E
\left[
\log
\widetilde{P}_N(B(U,\eps))
\right]
-
\Phi[B]
+
\mathcal{R}_N(x, \mathcal{Q}, U, B(U,\eps))
+
\mathcal{O}(\eps)
.
\end{align}
The random measure $\widetilde{P}_N$ satisfies the 
assumptions of Corollary
\ref{vektor-sk-generic-quenched-ldp-lower-bound-strict-convexity-everywhere}.
Indeed:
\begin{enumerate}
\item 
Due to representation \eqref{vektor-sk-i-n-grem-equals}, mapping
$I(\cdot)$ is differentiable
with respect to $\Lambda$. Henceforth assumption
\eqref{itm:as2:differentiability-of-the-rate-function} of the corollary is also
fulfilled.
\item
Let us note at first that, thanks to Proposition
\ref{vektor-sk-concentration-sk-free-energy-grem-free-energy}, we
have
$
\mathcal{D}(I)
=
\R^d
$.
Thus, the assumption \eqref{itm:as2:regularity-for-the-ldp-lower-bound} of
Corollary \ref{vektor-sk-generic-quenched-ldp-lower-bound-strict-convexity-everywhere} is
satisfied, as is condition
\eqref{vektor-sk-zero-in-interior-of-domain-of-I}.  
\end{enumerate}
Moreover, the assumptions of Lemma \ref{vektor-sk-log-laplace-trasform-concentration}
are satisfied:
\begin{enumerate}
\item 
The concentration of measure condition is satisfied due to  Proposition
\ref{vektor-sk-concentration-sk-free-energy-grem-free-energy}.
\item 
The tail decay is obvious since the family 
$
\{
\widetilde{P}_N
:
N \in \N
\}
$ 
has compact support. 
Namely, for all $N \in \N$, we have $\supp \widetilde{P}_N = \mathcal{U}$. Thus the
measure $\widetilde{Q}_{N,\Lambda}$ (cf.
\eqref{vektor-sk-tilde-q-n-lambda-measure-definition}) generated by
$\widetilde{P}_N$ has the same support.
Thus, $\supp \widetilde{Q}_{N,\Lambda} = \mathcal{U}$.
\item
The non-degeneracy is assured by Lemma
\ref{vektor-sk-restricted-free-enrgy-rough-lower-bound-lemma}.
\end{enumerate}
Hence, due to
\eqref{vektor-sk-expect-log-p-n-b-eps-equals-epext-log-tilde-p-b-eps-plus-remainder},
arguing in the same way as in Theorem \ref{vektor-sk-generic-quenched-ldp-lower-bound}, we arrive at
\eqref{vektor-sk-ldp-lower-bound-grem-lim-eps-0-N-infty}. 
Note that the 
$
N \uparrow +\infty
$ 
limit of 
$
\mathcal{R}_N(x, \mathcal{Q}, U, B(U,\eps))
$
exists, since in
\eqref{vektor-sk-expect-log-p-n-b-eps-equals-epext-log-tilde-p-b-eps-plus-remainder}
the limits of the other two $N$-dependent quantities exist due to
\cite{Guerra-Toninelli-Generalized-SK-2003}. The subsequent $\eps \downarrow +0$
limit of the remainder term exists due to the monotonicity.

Finally, taking the
supremum over $U\in\mathcal{V}$ in
\eqref{vektor-sk-ldp-lower-bound-grem-lim-eps-0-N-infty}, we get \eqref{vektor-sk-lower-bound-plus-remainder}. 
\end{proof}

\section{Guerra's comparison scheme}
\label{sec:as2:guerras-scheme}

In this section, we shall apply Guerra's comparison scheme (see the recent
accounts by \cite{GuerraReview2005,Talagrand2007a,AizenmanSimsStarr2006}) to the
SK model with multidimensional spins. However, we shall use also the ideas (and
the language) of \cite{AizenmanSimsStarr2003}. In particular, we shall use the
same local comparison functional \eqref{eq:remainder:comparison-functional} as in
the $\text{AS}^2$ scheme, see \eqref{eq:as2-phi-of-t-definition}. The section
contains the proofs of the upper
\eqref{eq:multidimensional-sk:free-energy-upper-bound} and lower
\eqref{eq:multidimensional-sk:vektor-sk-lower-bound-plus-remainder} bounds on the
free energy without Assumption~\ref{asm:multidimensional-sk:hadamard-squares}.
The proofs use the GREM-like Gaussian processes, RPCs as in the $\text{AS}^2$
scheme. We also obtain an analytic representation of the remainder term (which is
an artifact of this scheme) using the properties of the Bolthausen-Sznitman
coalescent.

\subsection{Multidimensional Guerra's scheme}
Let 
$
\xi 
= 
\xi(
x_1,\ldots,x_n
)
$
be an RPC process. Theorem~5.3 of \cite{AizenmanSimsStarr2006} guarantees
that there exists a rearrangement 
$
\xi 
=
\{
\xi(i)
\}_{
i \in \N
} 
$ 
of the 
$
\xi
$'s 
atoms in a decreasing order. Recall \eqref{eq:introduction:lexicographic-overlap} and define a (random)
\emph{limiting ultrametric overlap} 
$
q_\text{L}: \N^2 \to [0;n] \cap \Z
$ 
as follows 
\begin{align}
\label{eq:review:limiting-ultrametric-overlap}
q_\text{L}(i,j)
\equiv
1
+
\max
\{
k \in [0;n] \cap \Z
:
[\pi(i)]_k
=
[\pi(j)]_k
\}
,
\end{align}
where we use the convention that 
$
\max \emptyset = 0
$.
This overlap valuation induces a sequence of \emph{random partitions} of 
$
\N
$ into \emph{equivalence classes}. Namely, given a
$
k \in \N \cap [0;n]
$, 
we define, for any
$
i,j \in \N
$, the \emph{Bolthausen-Sznitman equivalence relation} as follows
\begin{align}
\label{eq:review:bolthausen-sznitman-equivalence}
i
\underset{k}{\sim}
j
\overset{\text{def}}{\Longleftrightarrow}
q_\text{L}(i,j) \geq k
.
\end{align}
Given $n \in \N$, assume that
$
x
$
and 
$
\mathcal{Q}
$
satisfy \eqref{eq:partition-of-the-unit-interval} 
and 
\eqref{eq:multidimesnional-sk:monotonicity-of-q-matrix-sequence},
respectively.
Recall the definitions of
the Gaussian processes $X$ and $A$ which satisfy \eqref{eq:as2:sk-with-multidimensional-spins-process} and
\eqref{eq:as2:a-process-definition}, respectively.
We consider, for $t \in [0;1]$, the following
interpolating Hamiltonian on the configuration space $\Sigma_N \times \mathcal{A}$
\begin{align}
\label{eq:remainder:interpolating-hamiltonian}
H_t(\sigma,\alpha)
\equiv
\sqrt{t} 
X(\sigma)
+
\sqrt{1-t}
A(\sigma,\alpha)
.
\end{align}
Given $\mathcal{U} \subset \symmetric^+(d)$, the Hamiltonian
\eqref{eq:remainder:interpolating-hamiltonian} in the usual way induces the following \emph{local free energy}
\begin{align}
\label{eq:as2-phi-of-t-definition}
\varphi_N(t,x,Q,\mathcal{U})
\equiv
\Phi_N(x,\mathcal{U})
\left[
H_t
\right]
,
\end{align}
where we use the same local comparison functional
\eqref{eq:remainder:comparison-functional} as in the $\text{AS}^2$ scheme.
Using  \eqref{eq:as2:local-free-energy}, we obtain then
\begin{align*}
\varphi(0,x,Q,\mathcal{U}) 
= 
\Phi_N(x,\mathcal{U})
[A]
\text{ and }
\varphi(1,x,Q,\mathcal{U}) 
= 
\Phi_N(x,\mathcal{U})
[X]
=
p_N(\mathcal{U})
.
\end{align*}
Now, we are going to disintegrate the Gibbs measure
defined on $
\mathcal{U}
\times 
\mathcal{A}
$ 
into two Gibbs measures acting on $\mathcal{U}$ and $\mathcal{A}$
separately. For this purpose we define the correspondent (random) \emph{local
free energy} on $\mathcal{U}$ as follows
\begin{align}
\label{eq:as2:psi-of-t-v-alpha}
\psi(t,x,Q,\alpha,\mathcal{U})
\equiv
\log
\int_{
\Sigma_N(\mathcal{U})
}
\exp
\left[
\beta \sqrt{N}
H_t(\sigma,\alpha)
\right]
\dd \mu^{\otimes N}(\sigma)
.
\end{align}
For $\alpha \in \mathcal{A}$, we can define the (random) \emph{local
Gibbs measure} 
$
\mathcal{G}(t,Q,\alpha,\mathcal{U})
\in
\mathcal{M}_1(\Sigma_N) 
$ 
by demanding that the following holds
\begin{align*}
\frac{
\dd 
\mathcal{G}(t,x,Q,\alpha,\mathcal{U})
}{
\dd 
\mu^{\otimes N} 
}
(\sigma)
\equiv
\I_{
\Sigma_N(\mathcal{U})
}
(\sigma)
\exp
\left[
\beta \sqrt{N}
H_t(\sigma,\alpha)
-
\psi(t,x,Q,\mathcal{U},\alpha)
\right]
.
\end{align*}
Let us define a certain reweighting of the RPC $\xi$ with the help of
\eqref{eq:as2:psi-of-t-v-alpha}.  We define the random point process 
$
\{
\tilde{\xi}
\}_{
\alpha
\in
\mathcal{A}
}
$ 
in the following way
\begin{align*}
\tilde{\xi}(\alpha)
\equiv
\xi(\alpha)
\exp
\left(
\psi(t,x,Q,\mathcal{U},\alpha)
\right)
.
\end{align*}
We also define the \emph{normalisation operation}
$
\mathcal{N}
:
\mathcal{M}_{\text{f}}(\mathcal{A})
\to
\mathcal{M}_1(\mathcal{A})
$ 
as
\begin{align*}
\mathcal{N}
\left(
\xi
\right)
(\alpha)
\equiv
\frac{
\xi(\alpha)
}{
\sum_{
\alpha'
\in
\mathcal{A}
}
\xi(\alpha')
}
.
\end{align*}
We introduce the \emph{local Gibbs measure}
$
\mathcal{G}(t,x,Q,\mathcal{U})
\in
\mathcal{M}_1
(
\mathcal{U} \times \mathcal{A}
)
$, 
for any 
$
\mathcal{V} 
\subset
\mathcal{U} \times \mathcal{A}
$, 
as follows 
\begin{align}
\label{eq:as2:desintegrated-gibbs-measure}
\mathcal{G}(t,x,Q,\mathcal{U})
\left[
\mathcal{V} 
\right]
\equiv
\sum_{
\alpha
\in
\mathcal{A}_n
}
\mathcal{N}(\tilde{\xi})(\alpha)
\mathcal{G}(t,x,Q,\alpha,\mathcal{U})
\left[
\mathcal{V}
\right]
.
\end{align}
Finally, we introduce, what shall call \emph{Guerra's remainder term}:
\begin{align}
\label{eq:lcture-03:the-remainder-term-definition}
\mathcal{R}(t,x,Q,\mathcal{U})
\equiv
-
\frac{\beta^2}{2}
\E
\left[
\mathcal{G}(t,x,Q,\mathcal{U})
\otimes
\mathcal{G}(t,x,Q,\mathcal{U})
\left[
\Vert
R(\sigma^{1},\sigma^{(2)})
-
Q(\alpha^{(1)},\alpha^{(2)})
\Vert_\text{F}^2
\right]
\right]
.
\end{align}
Note that \eqref{eq:lcture-03:the-remainder-term-definition}
coincides with \eqref{ab.200} after substituting \eqref{eq:as2:interpolating-as2-hamiltonian} with
\eqref{eq:as2:guerras-interpolation}.
\subsection{Local comparison}
We recall for completeness the following.
\begin{proposition}[Ruelle~\cite{Ruelle1987},
Bolthausen and Sznitman~\cite{BolthausenSznitman1998}]
\label{prp:review:expected-mutual-overlap-distribution}
For any 
$
k \in [1;n+1] \cap \N
$,
we have
\begin{align*}
\E
\left[
\mathcal{N}(\xi)
\otimes
\mathcal{N}(\xi)
\left\{
(\alpha^{(1)},\alpha^{(2)})
\in
\mathcal{A}^2
:
q_{
\text{L}
}
(\alpha^1,\alpha^2)
\leq
k
\right\}
\right]
=
x_{k}
.
\end{align*}
\end{proposition}
The results of Section~\ref{sec:as2:free-energy-upper-and-lower-bounds} can
be straightforwardly generalised to the comparison scheme based on
\eqref{eq:remainder:interpolating-hamiltonian}. 
Given $\eps, \delta > 0$ and $\Lambda \in \symmetric(d)$, define
\begin{align}
\label{eq:remainder:small-vicinity-of-u}
\mathcal{V}(\Lambda,\mathcal{U},\eps,\delta)
\equiv
\{
U^\prime
\in
\symmetric(d)
:
\Vert
U^\prime
-
U
\Vert_\text{F}
<
\eps
,
\langle
U^\prime
-
U
,
\Lambda
\rangle
<
\delta
\}
.
\end{align}
We now specialise to the case 
$
\mathcal{U} = \Sigma_N(\mathcal{V}(\Lambda,U,\eps,\delta))
$.
\begin{lemma}
\label{thm:as2:guerras-interpolation}
We have
\begin{align}
\label{eq:as2:guerras-interpolation-speed}
\frac{\partial}{\partial t}
\varphi_N(t,x,Q,\mathcal{V}(\Lambda,U,\eps,\delta))
=
&
\mathcal{R}(t,x,Q,\Sigma_N(\mathcal{A}(\Lambda,U,\eps,\delta)))
\nonumber
\\
&
-
\frac{\beta^2}{2}
\sum_{
k=1
}^{
n
}
x_k
\left(
\Vert
Q^{(k+1)}
\Vert_\text{F}^2
-
\Vert
Q^{(k)}
\Vert_\text{F}^2
\right)
+
\mathcal{O}(\eps)
.
\end{align}
\end{lemma}
\begin{proof}
This is an immediate consequence of
Proposition~\ref{prp:as2:gaussian-comparison-of-free-energy}.
Indeed, recalling that $Q(\alpha^{(1)},\alpha^{(1)}) = U$, and setting 
$
\mathcal{U}
\equiv
\Sigma(B(U,\eps))
$, 
we have
\begin{align}
\label{eq:remainder:d-phi-d-t}
\frac{\partial}{\partial t}
&
\varphi(t,x,Q,\mathcal{U})
\nonumber
\\
&
=
\frac{\beta^2}{2}
\E
\left[
\mathcal{G}(t,x,Q,\mathcal{U})
\otimes
\mathcal{G}(t,x,Q,\mathcal{U})
\left[
\Vert
R(\sigma^{(1)},\sigma^{(1)})
-
U
\Vert_\text{F}^2
-
\Vert
R(\sigma^{(1)},\sigma^{(2)})
-
Q(\alpha^{(1)},\alpha^{(2)})
\Vert_\text{F}^2
\right.
\right.
\nonumber
\\
&
\quad
\quad
-
\left.
\left.
\left(
\Vert
U
\Vert_\text{F}^2
-
\Vert
Q(\alpha^{(1)},\alpha^{(2)})
\Vert_\text{F}^2
\right)
\right]
\right]
\nonumber
\\
&
=
-
\frac{\beta^2}{2}
\E
\left[
\mathcal{G}(t,x,Q,\mathcal{U})
\otimes
\mathcal{G}(t,x,Q,\mathcal{U})
\left[
\Vert
R(\sigma^{(1)},\sigma^{(2)})
-
Q(\alpha^{(1)},\alpha^{(2)})
\Vert_\text{F}^2
\right]
\right]
\nonumber
\\
&
\quad
\quad
-
\frac{\beta^2}{2}
\E
\left[
\mathcal{G}(t,x,Q,\mathcal{U})
\otimes
\mathcal{G}(t,x,Q,\mathcal{U})
\left[
\Vert
U
\Vert_\text{F}^2
-
\Vert
Q(\alpha^{(1)},\alpha^{(2)})
\Vert_\text{F}^2
\right]
\right]
+
\mathcal{O}(\eps)
.
\end{align}
Using 
Proposition~\ref{prp:review:expected-mutual-overlap-distribution}, we get
\begin{align}
\label{eq:remainder:phi-of-b-term}
\frac{\beta^2}{2}
&
\E
\left[
\mathcal{G}(t,x,Q,\mathcal{U})
\otimes
\mathcal{G}(t,x,Q,\mathcal{U})
\left[
\Vert
U
\Vert_\text{F}^2
-
\Vert
Q(\alpha^{(1)},\alpha^{(2)})
\Vert_\text{F}^2
\right]
\right]
\nonumber
\\
&
=
\frac{\beta^2}{2}
\E
\left[
\mathcal{N}(\xi)
\otimes
\mathcal{N}(\xi)
\left[
\sum_{
k
=
q_{\text{L}}
(\alpha^{(1)},\alpha^{(2)})
}^n
\left(
\Vert
Q^{(k+1)}
\Vert_\text{F}^2
-
\Vert
Q^{(k)}
\Vert_\text{F}^2
\right)
\right]
\right]
\nonumber
\\
&
=
\frac{\beta^2}{2}
\sum_{k=1}^n
\left(
\Vert
Q^{(k+1)}
\Vert_\text{F}^2
-
\Vert
Q^{(k)}
\Vert_\text{F}^2
\right)
\E
\left[
\mathcal{N}(\xi)
\otimes
\mathcal{N}(\xi)
\{
k
\geq
q_{\text{L}}
(\alpha^{(1)},\alpha^{(2)})
\}
\right]
\nonumber
\\
&
=
\frac{\beta^2}{2}
\sum_{
k=1
}^{
n
}
x_k
\left(
\Vert
Q^{(k+1)}
\Vert_\text{F}^2
-
\Vert
Q^{(k)}
\Vert_\text{F}^2
\right)
.
\end{align}
Combining \eqref{eq:remainder:d-phi-d-t} and
\eqref{eq:remainder:phi-of-b-term}, we get
\eqref{eq:as2:guerras-interpolation-speed}

\end{proof}
\begin{lemma}
\label{lem:remainder:guerras-sum-rule}
We have
\begin{align}
\label{eq:as2:sum-rule}
p_N(
\Sigma_N(B(U,\eps))
)
=
&
\Phi_N(x,\Sigma_N(B(U,\eps)))
\left[
A
\right]
-
\frac{\beta^2}{2}
\sum_{
k=1
}^{
n
}
x_k
\left(
\Vert
Q^{(k+1)}
\Vert_\text{F}^2
-
\Vert
Q^{(k)}
\Vert_\text{F}^2
\right)
\nonumber
\\
&
+
\int_{0}^{1}
\mathcal{R}(t,x,Q,\Sigma_N(B(U,\eps))
\dd t
+
\mathcal{O}(\eps)
.
\end{align}
\end{lemma}
\begin{remark}
Note that the above lemma also holds if we substitute $B(U,\eps)$
with the smaller set $\mathcal{V}(\Lambda,U,\eps,\delta)$.
\end{remark}
\begin{proof}
The claim follows from \eqref{eq:as2:guerras-interpolation-speed}
by integration. 
\end{proof}

\begin{proposition}
\label{eq:remainder:guerras-upper-and-lower-bounds}
There exists $C = C(\Sigma,\mu) > 0$ such that,
for all $U \in \symmetric^+(d)$ as above, and all $\eps, \delta > 0$,
there exists an $\delta$-minimal Lagrange multiplier
$
\Lambda
=
\Lambda(U,\eps,\delta) 
\in 
\symmetric(d)
$
in \eqref{eq:as2:local-parisi-functional}
such that, for all $t \in [0;1]$, 
and all $(x,\mathcal{Q})$,
we have
\begin{align}
\label{eq:remainder:guerras-sum-rule-ldp-upper-bound}
p_N(
\Sigma_N(\mathcal{V}(\Lambda,U,\eps,\delta))
)
\leq
&
\inf_{
\Lambda
\in 
\symmetric(d)
}
f(x, \mathcal{Q}, U, \Lambda)
+
C(\eps + \delta)
\end{align}
and
\begin{align}
\label{eq:remainder:guerras-sum-rule-ldp-lower-bound}
\lim_{
N \uparrow +\infty
}
p_N(
\Sigma_N(B(U,\eps))
)
\geq
&
\inf_{
\Lambda
\in 
\symmetric(d)
}
f(x, \mathcal{Q}, U, \Lambda)
+
\lim_{
N \uparrow +\infty
}
\int_{0}^{1}
\mathcal{R}(t,x,Q,\Sigma_N(B(U,\eps)))
\dd t
\nonumber
\\
&
-
C(\eps + \delta)
.
\end{align}
\end{proposition}
\begin{remark}
The following upper bound also holds true. There exists $C = C(\Sigma,\mu)>0$,
such that, for any $\Lambda \in \symmetric(d)$,

\begin{align}
\label{eq:remainder:local-free-energy-lambda-dependent-upper-bound}
p_N(
\Sigma_N(B(U,\eps))
)
\leq
f(x, \mathcal{Q}, U, \Lambda)
+
C
\Vert
\Lambda
\Vert_{\text{F}}
\eps
.
\end{align}
\end{remark}
\begin{proof}
The result follows from Lemma~\ref{lem:remainder:guerras-sum-rule} by the same
arguments as in the proofs of Theorems~\ref{thm:as2:pressure-upper-bound} and
\ref{thm:as2:pressure-lower-bound}.
\end{proof}
\subsection{Free energy upper and lower bounds}
Similarly to the quenched LDP bounds for the 
$
\text{AS}^2
$ 
scheme in the SK model with multidimensional spins (see
Section~\ref{sec:sk-with-multidimensional-spins:gaertner-ellis}), we get the
quenched LDP bounds for Guerra's scheme in the same model without
Assumption~\ref{asm:multidimensional-sk:hadamard-squares} on $Q$.

Recall the definition of the local Parisi functional $f$ \eqref{eq:as2:local-parisi-functional}.
\begin{theorem}
\label{thm:multidimensional-sk:free-energy-upper-bound}
For any closed set $\mathcal{V} \subset \symmetric(d)$, we have 
\begin{align}
\label{eq:multidimensional-sk:free-energy-upper-bound}
p(\mathcal{V})
\leq
\sup_{
U
\in
\mathcal{V}
\cap
\mathcal{U}
}
\inf_{
(
x
,
\mathcal{Q}
,
\Lambda
)
}
f(x, \mathcal{Q}, \Lambda, U)
,
\end{align}
where the infimum runs over all $x$ satisfying
\eqref{eq:partition-of-the-unit-interval}, all $\mathcal{Q}$ satisfying 
\eqref{eq:multidimesnional-sk:monotonicity-of-q-matrix-sequence} and all $\Lambda \in
\symmetric(d)$.
\end{theorem}
\begin{proof}
The proof is identical to the one of Theorem~\ref{thm:as2:pressure-upper-bound}.
\end{proof}
Define the \emph{local limiting Guerra remainder term}  
$
\mathcal{R}(x, \mathcal{Q}, U)
$
as follows
\begin{align}
\label{eq:remainder:local-limiting-guerra-remainder}
\mathcal{R}(x, \mathcal{Q}, U)
\equiv
-
\lim_{
\eps
\downarrow
+0
}
\lim_{
N
\uparrow
+\infty
}
\int_{0}^{1}
\mathcal{R}(t,\Sigma_N(B(U,\eps)))
\dd t
\leq
0
.
\end{align}
The existence of the limits in
\eqref{eq:remainder:local-limiting-guerra-remainder} is proved similar to the
case of the $\text{AS}^2$ scheme, see the proof of
Theorem~\ref{thm:as2:pressure-lower-bound}.
\begin{theorem}
\label{thm:multidimensional-sk:free-energy-lower-bound}
For any open set $\mathcal{V} \subset \symmetric(d)$, we have
\begin{align}
\label{eq:multidimensional-sk:vektor-sk-lower-bound-plus-remainder}
p(
\mathcal{V}
)
\geq
\sup_{
U
\in
\mathcal{V}
\cap
\mathcal{U}
}
\inf_{
(
x
,
\mathcal{Q}
,
\Lambda
)
}
\left[
f(x, \mathcal{Q}, \Lambda, U)
+
\mathcal{R}(x, \mathcal{Q}, U)
\right]
,
\end{align}
where the infimum runs over all $x$ satisfying
\eqref{eq:partition-of-the-unit-interval}; all $\mathcal{Q}$ satisfying
\eqref{eq:multidimesnional-sk:monotonicity-of-q-matrix-sequence} and
all 
$
\Lambda 
\in
\symmetric(d)
$.
\end{theorem}
\begin{proof}
The proof is identical to the one of
Theorem~\ref{thm:as2:pressure-lower-bound}. The only new ingredient is
Lemma~\ref{thm:as2:guerras-interpolation} needed to recover Guerra's
remainder term \eqref{eq:lcture-03:the-remainder-term-definition}.
\end{proof}

\subsection{The filtered $d$-dimensional GREM}
\label{sec:remainder:the-filtered-d-dimensional-grem}
Given $U \in \symmetric^+(d)$ non-negative definite, 
denote by $\mathcal{Q}(U, d)$ the set of all càdlàg (right continuous with
left limits) $\symmetric^+(d)$-valued non-decreasing paths which end in matrix
$U$, i.e.,
\begin{align}
\label{eq:remainder:set-of-multidim-order-parameters}
\mathcal{Q}(U, d)
\equiv 
\{
\rho : [0;1] \to \symmetric^+(d)
\mid
\rho(0)=0
;
\rho(1)=U
;
\rho(t) 
\preceq
\rho(s)
,
\text{ for }
t \leq s 
;
\rho
\text{ is c\'adl\'ag}
\}
.
\end{align}
Define the natural
inverse $
\rho^{-1}: \Im \rho \to [0;1]
$
as 
\begin{align*}
\rho^{-1}(Q) 
\equiv
\inf\{
t \in [0;1] 
\mid
\rho(t) \succeq Q
\}
,
\end{align*}
where $\Im \rho \equiv \rho([0;1])$. 
Let $x \equiv \rho^{-1} \circ \rho \in \mathcal{Q}(1,1)$.

Let also $\mathcal{Q}^{\prime}(U,d) \subset \mathcal{Q}(U,d)$ be the space of
all piece-wise constant paths in $\mathcal{Q}(U,d)$ with finite (but arbitrary) 
number of jumps with an additional requirement that they have a jump at $x=1$. 
Given some $\rho \in \mathcal{Q}^\prime(U,d)$, we enumerate its jumps and define
the finite collection of matrices $\{Q^{(k)}\}_{k=0}^{n+1} \equiv \Im \rho
\subset \R^d$. This
implies that there exist 
$
\{
x_k
\}_{
k=0
}^{n+1}
\subset
\R
$
such that
\begin{align*}
0 \equiv x_0 < x_1 < \ldots < x_{n} < x_{n+1} \equiv 1
,
\\
0 \equiv Q^{(0)}
\preceq
Q^{(1)}
\preceq
Q^{(2)}
\preceq
\dots
\preceq
Q^{(n+1)}
\equiv
U
,
\end{align*}
where 
$
\rho(x_k) = Q^{(k)}
$.
Let us associate to $\rho \in \mathcal{Q}^\prime(U,d)$ a new path $\tilde{\rho}
\in \mathcal{Q}(U,d)$ which is obtained by the linear interpolation of the path
$\rho$. Namely, let
\begin{align*}
\tilde{\rho}(t)
\equiv
Q^{(k)}
+
(Q^{(k+1)}-Q^{(k)})
\frac{
t-x_k
}{
x_{k+1}-x_k
}
,
t \in [x_k; x_{k+1})
.
\end{align*}
Let $g: \R^d \to \R$ be a function satisfying
Assumption~\ref{asm:remainder:end-condition}. Let us introduce the
\emph{filtered $d$-dimensional GREM process} $W$. Let
\begin{align*}
W
\equiv
\left\{
\{
W_k(t,[\alpha]_k)
\}_{
t \in \R_+
}
:
\alpha \in \mathcal{A}
,
k \in [0;n] \cap \N
\right\}
\end{align*}
be the collection of independent (for different $\alpha$ and $k$) $\R^d$-valued
correlated Brownian motions satisfying
\begin{align*}
W_k(t,[\alpha]_k)
\sim
(Q^{(k+1)}-Q^{(k)})^{1/2}
W
\left(
\frac{
t-x_{k}
}{
x_{k+1}-x_{k}
}
\right)
,
\end{align*}
where $\{W(t)\}_{t \in \R_+}$ is the standard (uncorrelated) $\R^d$-valued
Brownian motion.
Now, for $k \in [0;n] \cap \N$, we define the $\R^d$-valued process 
$
\{
Y(t,\alpha)
\mid
\alpha \in \mathcal{A}
,
t \in [0;1]
\}
$
by
\begin{align*}
Y(t,\alpha)
\equiv
\sum_{
k=0
}^n
\I_{
[x_k;1]
}(t)
W_k(t \wedge x_{k+1},[\alpha]_k)
.
\end{align*}
\begin{lemma}
For $\alpha^{(1)},\alpha^{(2)} \in \mathcal{A}$, we
have 
\begin{align*}
\cov
\left[
Y(t_1,\alpha^{(1)})
,
Y(t_2,\alpha^{(2)})
\right]
=
\tilde{\rho}
\left(
t_1
\wedge 
t_2
\wedge 
x_{q_\text{L}(\alpha^{(1)},\alpha^{(2)})}
\right)
.
\end{align*}
\end{lemma}
\begin{proof}
The proof is straightforward.
\end{proof}
\begin{assumption}
\label{asm:remainder:end-condition}
Suppose that the function
$
g: \R^d \to \R
$
satisfies
$
g \in C^{(2)}(\R^d)
$
and, for any $c>0$,
we have
$
\int_{
\R^d
}
\exp
\left(
g(y)-c\Vert y \Vert_2^2
\right)
\dd y
<
\infty
$
and also
\begin{align}
\label{eq:remainder:end-condition}
\sup_{
y \in \R^d
}
\left(
\Vert
\nabla
g(y)
\Vert_2
+
\Vert
\nabla^2
g(y)
\Vert_2
\right)
<
+\infty
,
\end{align}
where $\nabla^2 g(y)$ denotes the matrix of second derivatives of the function
$g$ at $y \in \R^d$.
\end{assumption}
Assume $g$ satisfies the above assumption. Let 
$
f
\equiv
f_\rho
:
[0;1]
\times
\R^d
\to
\R
$
be the function satisfying the following (backward) recursive definition
\begin{align}
\label{vektor-sk-f-q-y-backward-recursive-def}
f(t, y)
\equiv
\begin{cases}
g(y)
,
&
t=1
,
\\
\frac{1}{x_k}
\log
\E
\left[
\exp
\left\{
x_k
f(
x_{k+1}
,
y
+
Y(x_{k+1},\alpha)
-
Y(t,\alpha)
)
\right\}
\right]
,
&
t
\in
[x_{k};x_{k+1})
,
\end{cases}
\end{align}
where 
$
k
\in
[0;n] \cap \N
$, $\alpha \in \mathcal{A}$ is arbitrary and fixed.
\begin{remark}
It is easy to recognise that the definition of $f$ is a continuous
``algorithmisation" of \eqref{eq:multidimensional-sk:x-0}. Namely,
$
X_k(x,\mathcal{Q}, \Lambda, U)
=
f(x_k, 0)
$,
where 
\begin{align}
\label{eq:remainder:multidim-sk-final-value}
f(1,y)
=
g(y) 
\equiv
\log
\int
_\Sigma
\exp
\left(
\sqrt{2}
\beta
\left\langle
y
,
\sigma
\right\rangle
+
\langle
\Lambda
\sigma
,
\sigma
\rangle
\right)
\dd \mu(\sigma)
.
\end{align}
\end{remark}

\subsection{A computation of the remainder term}
\label{sec:remainder-estimates:computing-the-remainder-term}

Recall the equivalence relation
\eqref{eq:review:bolthausen-sznitman-equivalence}. In words, the equivalence 
$
i
\underset{k}{\sim}
j
$
means that the atoms of the RPC $\xi$ with ranks $i$
and $j$ have the same ancestors up to the $k$-th generation. Varying the $k$ in
\eqref{eq:review:bolthausen-sznitman-equivalence}, we get
a family of equivalences on $\N$ which possesses important Markovian
properties, see \cite{BolthausenSznitman1998}.
\begin{lemma}
\label{lemma:vektor-sk-squares-of-rpc}
For all $k \in [0;n-1] \cap \N$,
we have
\begin{align}
\label{formula-rpc-i-neq-j}
\E
\left[
\sum_{
\substack{
i
\underset{k}{\sim}
j
\\
i
\underset{k+1}{\nsim}
j
}
}
\mathcal{N}
(\xi)(i)
\mathcal{N}
(\xi)(j)
\right]
=
x_{k+1}
-
x_{k}
,
\end{align}
and also
\begin{align}
\label{formula-rpc-i-eq-j}
\E
\Bigl[
\sum_{i}
\mathcal{N}
(\xi)(i)^2
\Bigr]
=
1-x_{n}
.
\end{align}
\end{lemma}
\begin{proof}
\begin{enumerate}
\item 
To prove \eqref{formula-rpc-i-neq-j} we notice that
\begin{align*}
\E
\left[
\sum_{
\substack{
i
\underset{k}{\sim}
j
\\
i
\underset{k+1}{\nsim}
j
}
}
\mathcal{N}
(\xi)(i)
\mathcal{N}
(\xi)(j)
\right]
&
=
\E
\left[
\vphantom{\sum}
\sum_{
i
\underset{k+1}{\nsim}
j
}
\mathcal{N}
(\xi)(i)
\mathcal{N}
(\xi)(j)
-
\sum_{
i
\underset{k}{\nsim}
j
}
\mathcal{N}
(\xi)(i)
\mathcal{N}
(\xi)(j)
\right]
\\
&
=
x_{k+1}
-
x_k
,
\end{align*}
where the last equality is due to
Proposition~\ref{prp:review:expected-mutual-overlap-distribution}.
\item
Similarly, \eqref{formula-rpc-i-eq-j} follows from the following observation
\begin{align*}
\E
\left[
\sum_{i}
\mathcal{N}^2
(\xi)(i)
\right]
&
=
\E
\left[
\vphantom{\sum}
\sum_{i,j}
\mathcal{N}
(\xi)(i)
\mathcal{N}
(\xi)(j)
-
\sum_{
i
\underset{n}{\nsim}
j
}
\mathcal{N}
(\xi)(i)
\mathcal{N}
(\xi)(j)
\right]
\\
&
=
1-x_{n}
,
\end{align*}
where the last equality is due to
Proposition~\ref{prp:review:expected-mutual-overlap-distribution}.
\end{enumerate}

\end{proof}
Note that, using the above notations, we readily have
\begin{align*}
A(\sigma,\alpha)
\sim
\left(
\frac{2}{N}
\right)^{1/2}
\sum_{i=1}^N
\langle
Y^{(i)}(1,\alpha)
,
\sigma_i
\rangle
,
\end{align*}
where 
$
\{
Y^{(i)}
\equiv
\{
Y^{(i)}(1,\alpha)
\}_{
\alpha \in \mathcal{A}
}
\}_{i=1}^N
$
are i.i.d. copies of 
$
\{Y(1,\alpha)\}_{
\alpha \in \mathcal{A}
}
$.
Consider the following weights
\begin{align*}
\tilde{\xi}^{(t)}
(\alpha)
\equiv
\xi(\alpha)
\exp
\left(
f(t,Y(t,\alpha))
\right)
.
\end{align*}
As in \cite{BolthausenSznitman1998}, the above weights induce the permutation
$
\tilde{\pi}^{(t)}: \N \to \mathcal{A}
$ 
such that, for all $i \in \N$, the following holds
\begin{align}
\label{eq:remainder:reordering}
\tilde{\xi}^{(t)}(\tilde{\pi}^{(t)}(i))
>
\tilde{\xi}^{(t)}(\tilde{\pi}^{(t)}(i+1))
.
\end{align}
In what follows, we shall use the short-hand notations 
$
\tilde{\xi}^{(t)}(i) 
\equiv 
\tilde{\xi}^{(t)}(\tilde{\pi}^{(t)}(i))
$,
$
\tilde{Y}^{(t)}(s,i)
\equiv
Y(s,\tilde{\pi}^{(t)}(i))
$
and
$
\tilde{Q}^{(t)}
\equiv
\{
\tilde{Q}^{(t)}(i,j)
\equiv
Q(
\tilde{\pi}^{(t)}(i)
,
\tilde{\pi}^{(t)}(j)
)
\}_{i,j \in \N}
$.
\begin{theorem}
\label{thm:remainder:some-limiting-grem-properties}
Given a discrete order parameter 
$
x
\in 
\mathcal{Q}^\prime(1,1)
$, we have
\begin{enumerate}
\item 
\emph{Independence \#1.}
The normalised RPC point process 
$
\mathcal{N}(\xi)
$
is independent from the corresponding
randomised limiting GREM overlaps
$q$.
\item
\emph{Independence \#2.}
The reordered filtered limiting GREM 
$
\tilde{Y}
$
is independent from the corresponding reordered weights
$
\tilde{\xi}
$.
\item
\emph{The reordering change of measure.}
Given 
$
I \Subset \N
$, 
let 
$
\nu_I(\cdot \vert Q)
$
be the joint distribution of
$
\{
Y(1,i)
\}_{
i \in I
}
$,
and
$
\tilde{\nu}_I(\cdot \vert Q)
$ 
be the joint distribution of 
$
\{
\tilde{Y}^{(1)}(1,i)
\}_{
i \in I
}
$ 
both conditional on 
$
Q
$.
Then
\begin{align}
\label{eq:remainder:reordering-change-of-measure}
\frac{
\dd 
\tilde{\nu}_I(\cdot \vert Q)
}{
\dd
\nu_I(\cdot \vert Q)
}
=
\prod_{
k=0
}^n
\prod_{
i
\in
\left(
I 
/
\underset{k}{\sim} 
\right)
}
\exp
\left(
x_k
\left\{
f(
x_{k+1}
,
Y(x_{k+1},i)
)
-
f_k(
x_k
Y(x_k,i)
)
\right\}
\right)
,
\end{align}
where the innermost product in the previous formula is taken over all
equivalence classes on the index set 
$
I
$ 
induced by the equivalence 
$
\underset{k}{\sim}
$.
\item
\label{the-averaging-property}
\emph{The averaging property.}
For all $s,t \in [0;1]$, we have 
\begin{align}
\label{eq:remainder:averaging-property}
\left(
\left\{
\xi^{(t)}(\alpha)
\right\}_{
\alpha \in \mathcal{A}
}
,
\tilde{Q}^{(t)}
\right)
\sim
\left(
\left\{
\xi^{(s)}(\alpha)
\right\}_{
\alpha \in \mathcal{A}
}
,
\tilde{Q}^{(s)}
\right)
.
\end{align}
\end{enumerate}
\end{theorem}
\begin{proof}
The proof is the same as in the case of the one-dimensional SK model, see
\cite{BolthausenSznitman1998,Arguin2006}. 
\end{proof}

Keeping in mind \eqref{eq:remainder:reordering-change-of-measure}, we define,
for $k \in [0;n-1] \cap \N$, the following random variables
\begin{align*}
T_k(\alpha)
\equiv
\exp
\left(
x_k
\left[
f(x_{k+1},Y(x_{k+1},\alpha))
-
f(x_k, Y(x_k, \alpha))
\right]
\right)
.
\end{align*}
Given
$
k \in [1;n] \cap \N
$,
assume that
$
\alpha^{(1)},\alpha^{(2)}
\in
\mathcal{A}
$
satisfy 
$
q_{\text{L}}
(\alpha^{(1)},\alpha^{(2)})
=
k
$.
We introduce, for notational convenience, the (random) measure 
$
\mu_k(t,\mathcal{U})
$ 
-- an
element of $
\mathcal{M}_1(\Sigma_N)
$ 
--
by demanding the following
\begin{multline}
\label{eq:lecture-04-mu-k-definition}
\mu_k
(t,\mathcal{U})
\left[
g
\right]
\equiv
\E
\left[
T_1(\alpha^{1})
\cdots
T_{k}(\alpha^{1})
T_{k+1}(\alpha^{1})
T_{k+1}(\alpha^{2})
\cdots
T_{n}(\alpha^{1})
T_{n}(\alpha^{2})
\right.
\\
\left.
\quad
\mathcal{G}(t,\alpha^{(1)},\mathcal{U})
\otimes
\mathcal{G}(t,\alpha^{(2)},\mathcal{U})
\left[
g
\right]
\right]
,
\end{multline}
where 
$
g:\mathcal{U}^2 \to \R
$
is an arbitrary measurable function such that
\eqref{eq:lecture-04-mu-k-definition} is finite. Using this notation, we can state the following lemma.
\begin{lemma}
\label{lem:remainder:reordering-undo-change-of-measure}
For any $i,j \in \N$, satisfying $i \underset{k}{\sim} j$, $i
\underset{k+1}{\nsim} j$, we have 
\begin{align}
\label{eq:remainder:mu-k-terms-oredered-unordered}
\E
\left[
\mathcal{G}(t,i,\mathcal{U})
\otimes
\mathcal{G}(t,j,\mathcal{U})
\left[
\Vert
R(\sigma^1,\sigma^2)
-
Q(i,j)
\Vert_\text{F}^2
\right]
\right]
=
\mu_k
(t,\mathcal{U})
\left[
\Vert
R(\sigma^1,\sigma^2)
-
Q^{(k)}
\Vert_\text{F}^2
\right]
.
\end{align}
\end{lemma}
\begin{proof}
This is a direct consequence of
\eqref{eq:remainder:reordering-change-of-measure} and the fact that under the
assumptions of the theorem $Q(i,j) = Q^{(k)}$.

\end{proof}
\begin{remark}
\label{rem:as2-mu-k-is-a-probability}
It is obvious from the previous theorem that $\mu_k$ is a probability measure.
\end{remark}

The main result of this subsection is an ``analytic
projection'' of the probabilistic RPC representation which integrates out the dependence on the RPC. Comparing to
\eqref{ab.200}, it has a more analytic flavor which will be exploited in the
remainder estimates (Section~\ref{sec:remainder:talagrands-remainder-estimates}). This is also a
drawback in some sense, since the initial beauty of the RPCs is lost.
\begin{theorem}
\label{thm:lectire-03-an-analytic-projection-the-v-term-without-rpc}
In the case of Guerra's interpolation \eqref{eq:as2:guerras-interpolation}, we
have
\begin{align}
\label{eq:vektor-sk-the-convenient-form-of-the-remainder}
\mathcal{R}(t,x,Q,\Sigma_N(B(U, \eps)))
=
&
\frac{1}{2}
\sum_{k=0}^{n-1}
(
x_{k+1}
-
x_{k}
)
\mu_k(t,\Sigma_N(B(U, \eps)))
\left[
\Vert
R(\sigma^1,\sigma^2)
-
Q^{(k)}
\Vert_\text{F}^2
\right]
\nonumber
\\
&
+
\mathcal{O}(\eps)
+
\mathcal{O}(1-x_{n})
,
\end{align}
as $\eps \to 0$ and $x_{n} \to 1$.
\end{theorem}

\begin{proof}
Recalling \eqref{eq:lcture-03:the-remainder-term-definition} and \eqref{eq:as2:desintegrated-gibbs-measure}, we
write 
\begin{align*}
\mathcal{R}(t,x,Q,\Sigma(U, \eps))
=
\frac{\beta^2}{2}
\E
\Bigl[
\sum_{
i, j
}
&
\mathcal{N}(\tilde{\xi})(i)
\mathcal{N}(\tilde{\xi})(j)
\\
&
\times
\mathcal{G}(t,x,Q,i,\mathcal{U})
\otimes
\mathcal{G}(t,x,Q,j,\mathcal{U})
\Bigl[
\Vert
R(\sigma^1,\sigma^2)
-
Q(i,j)
\Vert_\text{F}^2
\Bigr]
\Bigr]
.
\end{align*}
Using Theorem~\ref{thm:remainder:some-limiting-grem-properties}, we
arrive to
\begin{align*}
\mathcal{R}(t,x,Q,\Sigma(U, \eps))
=
\frac{\beta^2}{2}
\sum_{
i, j
}
&
\E
\left[
\mathcal{N}(\tilde{\xi})(i)
\mathcal{N}(\tilde{\xi})(j)
\right]
\\
&
\times
\E
\left[
\mathcal{G}(t,x,Q,i,\mathcal{U})
\otimes
\mathcal{G}(t,x,Q,j,\mathcal{U})
\left[
\Vert
R(\sigma^1,\sigma^2)
-
Q(i,j)
\Vert_\text{F}^2
\right]
\right]
.
\end{align*}
(We can interchange the summation and expectation since all summands are non-negative.)
The averaging property
(see Theorem~\ref{thm:remainder:some-limiting-grem-properties}) then gives
\begin{align}
\label{eq:remainder:r-t-u-eps-1}
\mathcal{R}(t,\Sigma(U, \eps))
=
\frac{\beta^2}{2}
\sum_{
i, j
}
\E
\left[
\mathcal{N}(\xi)(i)
\mathcal{N}(\xi)(j)
\right]
\E
\left[
\mathcal{G}(t,i,\mathcal{U})
\otimes
\mathcal{G}(t,j,\mathcal{U})
\left[
\Vert
R(\sigma^1,\sigma^2)
-
Q(i,j)
\Vert_\text{F}^2
\right]
\right]
.
\end{align}
For each $k \in [1;n-1] \cap \N$, we fix any indexes $i_0, i_0^{(k)}, j_0^{(k)}
\in \N$ such that $i_0^{(k)} \underset{k}{\sim} j_0^{(k)}$ and $ i_0^{(k)}
\underset{k+1}{\nsim} j_0^{(k)}$. Rearranging the terms in
\eqref{eq:remainder:r-t-u-eps-1}, we get
\begin{align}
\label{eq:remainder:r-t-u-eps}
\mathcal{R}(t,\Sigma(U, \eps))
&
=
\frac{\beta^2}{2}
\sum_{
k=1
}^{
n
}
\E
\left[
\mathcal{G}(t,i_0^{(k)},\mathcal{U})
\otimes
\mathcal{G}(t,j_0^{(k)},\mathcal{U})
\left[
\Vert
R(\sigma^1,\sigma^2)
-
Q^{(k)}
\Vert_\text{F}^2
\right]
\right]
\nonumber
\\
&
\quad\quad
\times
\sum_{
\substack{
i
\underset{k}{\sim}
j
\\
i
\underset{k+1}{\nsim}
j
}
}
\E
\left[
\mathcal{N}(\xi)(i)
\mathcal{N}(\xi)(j)
\right]
\nonumber
\\
&
\quad
+
\frac{\beta^2}{2}
\E
\left[
\mathcal{G}(t,i_0,\mathcal{U})
\otimes
\mathcal{G}(t,i_0,\mathcal{U})
\left[
\Vert
R(\sigma^1,\sigma^2)
-
U
\Vert_\text{F}^2
\right]
\right]
\sum_{
i
}
\E
\left[
\mathcal{N}(\xi)(i)^2
\right]
.
\end{align}
Finally, applying Lemmata \ref{lemma:vektor-sk-squares-of-rpc} and
\ref{lem:remainder:reordering-undo-change-of-measure} to
\eqref{eq:remainder:r-t-u-eps}, we arrive at
\eqref{eq:vektor-sk-the-convenient-form-of-the-remainder}.

\end{proof}
\section{The Parisi functional in terms of differential equations} 
\label{sec:remainder:the-parisi-functional-in-terms-of-differential-equations}
In this section, we study the properties of the multidimensional Parisi
functional. We derive the multidimensional version of the Parisi PDE. This
allows to represent the Parisi functional as a solution of a PDE evaluated at
the origin. We also obtain a variational representation of the Parisi
functional in terms of a HJB equation for a linear problem of diffusion
control. As a by-product, we arrive at the strict convexity of the Parisi
functional in 1-D which settles a problem of uniqueness of the optimal
Parisi order parameter posed by
\cite{TalagrandParisiMeasures2006a,PanchenkoParisiMeasures2004}.
\begin{lemma}
\label{lem:remainder:hopf-cole-validity}
Consider the function 
$
B: \R^d \times \R_+ \to \R
$ 
defined as
\begin{align*}
B(y,t)
\equiv
\frac{1}{x}
\log
\E
\left[
\exp
\left\{
x
f(y+z(t))
\right\}
\right]
,
\end{align*}
where $f: \R^d \to \R$ satisfies Assumption~\ref{asm:remainder:end-condition}
and $\{z(t)\}_{t \in [0;1]}$ is a Gaussian $\R^d$-valued process with 
$
\cov
\left[
z(t)
\right]
\equiv
Q(t)
\in 
\symmetric(d)
$ 
such that $Q(t)_{u,v}$ is differentiable, for all $u,v$.
Then
\begin{align}
\label{eq:remainder:b-differential-equation}
\partial_t
B(y,t)
=
\frac{1}{2}
\sum_{u,v=1}^d
\dot{Q}_{u,v}(t)
\left(
\partial^2_{
y_u y_v
}
B(y,t)
+
x
\partial_{
y_u
}
B(y,t)
\partial_{
y_v
}
B(y,t)
\right)
,
\quad
(t,y)
\in
(0;1)
\times
\R^d
.
\end{align}
In particular, the function $B$ is differentiable with respect to the
$t$-variable on $(0;1)$ and $C^2(\R^d)$ with respect to the $y$-variable.
\end{lemma}
\begin{proof}
Denote 
$
Z 
\equiv 
\E
\left[
\ee^{
x
f(y+z(t))
}
\right]
$.
By 
\cite[Lemma~A.1]{AizenmanSimsStarr2006}, we have
\begin{align*}
\partial_t B(y,t)
=
\frac{1}{
2 x
}
\left(
\frac{1}{Z}
\E
\Bigl[
\sum_{u,v=1}^d
\dot{Q}_{u,v}(t)
\partial_{z_u z_v}^2
\ee^{
x f(z)
}
\vert_{
z = y+z(t)
}
\Bigr]
\right)
.
\end{align*}
A straightforward calculation then gives
\begin{align}
\label{eq:remainder:b-t-derivative}
\partial_t B(y,t)
=
\frac{1}{
2 x
}
\left(
\frac{1}{Z}
\E
\Bigl[
\sum_{u,v=1}^d
\dot{Q}_{u,v}(t)
\left(
x^2 
\partial_{z_u} f(z)
\partial_{z_v} f(z)
+
x
\partial^2_{z_u z_v} f(z)
\right)
\ee^{x f(z)}
\vert_{
z = y + z(t)
}
\Bigr]
\right)
.
\end{align}
We also have
\begin{align}
\label{eq:remainder:b-y-derivative}
\partial_{y_u}
B(y,t)
=
\frac{1}{
x Z
}
\E
\left[
x
\ee^{
x f(z)
}
\partial_{z_u}
f(z)
\vert_{
z = y + z(t)
}
\right]
,
\end{align}
and
\begin{align}
\label{eq:remainder:b-second-y-derivative}
\partial^2_{y_u y_v}
B(y,t)
=
\frac{1}{x}
\Bigl(
&
\frac{1}{Z}
\E
\left[
\ee^{
x f(z)
}
\left(
x^2
\partial_{z_u}
f(z)
\partial_{z_v}
f(z)
+
\partial^2_{z_u z_v}
f(z)
\right)
\vert_{
z = y + z(t)
}
\right]
\nonumber
\\
&
-
\frac{1}{Z^2}
\E
\left[
x
\ee^{
x f(z)
}
\partial_{z_u}
f(z)
\vert_{
z = y + z(t)
}
\right]
\E
\left[
x
\ee^{
x f(z)
}
\partial_{z_v}
f(z)
\vert_{
z = y + z(t)
}
\right]
\Bigr)
.
\end{align}
Combining \eqref{eq:remainder:b-t-derivative},
\eqref{eq:remainder:b-y-derivative} and
\eqref{eq:remainder:b-second-y-derivative}, we get \eqref{eq:remainder:b-differential-equation}.
\end{proof}
\begin{proposition}
\label{proposition:vektor-sk-parisis-pde-descrete-case}
Denote $D \equiv \bigcup_{k=0}^n (x_k;x_{k+1})$.
The function $f = f_\rho$ defined in
\eqref{vektor-sk-f-q-y-backward-recursive-def} satisfies the final-value problem for the controlled
semi-linear parabolic Parisi-type PDE
\begin{align}
\label{vektor-sk-parisis-pde}
\begin{cases}
\partial_t
f(y,t)
+
\frac{1}{2}
\sum_{u,v=1}^d
\frac{\dd}{\dd t}
\tilde{\rho}_{u,v}(t)
\left(
\partial^2_{
y_u y_v
}
f(y,t)
+
x(t)
\partial_{
y_u
}
f(y,t)
\partial_{
y_v
}
f(y,t)
\right)
=
0
,
&
(t,y)
\in 
D \times \R^d
,
\\
f(1,y)=g(y)
,
&
y \in \R^d
,
\\
f(y,x_k-0)
=
f(y,x_k+0)
,
\quad
k \in [1;n] \cap \N
,
&
y \in \R^d
.
\end{cases}
\end{align}
Note that
$
\frac{\dd}{\dd t}
\tilde{\rho}(t)
=
\frac{
Q^{(k+1)}-Q^{(k)}
}{
x_{k+1}-x_k
}
$,
for $t \in (x_k; x_{k+1})$.
\end{proposition}
\begin{proof}
A successive application of Lemma~\ref{lem:remainder:hopf-cole-validity} to
\eqref{vektor-sk-f-q-y-backward-recursive-def} on the intervals $D$
starting from $(x_n;1)$ gives \eqref{vektor-sk-parisis-pde}.
\end{proof}
\begin{remark}
Note that a straightforward inspection of
\eqref{vektor-sk-f-q-y-backward-recursive-def},  using
\eqref{eq:remainder:b-t-derivative}, \eqref{eq:remainder:b-y-derivative} and
\eqref{eq:remainder:b-second-y-derivative}, shows that the function $f$ defined
in \eqref{vektor-sk-f-q-y-backward-recursive-def} is $C^1(D) \cap C([0;1])$ with respect to the $t$-variable and $C^2(\R^d)$ with respect to the $y$-variable.
\end{remark}
\begin{lemma}
Given 
$
\rho \in \mathcal{Q}^\prime(U,d)
$,
the function \eqref{vektor-sk-f-q-y-backward-recursive-def} 
satisfies the following:
\begin{align}
\label{eq:remainder:parisi-rpc-representation}
f_\rho(0,0)
=
\E
\Bigl[
\log
\sum_{
\alpha \in \mathcal{A}
}
\xi(\alpha)
\exp
\left\{
g(Y(1,\alpha))
\right\}
\Bigr]
.
\end{align}
\end{lemma}
\begin{proof}
This is an immediate consequence of the RPC averaging property 
\eqref{eq:remainder:averaging-property}.
\end{proof}
\begin{lemma}
\label{lem:remainder:parisi-formula-differentiation}
\mbox{}
\begin{enumerate}
\item 
Given $k \in [1;n] \cap \N$ and a non-negative definite matrix $Q \in
\symmetric(d)$, 
we have
\begin{align}
\label{eq:remainder:parisi-formula-differentiation}
\partial_{Q^{(k)} \leadsto Q}
f_\rho
(0,0)
=
-
\frac{1}{2}
(x_k - x_{k-1})
\E
\left[
\langle
Q
,
M
\rangle
\right]
,
\end{align}
where $M \in \R^{d \times d}$ is defined as
\begin{align*}
M_{u,v}
\equiv
&
T_1(\alpha^{(1)})
\cdots
T_{k}(\alpha^{(1)})
T_{k+1}(\alpha^{(1)})
T_{k+1}(\alpha^{(2)})
\cdots
T_{n}(\alpha^{(1)})
T_{n}(\alpha^{(2)})
\\
&
\partial_{
z_{u}
}
g(z)
\vert_{
z=Y(1,\alpha^{(1)})
}
\partial_{
z_{v}
}
g(z)
\vert_{
z=Y(1,\alpha^{(2)})
}
\end{align*}
with 
$
q_\text{L}
(\alpha^{(1)},\alpha^{(2)})
=
k
$. 
Moreover, \eqref{eq:remainder:parisi-formula-differentiation} does not depend
on the choice of 
$
\alpha^{(1)},\alpha^{(2)}
\in \mathcal{A}
$
but only on $k$.
\item 
Given a non-negative definite matrix 
$
Q \in
\symmetric(d)
$, 
we have
\begin{align}
\label{eq:remainder:parisi-formula-differentiation-wrt-variance}
\partial_{U \leadsto Q}
f_\rho
(0,0)
=
\frac{1}{2}
\E
\left[
\langle
Q
,
M^\prime
\rangle
\right]
,
\end{align}
where 
$
M^\prime
\in
\symmetric(d)
$ 
is satisfies
\begin{align*}
M^\prime_{u,v}
=
T_1(\alpha)
\cdots
T_n(\alpha)
\Bigl(
\partial^2_{
z_u
z_v
}
g(z)
+
\partial_{
z_u
}
g(z)
\partial_{
z_v
}
g(z)
\Bigr)
\Big\vert_{
z=Y(1,\alpha)
}
+
\mathcal{O}(1-x_n)
,
\end{align*}
as $x_n \to 1$. Note that
\eqref{eq:remainder:parisi-formula-differentiation-wrt-variance} obviously does
not depend on the choice of $\alpha \in \mathcal{A}$.
\end{enumerate}
\end{lemma}
\begin{proof}
Applying \cite[Lemma~A.1]{AizenmanSimsStarr2006} to
\eqref{eq:remainder:parisi-rpc-representation}, we obtain
\begin{align*}
\partial_s
&
\E
\Bigl[
\log
\sum_{
\alpha \in \mathcal{A}
}
\exp
\left\{
g(Y(1,\alpha))
\right\}
\Big\vert_{
Q^{(k)}
=
Q^{(k)} + s Q
}
\Bigr]
\\
&
=
\frac{1}{2}
\E
\Bigl[
\sum_{
u,v=1
}^d
\mathcal{N}(\tilde{\xi})
\otimes
\mathcal{N}(\tilde{\xi})
\Bigl[
\partial_s
\left(
Q(\alpha^{(1)},\alpha^{(2)})_{u,v}
\vert_{
Q^{(k)}
=
Q^{(k)} + s Q
}
\right)
\\
&
\quad\quad
\Bigl\{
\I_{
\alpha^{(1)} 
=
\alpha^{(2)} 
}
(\alpha^{(1)},\alpha^{(2)})
\Bigl(
\partial^2_{
z_u
z_v
}
g(z)
+
\partial_{
z_u
}
g(z)
\partial_{
z_v
}
g(z)
\Bigr)
\Big\vert_{
z=Y(1,\alpha^{(1)})
}
\\
&
\quad\quad
-
\partial_{
z_{u}
}
g(z)
\vert_{
z=Y(1,\alpha^{(1)})
}
\partial_{
z_{v}
}
g(z)
\vert_{
z=Y(1,\alpha^{(2)})
}
\Bigr\}
\Big\vert_{
Q^{(k)}
=
Q^{(k)} + s Q
}
\Bigr]
\Bigr]
.
\end{align*}
Note that
\begin{align*}
\partial_s
\left(
Q(\alpha^{(1)},\alpha^{(2)})_{u,v}
\vert_{
Q^{(k)}
=
Q^{(k)} + s Q
}
\right)
=
\begin{cases}
Q_{u,v}
,
&
q_\text{L}(\alpha^{(1)},\alpha^{(2)})
=
k
,
\\
0
,
&
q_\text{L}(\alpha^{(1)},\alpha^{(2)})
\neq
k
.
\end{cases}
\end{align*}
\begin{enumerate}
\item 
Define $M(\alpha^{(1)},\alpha^{(2)}) \in \R^{d \times d}$ as
\begin{align*}
M(\alpha^{(1)},\alpha^{(2)})_{u,v}
\equiv
\partial_{
z_{u}
}
g(z)
\vert_{
z=Y(1,\alpha^{(1)})
}
\partial_{
z_{v}
}
g(z)
\vert_{
z=Y(1,\alpha^{(2)})
}
.
\end{align*}
Hence, we arrive at
\begin{align*}
\partial_{Q^{(k)} \leadsto Q}
f_\rho
(0,0)
=
-
\frac{1}{2}
\E
\Bigl[
\sum_{
\alpha^{(1)}
\alpha^{(2)}
\in
\mathcal{A}
}
\I_{
q_\text{L}(\alpha^{(1)},\alpha^{(2)})
=
k
}
\xi(\alpha^{(1)})
\xi(\alpha^{(2)})
(\alpha^{(1)},\alpha^{(2)})
\langle
Q
,
M(\alpha^{(1)},\alpha^{(2)})
\rangle
\Bigr]
.
\end{align*}
The proof is concluded similarly to the proof of
Theorem~\ref{thm:lectire-03-an-analytic-projection-the-v-term-without-rpc} by
using the properties of the RPC
(Theorem~\ref{thm:remainder:some-limiting-grem-properties} and
Lemma~\ref{lemma:vektor-sk-squares-of-rpc}).
\item
The proof is the same as in (1).
\end{enumerate}
\end{proof}
The following is a multidimensional version of
\cite[Lemma~4.3]{TalagrandParisiFormula2006}.
\begin{lemma}
For any $\alpha \in \mathcal{A}$, we have
\begin{enumerate}
\item 
\begin{align*}
\partial_{x_k}
f_\rho(0,0)
\vert_{
x_k = x_{k-1} 
}
=
\frac{1}{x_{k-1}}
\E
\Bigl[
&
T_1(\alpha)
\cdots
T_{k-2}(\alpha)
T_{k-1}(\alpha)
\vert_{
x_k = x_{k-1}
}
\\
&
\Bigl(
\E
\Bigl[
f(x_{k+1}, Y(x_{k+1}, \alpha))
T_k(\alpha)
\vert_{
x_k = x_{k-1}
}
\Bigr]
-
f(x_{k}, Y(x_{k}, \alpha))
\Bigr)
\Bigr]
.
\end{align*}
\item
Let 
$M \in \symmetric(d)$ with
$
M_{u,v} \equiv \partial_{z_u}
f(x_{k}, Y(x_{k}, \alpha))
\partial_{z_v}
f(x_{k}, Y(x_{k}, \alpha))
$, then
\begin{align*}
\partial^2_{Q^{(k)} \leadsto Q, x_k}
f_\rho(0,0)
=
\frac{1}{2}
\E
\Bigl[
T_1(\alpha)
\cdots
T_{k-2}(\alpha)
\langle
Q
,
M
\rangle
\Bigr]
.
\end{align*}
\end{enumerate}

\end{lemma}
\begin{proof}
This proof is the same as in \cite{TalagrandParisiFormula2006}.
\end{proof}
We now generalise the PDE \eqref{vektor-sk-parisis-pde}.
Given a piece-wise continuous $x \in \mathcal{Q}(1,1)$ and 
$
Q \in
\mathcal{Q}(U,d)
$, 
consider the following terminal value problem
\begin{align}
\label{eq:remainder:generalised-parisi-pde}
\begin{cases}
\partial_t f
+
\frac{1}{2}
\left(
\langle
\dot{Q}
,
\nabla^2
f
\rangle
+
x
\langle
\dot{Q}
\nabla
f
,
\nabla
f
\rangle
\right)
=
0
,
&
(y,t)
\in
\R^d \times (0,1)
,
\\
f(y,1)
=
g(y)
.
\end{cases}
\end{align}
We say that $f \in C([0;1] \times \R^d \to \R)$ is a piece-wise viscosity
solution of \eqref{eq:remainder:generalised-parisi-pde}, if there exists the partition of the
unit segment $0 =: x_0 < x_1 < \ldots < x_{n+1} \equiv 1$ such that, for each $k \in [0,n] \cap \N$, $f: (x_k;x_{k+1}) \times \R^d \to \R$ is a
viscosity solution (see, e.g.,  \cite{Briand2007})
of
\begin{align*}
\begin{cases}
\partial_t f
+
\frac{1}{2}
\left(
\langle
\dot{Q}
,
\nabla^2
f
\rangle
+
x
\langle
\dot{Q}
\nabla
f
,
\nabla
f
\rangle
\right)
=
0
,
&
(y,t)
\in
\R^d \times (x_k,x_{k+1})
,
\\
f(y,x_{k+1}+0)
=
f(y,x_{k+1}-0)
,
\\
f(y,1)
=
g(y)
.
\end{cases}
\end{align*}
\begin{proposition}
\label{eq:remainder:parisi-functional-lipschitzianity}
For any
$
\rho^{(1)}, \rho^{(2)} 
\in 
\mathcal{Q}^\prime(U,d)
$,
we have
\begin{align*}
\vert
f_{
\rho^{(1)}
}
(0,0)
-
f_{
\rho^{(2)}
}
(0,0)
\vert
\leq
\frac{C}{2}
\int_0^1
\Vert
\rho^{(1)}(t)
-
\rho^{(2)}(t)
\Vert_{\text{F}}
\dd t
,
\end{align*}
where 
$
C 
=
C(\Sigma)
\equiv 
\E
\left[
\Vert
M
\Vert_{\text{F}}
\right]
$.
\end{proposition}
\begin{proof}
This is an adaptation of the proof of
\cite[Theorem~3.1]{TalagrandParisiMeasures2006a} to the  multidimensional case.
Assume without loss of generality that the paths $\rho^{(1)}$ and $\rho^{(2)}$
have same jump times $\{x_k\}_{k=0}^{n+1}$. Denote the corresponding overlap
matrices as $\{Q^{(1,k)}\}_{k=0}^{n+1}$ and $\{Q^{(2,k)}\}_{k=0}^{n+1}$. Given
$s \in [0;1]$, define the new path $\rho(s) \in \mathcal{Q}^{\prime}(U,d)$ by
assuming that it has the same jump times $\{x_k\}_{k=0}^{n+1}$ as the paths
$\rho^{(1)}, \rho^{(2)}$ and defining its overlap matrices as $Q^{(k)}(s) \equiv s
Q^{(1,k)} + (1-s) Q^{(2,k)}$. On the one hand, we readily have
\begin{align*}
\int_0^1
\Vert
\rho^{(1)}(t)
-
\rho^{(2)}(t)
\Vert_{\text{F}}
\dd t
=
\sum_{k=1}^n
(x_k - x_{k-1})
\Vert
Q^{(1,k)}
-
Q^{(2,k)}
\Vert_{\text{F}}
.
\end{align*}
On the other hand, using
Lemma~\ref{lem:remainder:parisi-formula-differentiation}, we have
\begin{align*}
\vert
\partial_s
f_{
\rho(s)
}
(0,0)
\vert
\leq
\frac{C}{2}
\sum_{k=1}^n
(x_k-x_{k-1})
\Vert
Q^{(1,k)}
-
Q^{(2,k)}
\Vert_{\text{F}}
.
\end{align*}
Finally, we have
\begin{align*}
\vert
f_{
\rho^{(1)}
}
(0,0)
-
f_{
\rho^{(2)}
}
(0,0)
\vert
\leq
\int_0^1
\vert
\partial_s
f_{
\rho(s)
}
(0,0)
\vert
\dd s
.
\end{align*}
Combining the last three formulae, we get the theorem.
\end{proof}
\begin{remark}
Note that using the same argument and notations as in the previous theorem we
get that, for any $(y,t) \in \R^d \times [0;1]$,
\begin{align*}
\vert
f_{
\rho^{(1)}
}
(y,t)
-
f_{
\rho^{(2)}
}
(y,t)
\vert
\leq
\frac{C(\Sigma)}{2}
\int_t^1
\Vert
\rho^{(1)}(s)
-
\rho^{(2)}(s)
\Vert_{\text{F}}
\dd s
.
\end{align*}
\end{remark}
\begin{remark}
Note that we can associate to each $\rho \in \mathcal{Q}(U,d)$
a $\symmetric^+(d)$-valued countably additive vector measure 
$
\nu_\rho
\in
\mathcal{M}([0;1], \symmetric^+(d))
$ 
by the following standard
procedure. Given $[a;b) \subset [0;1]$, define 
\begin{align*}
\nu_\rho([a;b))
\equiv
\rho(b)-\rho(a)
\end{align*}
and then extend the measure, e.g., to all Borell subsets of $[0;1]$.
\end{remark}
\begin{theorem}
\label{thm:remainder:extension-by-continuity-of-the-parisi-functional}
Given $U \in \symmetric^+(d)$, we have
\begin{enumerate}
\item 
The set $\mathcal{Q}(U,d)$ is compact under the topology induced by
the following norm
\begin{align}
\label{eq:remainder:order-parameter-space-norm}
\Vert
\rho
\Vert
\equiv
\int_0^1
\Vert
\rho(t)
\Vert_{\text{F}}
\dd t
,
\quad
\rho \in \mathcal{Q}(U,d)
.
\end{align}
\item
The functional 
$
\mathcal{Q}^\prime(U,d)
\ni
\rho 
\mapsto
f_\rho(0,0)
$
is Lipschitzian and can be uniquely extended by continuity to the whole
$\mathcal{Q}(U,d)$. 
\end{enumerate}
\end{theorem}
\begin{proof}
\begin{enumerate}
\item
The topology induced by the norm
\eqref{eq:remainder:order-parameter-space-norm} coincides with the topology of
weak convergence of the above-defined vector measures. Since $\mathcal{Q}(U,d)$
is a bounded set, it is compact in the weak topology.
\item 
This is an immediate consequence of
Proposition~\ref{eq:remainder:parisi-functional-lipschitzianity}.
\end{enumerate}
\end{proof}
In the next result, we summarise some results on the PDE
\eqref{eq:remainder:generalised-parisi-pde} for the non-discrete parameters,
cf. Proposition~\ref{proposition:vektor-sk-parisis-pde-descrete-case}.
\begin{theorem}
\label{thm:remainder:piece-wise-viscosity-solutions-to-parisis-pde}
\mbox{}
\begin{enumerate}
\item
\emph{Existence.}
Assume that $Q$ is in $\mathcal{Q}(U,d)$ and is piece-wise $C^{(1)}$. Assume
also that $x$ is in  $\mathcal{Q}(1,1)$ and is piece-wise continuous. Then the
terminal value problem \eqref{eq:remainder:generalised-parisi-pde} has a unique continuous, piece-wise viscosity solution 
$
f_{Q,x} \in C([0;1] \times \R^d)
$.
\item
\emph{Monotonicity with respect to $x$.}
Assume $Q \in \mathcal{Q}(U,d)$. Assume also that $x^{(1)},x^{(2)} \in
\mathcal{Q}(1,1)$ are such that $x^{(1)}(t) \leq x^{(2)}(t)$, almost everywhere
for $t \in [0;1]$. Let $f_{Q,x^{(1)}}$ and $f_{Q, x^{(2)}}$ be the
corresponding solutions of \eqref{eq:remainder:generalised-parisi-pde}. Then $f_{Q,x^{(1)}} \leq f_{Q,
x^{(2)}}$.
\item
\emph{Monotonicity with respect to $g$.} Assume $g_1, g_2: \R^d \to \R$ satisfy
Assumption~\ref{asm:remainder:end-condition} and also $g_1 \leq g_2$ almost
everywhere. Let $f_{g_1}, f_{g_2}: \R^d \times [0;1] \to \R$ be the
corresponding solutions of \eqref{eq:remainder:generalised-parisi-pde} with
$g = g_1$, $g = g_2$, respectively. Then 
$
f_{g_1} \leq f_{g_2}
$.
\end{enumerate}
\end{theorem}
\begin{proof}
\begin{enumerate}
\item 
Due to the assumptions, the diffusion matrix 
$
\dot{Q}(t)
=
\dot{\rho}(t)
$ 
in \eqref{eq:remainder:generalised-parisi-pde}
is non-negative definite. Applying 
\cite[Proposition~8]{Briand2007} to the PDE
\eqref{eq:remainder:generalised-parisi-pde} successively on the intervals
$[x_k;x_{k+1})$, where the 
$
\dot{\rho}
$ 
is continuous, gives the existence of the solutions in viscosity sense and,
moreover, gives their continuity. Uniqueness is ensured by
\cite[Theorem~1.1]{DaLioLey2006}.
\item
By the approximation argument (cf.
Theorem~\ref{thm:remainder:extension-by-continuity-of-the-parisi-functional}),
it is enough to assume that $x^{(1)},x^{(2)} \in \mathcal{Q}^\prime(1,1)$ and
$Q \in \mathcal{Q}^\prime(U,d)$.
Then Proposition~\ref{proposition:vektor-sk-parisis-pde-descrete-case} gives
the existence of the corresponding piece-wise classical solutions of
\eqref{eq:remainder:generalised-parisi-pde}: $f_{Q,x^{(1)}}, f_{Q,x^{(2)}}$. 
These solutions are obviously also the (unique) piece-wise viscosity solutions
of \eqref{eq:remainder:generalised-parisi-pde}. The comparison result
\cite[Theorem~5]{Briand2007} and the non-linear Feynman-Kac formula
\cite[Proposition~8 ]{Briand2007} give then the claim.
\item
This can be seen either from the representation
\eqref{eq:remainder:parisi-rpc-representation} and an approximation argument, or
exactly as in (2) by invoking the results of \cite{Briand2007}.
\end{enumerate}
\end{proof}

\subsection{The Parisi functional}
We consider now a specific terminal condition in the system
\eqref{vektor-sk-parisis-pde} given in
\eqref{eq:remainder:multidim-sk-final-value}.

Given $\rho \in \mathcal{Q}(U,d)$, let $f_\rho : [0;1]
\times \R^d \to \R$ be the value of (the continuous extension onto
$\mathcal{Q}(U,d)$  of) the solution of \eqref{vektor-sk-parisis-pde} with the
specific terminal condition given by \eqref{eq:remainder:multidim-sk-final-value}. Following the
ideas in the physical literature, we now define the \emph{Parisi functional} 
$
\mathcal{P}(\beta, \rho, \Lambda): 
\R_+ \times \mathcal{Q}^{\prime}(U, d) \times \symmetric^+(d) \times
\symmetric(d) \to \R $
in as
\begin{align}
\label{vektor-sk-pasisis-functional}
\mathcal{P}(\beta, \rho, \Lambda)
\equiv 
f_\rho(0,0)
-
\frac{
\beta^2
}{
2
}
\int_{0}^1
x(t)
\dd
\left(
\Vert
\rho(t)
\Vert_\text{F}^2
\right)
-
\langle
U,
\Lambda
\rangle
.
\end{align}
The integral in
\eqref{vektor-sk-pasisis-functional} is understood in the usual
Lebesgue-Stiltjes sense.
\begin{remark}
Note that the path integral term in \eqref{vektor-sk-pasisis-functional} equals
$f(0,0)$, where $f(t,y)$ is the solution of
\eqref{eq:remainder:generalised-parisi-pde} with the following boundary condition
\begin{align*}
g(y)
\equiv
\beta
\langle
y
,
\I
\rangle
=
\beta
\sum_{u=1}^d
y_u
,
\quad
y \in \R^d
.
\end{align*}
\end{remark}
Obviously
$\mathcal{Q}^\prime(d)$ is dense in $\mathcal{Q}(d)$.
\begin{theorem}
We have
\begin{align}
\label{eq:multidim-sk:upper-bound-pde}
p(\beta)
\leq
\sup_{
\substack{
U
\in 
\symmetric^+(d)
}
}
\inf_{
\substack{
\rho
\in
\mathcal{Q}^\prime(U,d)
\\
\Lambda
\in
\symmetric(d)
}
}
\mathcal{P}(\beta, \rho, \Lambda)
.
\end{align}
\end{theorem}
\begin{proof}
The bound
\eqref{eq:multidim-sk:upper-bound-pde} is a straightforward consequence of
Theorem~\ref{thm:multidimensional-sk:free-energy-upper-bound}. 
\end{proof}

\subsection{On strict convexity of the Parisi functional and its variational
representation} 
\label{sec:remainder:strict-convexity-of-the-parisi-functional-and-the-variational-representation}

In this subsection, we derive a variational representation for Parisi's
functional. As a consequence, for $d=1$, we prove that the functional is
strictly convex with respect to the $x \in \mathcal{Q}(1,1)$, if the terminal condition
$g$ (cf. \eqref{eq:remainder:generalised-parisi-pde}) is strictly
convex and increasing. This result is related to the problem of strict
convexity of the Parisi functional in the case of the SK model.

Let 
$
W \equiv \{ W(s) \}_{s \in \R_+}
$ 
be the standard $\R^d$-valued Brownian motion
and let 
$
\{
\mathcal{F}_t
\}_{t \in \R_+}
$ be the correspondent filtration.
Define
\begin{align*}
\mathcal{U}[t;T]
\equiv
\{
u: [t;T] \to \R^d
\mid
\text{
$u$ is $\{\mathcal{F}_t\}_{t \in \R_+}$ progressively measurable
}
\}
.
\end{align*}
Given 
$
u \in \mathcal{U}[t;1]
$, 
$
Q \in \mathcal{Q}(U,d)
$
and
$
x \in \mathcal{Q}(1,1)
$,
consider the following $\R^d$-valued and adapted to 
$
\{\mathcal{F}_t\}_{t \in \R^+}
$
diffusion
\begin{align*}
Y^{(Q,x,u,t,y)}(s)
=
y
-
\int_t^s
\left(
x(s) 
\dot{Q}(s)
\right)^{1/2}
u(s)
\dd s
+
\int_t^s
\left(
\dot{Q}(s)
\right)^{1/2}
\dd W(s)
,
\quad
s \in [t;1]
.
\end{align*}
Given some function $g: \R^d \to \R$ satisfying
Assumption~\ref{asm:remainder:end-condition}, define 
$
f_{Q,x}: \R^d \times [0;1] \to \R
$
as 
\begin{align}
\label{eq:remainder-variational-representation-of-parisis-functional}
f_{Q,x}(y,t)
\equiv
\sup_{
u \in \mathcal{U}[t;1]
}
\E
\left[
g(
Y^{(Q,x,u,t,y)}(1)
)
-
\frac{1}{2}
\int_t^1
\Vert
u(s)
\Vert_2^2
\dd s
\right]
.
\end{align}
\begin{proposition}
\label{eq:remainder:strict-convexity-of-the-controlled-diffusion}
Let $d=1$. If $g$ is strictly convex and increasing, then
the functional $\mathcal{Q}(1,1) \ni x \mapsto f_{Q,x}$ is strictly convex.
\end{proposition}
\begin{proof}
We have
\begin{align*}
Y^{(Q,x,u,t,y)}(1)
=
y
-
\int_t^1
\left(
x(s) 
\dot{Q}(s)
\right)^{1/2}
u(s)
\dd s
+
\int_t^1
\left(
\dot{Q}(s)
\right)^{1/2}
W(s)
.
\end{align*}
By an approximation argument, it is enough to prove the strict convexity for the
continuous $x_1, x_2 \in \mathcal{Q}(1,1)$ ($x_1 \neq x_2$). For any $\gamma
\in (0;1)$, we have
\begin{align}
\label{eq:remainder:y-convexity}
Y^{(Q,\gamma x_1 + (1-\gamma) x_2,u,t,y)}(1)
&
=
-
\int_t^1
\left(
\gamma x_1 + (1-\gamma) x_2 
\dot{Q}(s)
\right)^{1/2}
u(s)
\dd s
+
\int_t^1
\left(
\dot{Q}(s)
\right)^{1/2}
W(s)
\nonumber
\\
&
<
-
\gamma
\int_t^1
\left(
x_1 
\dot{Q}(s)
\right)^{1/2}
u(s)
\dd s
-
(1-\gamma)
\int_t^1
\left(
x_2
\dot{Q}(s)
\right)^{1/2}
u(s)
\dd s
\nonumber
\\
&
\quad
+
\int_t^1
\left(
\dot{Q}(s)
\right)^{1/2}
W(s)
\nonumber
\\
&
=
\gamma
Y^{(Q,x_1,u,t,y)}(1)
+
(1-\gamma)
Y^{(Q,x_2,u,t,y)}(1)
,
\end{align}
where the strict inequality above is due to the strict concavity of the square
root function. The strict convexity and monotonicity of $g$ combined with the
representation \eqref{eq:remainder:y-convexity} implies that
\eqref{eq:remainder-variational-representation-of-parisis-functional} is
strictly convex as a function of $x$, since a supremum of a family of convex
functions is convex.
\end{proof}
\begin{proposition}
\label{prop:remainder:hjb-for-value-function}
Given a piece-wise continuous $x \in \mathcal{Q}(1,1)$ and a 
$
Q \in
\mathcal{Q}(U,d)
$ 
which is piece-wise in $C^1(0;1)$, the function 
$
f_{Q,x}:
\R^d \times [0;1] \to \R
$ 
defined by
\eqref{eq:remainder-variational-representation-of-parisis-functional} is a
unique, continuous, piece-wise viscosity solution of the following terminal
value problem
\begin{align*}
\begin{cases}
\partial_t f
+
\frac{1}{2}
\left(
\langle
\dot{Q}
,
\nabla^2
f
\rangle
+
x
\langle
\dot{Q}
\nabla
f
,
\nabla
f
\rangle
\right)
=
0
,
&
(y,t)
\in
\R^d \times (0,1)
,
\\
f(y,1)
=
g(y)
.
\end{cases}
\end{align*}
\end{proposition}
\begin{proof}
In a way similar to the proof of
Theorem~\ref{thm:remainder:piece-wise-viscosity-solutions-to-parisis-pde}, we
successively use \cite[Theorem~2.1]{DaLioLey2006} on the intervals
$(x_k;x_{k+1})$, where the data of the PDE are continuous.
\end{proof}

\begin{theorem}
\label{thm:remainder:strict-convexity-parisi-functional-1d}
Assume $d=1$. Suppose also that $g$ satisfies the assumptions of
Proposition~\ref{eq:remainder:strict-convexity-of-the-controlled-diffusion}.
For any $u \in \R$, the generalised Parisi functional given by
\eqref{vektor-sk-pasisis-functional} with $f_\rho(0,0)$ corresponding to the
terminal condition $g$ is strictly convex on $Q(u,1)$. Consequently, there
exists a unique optimising order parameter.
\end{theorem}
\begin{proof}
In 1-D, we can choose the coordinates such that $Q \equiv U t$, on $[0;1]$.
Consequently, $\dot{Q} \equiv U \equiv \text{const}$ on $[0;1]$. Hence, it is
enough check the strict convexity with respect to $x \in \mathcal{Q}(1,1)$. The result follows by approximation in the
norm \eqref{eq:remainder:order-parameter-space-norm} of an arbitrary pair of
different elements of $\mathcal{Q}(U,d)$ by a pair of
elements of $\mathcal{Q}^\prime(U,d)$ and
Propositions~\ref{proposition:vektor-sk-parisis-pde-descrete-case},
\ref{eq:remainder:strict-convexity-of-the-controlled-diffusion} and \ref{prop:remainder:hjb-for-value-function}.
\end{proof}
\begin{remark}
Due to the monotonicity assumption on $g$,
Theorem~\ref{thm:remainder:strict-convexity-parisi-functional-1d} does not
cover the case of the SK model, where the terminal value $g$ is given by
\eqref{eq:remainder:multidim-sk-final-value}.
\end{remark}
\subsection{Simultaneous diagonalisation scenario}
\label{sec:remainder:simultaneous-diagonalisation-scenario}

In the setups with highly symmetric state spaces $\Sigma_N$ (such as the
spherical spin models of \cite{PanchenkoTalagrandMultipleSKModels2006} or
the Gaussian spin models, see
Section~\ref{sec:gaussian-spins:gaussian-spins} below), less complex order
parameter spaces as $Q(U,d)$ suffice.

Given some orthogonal matrix 
$
O \in \mathcal{O}(d)
$, 
we briefly discuss the case 
$
\rho \in \mathcal{Q}_{\text{diag}}(U,O,d)
$, 
where
\begin{align*}
\mathcal{Q}_{\text{diag}}(U,O,d)
\equiv
\{
\rho
\in 
\mathcal{Q}(U,d)
\mid
\text{for all $t \in [0;1]$, the matrix $O\rho(t)O^*$
is diagonal}
\}
.
\end{align*}
The space
$\mathcal{Q}_{\text{diag}}(U,O,d)$ is obviously isomorphic to the space of
``paths'' with the non-decreasing coordinate functions in $\R^d$, starting from
the origin and ending at $u$, i.e.,
\begin{align*}
\bar{\mathcal{Q}}(u, d)
\equiv 
\{
\rho : [0;1] \to \R^{d}
\mid
\bar{\rho}(0)=0
;
\bar{\rho}(1)=u
;
\bar{\rho}(t) 
\preceq
\bar{\rho}(s)
,
\text{ for }
t \leq s 
;
\bar{\rho}
\text{ is c\'adl\'ag}
\}
,
\end{align*}
where  $u = O U O^* \in \R^d$. The isomorphism is then given by
\begin{align}
\label{eq:remainder:diagonalisation-isomorphism}
\bar{\mathcal{Q}}(u,d) 
\ni
\bar{\rho} 
\mapsto
O \rho O^*
\in
\mathcal{Q}_{\text{diag}}(U,O,d)
.
\end{align}

\section{Remainder estimates}
\label{sec:remainder:talagrands-remainder-estimates}

In this section, we partially extend Talagrand's remainder estimates to the
multidimensional setting. 
Due to
Proposition~\ref{eq:remainder:guerras-upper-and-lower-bounds}, to prove the
validity of Parisi's formula it is enough to show that all the $\mu_k$ terms in \eqref{eq:vektor-sk-the-convenient-form-of-the-remainder} almost
vanish for the almost optimal parameters of the optimisation problem
in \eqref{eq:multidimensional-sk:free-energy-upper-bound}. This can be done if the free energy of two coupled replicas of the system
\eqref{eq:lecture-04-special-free-energy-like-quantity} is strictly smaller
than twice the free energy of the uncoupled single system
\eqref{eq:as2-phi-of-t-definition}, see
inequality \eqref{thm:as2-a-sufficient-condition-for-mu-k-to-vanish}.
However, the systems involved in
\eqref{thm:as2-a-sufficient-condition-for-mu-k-to-vanish} are effectively at
least as complex as the SK model itself. In
Section~\ref{upper-bounds-on-phi-2-guerras-scheme-revisited}, we again apply
Guerra's scheme to obtain the upper bounds on \eqref{eq:lecture-04-special-free-energy-like-quantity} in terms of the free energy of the corresponding comparison GREM-inspired model. One might then
hope that by a careful choice of the comparison model one can prove inequality \eqref{thm:as2-a-sufficient-condition-for-mu-k-to-vanish}. In
Sections~\ref{sec:remainder:adjustment-of-the-upper-bounds-on-varphi} and
\ref{sec:a-priori-estimates}, we formulate some conditions on the comparison
system which would suffice to get inequality
\eqref{thm:as2-a-sufficient-condition-for-mu-k-to-vanish}, giving, hence, the
conditional proof of the Parisi formula, see 
Theorem~\ref{thm:remainder:the-conditional-parisi-formula}.

\subsection{A sufficient condition for $\mu_k$-terms to vanish}

In this subsection, we are going to establish a sufficient condition for the
measures $\mu_k$ to vanish. This condition states roughly the following. Whenever
the free energy of a certain replicated system uniformly in $N$
strictly less then twice the free energy of the single system, the measure
$\mu_k$ vanishes in $N \to +\infty$ limit (see
Lemma~\ref{thm:as2-a-sufficient-condition-for-mu-k-to-vanish}).

Keeping in mind the definition of $\mu_k$ (cf. 
\eqref{eq:lecture-04-mu-k-definition})
and of the Hamiltonian $H_t(\sigma,\alpha)$ (cf.
\eqref{eq:remainder:interpolating-hamiltonian}), we define, for $\alpha^{(1)},\alpha^{(2)} \in \mathcal{A}^{(2),k}$, the
corresponding replicated Hamiltonian as
\begin{align}
\label{eq:as2-double-replicated-rest-term-hamiltonian}
H^{(2)}_t(\sigma^{(1)},\sigma^{(2)},\alpha^{(1)},\alpha^{(2)})
\equiv
H_t(\sigma^{(1)}, \alpha^{(1)})
+
H_t(\sigma^{(2)}, \alpha^{(2)})
.
\end{align}
\begin{remark}
We note here that the distribution of the Hamiltonian $H_t(k,\sigma^{(1)},\sigma^{(2)})$
depends only on $k$ and not on the choice of the indices
$\alpha^{(1)},\alpha^{(2)} \in \mathcal{A}^{(2),k}$.
\end{remark}
\begin{remark}
The superscript $(2)$ in \eqref{eq:as2-double-replicated-rest-term-hamiltonian}
(and in what follows) indicates that the quantity is related to the twice replicated objects.
\end{remark}
Define
\begin{align*}
\mathcal{A}^{(2),k}
\equiv
\{
(\alpha^{(1)},\alpha^{(2)})
\in
\mathcal{A}^2
:
q_{\text{L}}(\alpha^{(1)},\alpha^{(2)})=k
\}
.
\end{align*}
Additionally, for any
$
\mathcal{V}
\subset 
\Sigma(B(U,\eps))^2
$
and any suitable Gaussian process, 
\begin{align*}
\{
F(\sigma^{(1)},\sigma^{(2)},\alpha^{(1)},\alpha^{(2)})
:
\sigma^{(1)},\sigma^{(2)}
\in
\Sigma_N
,
\alpha^{(1)},\alpha^{(2)}
\in
\mathcal{A}
\}
,
\end{align*}
we define the \emph{local remainder comparison functional} as
\begin{align}
\label{eq:as2-phi-2-k-x}
\Phi^{(2),k,x}_{
\mathcal{V}
}
\left[
F
\right]
\equiv
\frac{1}{N}
\E
\Bigl[
\log
&
\iint_{
\mathcal{V}
}
\iint_{
\mathcal{A}^{(2),k}
}
\exp
\left\{
\beta
\sqrt{N}
F(\sigma^{(1)},\sigma^{(2)},\alpha^{(1)},\alpha^{(2)})
\right\}
\nonumber
\\
&
\dd 
\mu^{\otimes N}(\sigma^{(1)})
\dd 
\mu^{\otimes N}(\sigma^{(2)})
\dd \xi(\alpha^{(1)})
\dd \xi(\alpha^{(2)})
\Bigr]
.
\end{align}
Define
\begin{align}
\label{eq:lecture-04-special-free-energy-like-quantity}
\varphi^{(2)}_{N}
(k,t,x,Q,\mathcal{V})
\equiv
\Phi^{(2),k}_{
\mathcal{V}
}
\left[
H^{(2)}_t
\right]
.
\end{align}
\begin{lemma}
Recalling the definition \eqref{eq:as2-phi-of-t-definition}, for any
$
\mathcal{V}
\subset
\Sigma(B(U,\eps))^2
$, we
have 
\begin{align}
\label{eq:as2-non-strict-inequality-for-the-coupled-rest-term-free-energy}
\varphi^{(2)}_N
(
k
,
t
,
x,Q,
\mathcal{V}
) 
\leq
\varphi^{(2)}_N
(
k
,
t
,
x,Q,
\Sigma(B(U,\eps))^2
)
=
2 \varphi_N
(
t
, 
x,Q,
\Sigma(B(U,\eps))
)
.
\end{align}
\end{lemma}
\begin{proof}
The first inequality in
\eqref{eq:as2-non-strict-inequality-for-the-coupled-rest-term-free-energy} is obvious, since the expression under the integral in 
\eqref{eq:as2-phi-2-k-x} is positive. 
The equality in
\eqref{eq:as2-non-strict-inequality-for-the-coupled-rest-term-free-energy} is
an immediate consequence of the RPC averaging property 
\eqref{eq:remainder:averaging-property}.
\end{proof}
In what follows, we shall be looking for the sharper (in particular, 
\emph{strict}) versions of the inequality
\eqref{eq:as2-non-strict-inequality-for-the-coupled-rest-term-free-energy}
because of the following observation due to Talagrand
\cite{TalagrandParisiFormula2006}.
\begin{lemma}
\label{thm:as2-a-sufficient-condition-for-mu-k-to-vanish}
Fix an arbitrary
$
\mathcal{V}
\subset
\Sigma_N(B(U,\eps))^2
$.
Suppose that, for some $\eps > 0$, the following inequality holds 
\begin{align}
\label{eq:as2-strict-inequality-for-the-coupled-rest-term-free-energy}
\varphi^{(2)}_N(k,t,x,Q,\mathcal{V})
\leq
2 \varphi_N
(
t
, 
x,Q,
\Sigma_N(B(U,\eps))
)
-
\eps
.
\end{align}
Then, for some $K>0$, we have
\begin{align*}
\mu_k
(
\mathcal{V}
)
\leq
K \exp 
\left(
-
\frac{N}{K}
\right)
.
\end{align*}
\end{lemma}
\begin{proof}
The proof is based on Theorem
\ref{vektor-sk-concentration-of-measure-for-free-energies} and follows the
lines of \cite[Lemma~7]{panchenko-free-energy-generalized-sk-2005}.
\end{proof}

\subsection{Upper bounds on $\varphi^{(2)}$: Guerra's scheme revisited}
\label{upper-bounds-on-phi-2-guerras-scheme-revisited}
In this subsection, we shall develop a mechanism to obtain upper bounds on 
$\varphi^{(2)}$ defined in
\eqref{eq:lecture-04-special-free-energy-like-quantity}. This will be achieved
in the full analogy to Guerra's scheme by using a suitable Gaussian
comparison system.

\label{sec:a-comparison-trajectory-on-the-doubled-path-space}
Given $U \in \symmetric^+(d)$, we say that $V \in \R^{d \times d}$ is an
\emph{admissible mutual overlap matrix for $U$}, if
\begin{align}
\label{eq:as2-mathfrak-u-definition}
\mathfrak{U}
\equiv
\begin{bmatrix}
U 
& 
V
\\
V^* 
&
U
\end{bmatrix}
\in
\symmetric^+(2d)
.
\end{align}
Furthermore, define 
\begin{align*}
\mathcal{V}(U)
\equiv
\{
V \in \R^{d \times d}
:
\text{$V$ is an admissible mutual overlap matrix for $U$}
\}
.
\end{align*}

Hereinafter without further notice we assume that $\mathfrak{U} \in
\symmetric^+(2d)$ has the form \eqref{eq:as2-mathfrak-u-definition}, where $V$ is some admissible mutual
overlap matrix for $U$.

Let 
$
\mathfrak{Q} \in \mathcal{Q}(\mathfrak{U},2d)
$.
Let 
$
\mathfrak{x} 
\equiv
\{
\mathfrak{x}_l 
\in
[0;1]
\}_{l=1}^{\mathfrak{n}}
$
be the ``jump times'' of the path $\varrho$. We assume that the ``times'' are
increasingly ordered, i.e., 
\begin{align*}
0
=
\mathfrak{x}_0
<
\mathfrak{x}_1
<
\ldots
<
\mathfrak{x}_{\mathfrak{n}}
<
\mathfrak{x}_{\mathfrak{n}+1}
=
1
.
\end{align*}
Consider the following collection of matrices
\begin{align*}
\mathfrak{Q}
\equiv
\{
\mathfrak{Q}_l
\equiv
\mathfrak{Q}(\mathfrak{x}_l)
\subset 
\symmetric^+(2d)
\}_{l=0}^{\mathfrak{n}+1}
.
\end{align*}
We obviously then have
\begin{align}
\label{eq:remainder:remainder-comparison-matrices-monotonicity}
0
=
\mathfrak{Q}^{(0)}
\prec
\mathfrak{Q}^{(1)}
\prec
\ldots
\prec
\mathfrak{Q}^{({\mathfrak{n}})}
\prec
\mathfrak{Q}^{({\mathfrak{n}+1})}
=
\mathfrak{U}
.
\end{align}
Such a path $\mathfrak{Q}$ induces in the usual way the ``doubled''
GREM overlap kernel 
$
\label{eq:as2-mathfrak-q-kernel-definition}
\mathfrak{Q}
\equiv
\{
\mathfrak{Q}(\alpha^{(1)},\alpha^{(2)})
\in 
\symmetric^+(2d) 
\mid
\alpha^{(1)},\alpha^{(2)}
\in
\mathcal{A}_{\mathfrak{n}}
\}
$,
defined as
\begin{align*}
\mathfrak{Q}(\alpha^{(1)},\alpha^{(2)})
\equiv
\mathfrak{Q}^{(q_\text{L}(\alpha^{(1)},\alpha^{(2)}))}.
\end{align*}
We also need the $d \times d$ submatrices of the above overlap such that
\begin{align}
\label{eq:remainder:mathfrak-q-submatrices}
\mathfrak{Q}(\alpha^{(1)},\alpha^{(2)})
=
\begin{bmatrix}
\mathfrak{Q}\vert_{11}(\alpha^{(1)},\alpha^{(2)})
& 
\mathfrak{Q}\vert_{12}(\alpha^{(1)},\alpha^{(2)})
\\
\mathfrak{Q}\vert_{12}(\alpha^{(1)},\alpha^{(2)})^*
&
\mathfrak{Q}\vert_{22}(\alpha^{(1)},\alpha^{(2)})
\end{bmatrix}
.
\end{align}
\begin{remark}
For $\sigma^{(1)}  \sigma^{(2)} \in \Sigma_N$, we shall use the notation
$
\sigma^{(1)} 
\shortparallel 
\sigma^{(2)} 
\in 
\left(
\R^{2d}
\right)^N
$
to denote the vector obtained by
the following concatenation of the vectors $\sigma^{(1)}$ and $\sigma^{(2)}$
\begin{align*}
\sigma^{(1)} 
\shortparallel 
\sigma^{(2)} 
\equiv
\left(
\sigma^{(1)}_i
\sigma^{(2)}_i
\in 
\Sigma
\times
\Sigma
\subset
\R^{2d}
\right)_{
i=1
}^{
N
}
.
\end{align*}
\end{remark}
Let us observe that the process 
\begin{align*}
X^{(2)}
\equiv
\left\{
X^{(2)}(\tau)
=
X(\sigma^{(1)})
+
X(\sigma^{(1)})
\mid
\tau
=
\sigma^{(1)}
\shortparallel 
\sigma^{(2)}
;
\sigma^{(1)},\sigma^{(2)}
\in
\Sigma_N
\right\}
\end{align*}
is actually an instance of the $2d$-dimensional Gaussian process defined in
\eqref{eq:as2:sk-with-multidimensional-spins-process}. Hence, it has the
following correlation structure, for $\tau^1,\tau^2 \in \Sigma^{(2)}_N$,
\begin{align*}
\cov
\left[
X^{(2)}(\tau^1)
,
X^{(2)}(\tau^2)
\right]
=
\Vert
R^{(2)}
(
\tau^1
,
\tau^2
)
\Vert_{\text{F}}^2
.
\end{align*}
The path $\varrho$ induces also the following two new
(independent of everything before) comparison process
$
Y^{(2)} 
\equiv 
\left\{ 
Y^{(2)}(\alpha)
\in 
\R^{2d} 
\mid
\alpha
\in
\mathcal{A}_{\mathfrak{n}}
\right\}
$,
with the following correlation structures
\begin{align*}
\cov
\left[
Y^{(2)}(\alpha^{(1)})
,
Y^{(2)}(\alpha^{(2)})
\right]
=
\mathfrak{Q}(\alpha^{(1)},\alpha^{(2)})
\in
\symmetric^+(d)
.
\end{align*}
As usual, let $\{Y^{(2)}_i\}_{i=1}^N$ be the independent copies of $Y^{(2)}$. 
For the purposes of new Guerra's scheme we define a GREM-like
process (cf. \eqref{eq:as2:a-process-definition})
\begin{align*}
A^{(2)}
=
\{
A^{(2)}
(\tau,\alpha)
:
\tau
=
\sigma^{(1)}
\shortparallel
\sigma^{(2)}
;
\sigma^{(1)}
,
\sigma^{(2)} 
\in 
\Sigma_N 
;
\alpha
\in 
\mathcal{A}_{
\mathfrak{n}
}
\}
\end{align*}
as
\begin{align*}
A^{(2)}
(\tau,\alpha)
\equiv
\left(
\frac{2}{N}
\right)^{1/2}
\smash{
\sum_{i=1}^N
}
\vphantom{\sum}
\langle
Y^{(2)}_i(\alpha)
,
\tau_i
\rangle
.
\end{align*}
We fix some $t \in [0;1]$.
We would now like to apply Guerra's scheme to the comparison functional
\eqref{eq:as2-phi-2-k-x} and the following two processes
\begin{align*}
\left\{
H^{(2)}_t(\sigma^{(1)},\sigma^{(2)},\alpha)
\right\}_{
\sigma^{(1)},\sigma^{(2)} \in \Sigma_N
,
\alpha
\in
\mathcal{A}
}
,
\left\{
\sqrt{t}
A^{(2)}(\sigma^{(1)} \shortparallel \sigma^{(2)},\alpha)
\right\}_{
\sigma^{(1)},\sigma^{(2)} \in \Sigma_N
,
\alpha
\in
\mathcal{A}
}
.
\end{align*}
These two processes are, respectively, the counterparts of the processes
$X(\sigma)$ and $A(\sigma,\alpha)$ in Guerra's scheme.

Consider a path
$
\widetilde{Q}
\in
\mathcal{Q}^\prime(U,d)
$ 
with the following jumps
\begin{align*}
0
=:
\widetilde{Q}^{(0)}
\prec
\widetilde{Q}^{(1)}
\prec
\ldots
\prec
\widetilde{Q}^{(\mathfrak{n})}
\prec
\widetilde{Q}^{(\mathfrak{n}+1)}
.
\end{align*}
Let 
$
\widetilde{A}
\equiv
\left\{ 
\widetilde{A}(\sigma,\alpha)
:
\sigma
\in
\Sigma_N
;
\alpha
\in
\mathcal{A}_{
\mathfrak{n}
}
\right\}
$
be a Gaussian process (independent of all random objects around) with the
following covariance structure
\begin{align*}
\E
\left[
\widetilde{A}(\sigma^{(1)},\alpha^{(1)})
\widetilde{A}(\sigma^{(2)},\alpha^{(2)})
\right]
=
2
\langle
R(\sigma^{(1)},\sigma^{(2)})
,
\widetilde{Q}(\alpha^{(1)},\alpha^{(2)})
\rangle
.
\end{align*}
For notational convenience, we introduce also the
following process
\begin{align}
\label{eq:remainder:replicated-comparison-process-a-2-tilde}
\widetilde{A}^{(2)}
(
\sigma^{(1)}
\shortparallel
\sigma^{(2)}
,
\alpha^{(1)},
\alpha^{(2)}
)
\equiv
\widetilde{A}(\sigma^{(1)},\alpha^{(1)})
+
\widetilde{A}(\sigma^{(2)},\alpha^{(2)})
.
\end{align}
Recalling the replicated Hamiltonian
\eqref{eq:as2-double-replicated-rest-term-hamiltonian} and following 
Guerra's scheme, we introduce, for $s \in [0;1]$, the following interpolating
Hamiltonian
\begin{align}
\label{eq:as2-doubled-ass-hamiltonian}
H^{(2)}_{t,s}(\sigma^{(1)},\sigma^{(2)},\alpha^{(1)},\alpha^{(2)})
\equiv
&
\sqrt{st}
X^{(2)}(
\sigma^{(1)}
\shortparallel 
\sigma^{(2)}
)
+
\sqrt{(1-s)t}
A^{(2)}(\sigma^{(1)} \shortparallel \sigma^{(2)},\alpha^{(1)})
\nonumber
\\
&
+
\sqrt{1-t}
\widetilde{A}^{(2)}
(
\sigma^{(1)}
\shortparallel
\sigma^{(2)},
\alpha^{(1)},
\alpha^{(2)}
)
.
\end{align}
Given $\eps, \delta > 0$ and $\mathfrak{L} \in \symmetric(2d)$, define (cf.
\eqref{eq:remainder:small-vicinity-of-u})
\begin{align*}
\mathcal{V}^{(2)}(\mathfrak{L},\mathfrak{U},\eps,\delta)
\equiv
\{
\mathfrak{U}^\prime
\in
\symmetric^+(2d)
:
\Vert
\mathfrak{U}^\prime
-
\mathfrak{U}
\Vert_\text{F}
<
\eps
,
\langle
\mathfrak{U}^\prime
-
\mathfrak{U}
,
\mathfrak{L}
\rangle
<
\delta
\}
.
\end{align*}
We consider the following set of the local configurations
\begin{align}
\label{eq:as2-sigma-2-n-v-u-delta-eps}
\Sigma^{(2)}_N(\mathfrak{L},\mathfrak{U},\eps,\delta)
\equiv
\left\{
(\sigma^{(1)},\sigma^{(2)})
\in
\Sigma_N \times \Sigma_N
:
R^{(2)}_N(
\sigma^{(1)}\shortparallel\sigma^{(2)}
,
\sigma^{(1)}\shortparallel\sigma^{(2)}
)
\in
\mathcal{V}^{(2)}(\mathfrak{L},\mathfrak{U},\eps,\delta)
\right\}
.
\end{align}
Note that 
$
\Sigma^{(2)}_N(\mathfrak{L},\mathfrak{U},\eps,\delta) 
\subset
\Sigma_N(B(U,\eps))^2
$. 
We consider also the RPC $\zeta = \zeta(\mathfrak{x})$
generated by the vector $\mathfrak{x}$ and, for any suitable Gaussian process
\begin{align*}
F
\equiv
\{
F(\sigma^{(1)},\sigma^{(2)},\alpha^{(1)},\alpha^{(2)})
\mid
\sigma^{(1)},\sigma^{(2)}
\in
\Sigma_N
;
\alpha^{(1)},\alpha^{(2)}
\in
\mathcal{A}_n
\}
,
\end{align*}
define the corresponding local comparison functional (cf.
\eqref{eq:as2-phi-2-k-x}) as follows 
\begin{align*}
\Phi^{(2),k,\mathfrak{x}}_{
\mathcal{V}
}
\left[
F
\right]
\equiv
\frac{1}{N}
\E
\Bigl[
\log
&
\iint_{
\mathcal{V}
}
\iint_{
\mathcal{A}^{(2),k}
}
\exp
\left\{
\beta
\sqrt{N}
F(\sigma^{(1)},\sigma^{(2)},\alpha^{(1)},\alpha^{(2)})
\right\}
\\
&
\dd 
\mu^{\otimes N}(\sigma^{(1)})
\dd 
\mu^{\otimes N}(\sigma^{(2)})
\dd \zeta(\alpha^{(1)})
\dd \zeta(\alpha^{(2)})
\Bigr]
.
\end{align*}
Define the corresponding local free energy-like quantity as
(cf. \eqref{eq:as2-phi-of-t-definition})
\begin{align}
\label{eq:remainder:localised-free-energy-like-quantity}
\chi(s,t,k,\mathfrak{x},\mathfrak{Q},\widetilde{\mathfrak{Q}},\Sigma^{(2)}_N(\mathfrak{L},\mathfrak{U},\eps,\delta))
\equiv
\Phi^{(2),k,\mathfrak{x}}_{
\Sigma^{(2)}_N(\mathfrak{L},\mathfrak{U},\eps,\delta)
}
\left[
H^{(2)}_{t,s}
\right]
.
\end{align}
To lighten the notation, we indicate hereinafter only the dependence of $\chi$
on $s$.
Denote
\begin{align*}
B^{\mathfrak{x},\mathfrak{Q}}
\equiv
\frac{t \beta^2}{2}
\sum_{
l=1
}^{
\mathfrak{n}
}
\mathfrak{x}_l
\left(
\Vert
\mathfrak{Q}^{(l+1)}
\Vert_{\text{F}}^2
-
\Vert
\mathfrak{Q}^{(l)}
\Vert_{\text{F}}^2
\right)
.
\end{align*}
\begin{lemma}
There exists $C = C(\Sigma) > 0$ such that, for any $\mathfrak{U}$ as above, we
have
\begin{align}
\label{eq:as2-chi-of-s-derivative}
\frac{
\partial
}{
\partial s
}
\chi(s,t,k,\mathfrak{x},\mathfrak{Q},\widetilde{\mathfrak{Q}},\Sigma^{(2)}_N(\mathfrak{L},\mathfrak{U},\eps,\delta))
\leq
-
B^{\mathfrak{x},\mathfrak{Q}}
+
C\eps
,
\end{align}
Consequently,
\begin{align}
\varphi^{(2)}_{N}
(k,t,x,Q,\Sigma^{(2)}_N(\mathfrak{L},\mathfrak{U},\eps,\delta))
\leq
&
\Phi^{(2),k,\mathfrak{x}}_{
\Sigma^{(2)}_N(\mathfrak{L},\mathfrak{U},\eps,\delta)
}
\Bigl[
\sqrt{t}
A^{(2)}(\sigma^{(1)} \shortparallel \sigma^{(2)},\alpha^{(1)})
\nonumber
\\
&
+
\sqrt{1-t}
\widetilde{A}^{(2)}
(
\sigma^{(1)}
\shortparallel
\sigma^{(2)}
,
\alpha^{(1)},
\alpha^{(2)}
)
\Bigr]
-
B^{\mathfrak{x},\mathfrak{Q}}
+
C\eps
.
\end{align}
\end{lemma}

\begin{proof}
The idea is the same as in the proof of Theorem
\ref{thm:as2:guerras-interpolation} and is based on
Proposition~\ref{prp:as2:gaussian-comparison-of-free-energy}. Since we are
considering the localised free energy-like quantities
\eqref{eq:remainder:localised-free-energy-like-quantity}, the variance terms
induced by the interpolation \eqref{eq:as2-doubled-ass-hamiltonian} in
\eqref{eq:as2-scheme-comparison-formula} cancel out (up to the correction $\mathcal{O}(\eps)$) and we are left with the non-positive
contribution of the covariance terms.
\end{proof}
Given $\mathfrak{L} \in \symmetric(2d)$, we consider the following stencil of the Legendre transform
\begin{align}
\label{eq:remainder:replicated-legendre-transofrorm}
\widetilde{\Phi}^{(2),k,\mathfrak{x},\mathfrak{L}}
\left[
F
\right]
\equiv
-
\langle
\mathfrak{L}
,
\mathfrak{U}
\rangle
-
B^{\mathfrak{x},\mathfrak{Q}}
+
\frac{1}{N}
\E
[
\log
&
\iint_{
\Sigma_N^2
}
\iint_{
\mathcal{A}^{(2),k}
}
\exp
\{
\beta
\sqrt{N}
F(\sigma^{(1)},\sigma^{(2)},\alpha^{(1)},\alpha^{(2)})
\nonumber
\\
&
\quad
+
\langle
\mathfrak{L}
(
\sigma^{(1)}
\shortparallel
\sigma^{(2)}
)
,
\sigma^{(1)}
\shortparallel
\sigma^{(2)}
\rangle
\}
\nonumber
\\
&
\left.
\dd 
\mu^{\otimes N}(\sigma^{(1)})
\dd 
\mu^{\otimes N}(\sigma^{(2)})
\dd \zeta(\alpha^{(1)})
\dd \zeta(\alpha^{(2)})
\right]
.
\end{align}
\begin{definition}
Let $F: \symmetric(2d) \to \R$.
Given $\delta > 0$, we call 
$
\mathfrak{L}^{(0)} \in \symmetric(2d)
$
\emph{$\delta$-minimal
for $F$}, if
\begin{align*}
F(\Lambda^{(0)})
\leq
\inf_{
\Lambda \in \symmetric(2d)
}
F(\Lambda)
+
\delta
.
\end{align*}
\end{definition}
\begin{lemma}
\label{lem:as2:rest-term-ldp-upper-bound}
There exists $C = C(\Sigma) > 0$ such that,
for all $\mathfrak{U}$ and 
$
\mathfrak{Q} \in \mathcal{Q}^\prime(\mathfrak{U},2d)
$ 
as
above, all $\eps, \delta > 0$, there exists a $\delta$-minimal Lagrange multiplier
$
\mathfrak{L} 
=
\mathfrak{L}(\mathfrak{U},\eps,\delta) 
\in 
\symmetric(2d)
$
for \eqref{eq:remainder:replicated-legendre-transofrorm}
such that, for all $k \in [1;n] \cap \N$,
all $t \in [0;1]$, 
and all $(x,\mathcal{Q})$,
we have
\begin{align}
\varphi^{(2)}_{N}
(k,t,x,Q,\Sigma^{(2)}_N(\mathfrak{L},\mathfrak{U},\eps,\delta))
\leq
&
\inf_{
\mathfrak{L} \in \symmetric(2d)
}
\widetilde{\Phi}^{(2),k,\mathfrak{x},\mathfrak{L}}
\left[
\sqrt{t}
A^{(2)}(\sigma^{(1)} \shortparallel \sigma^{(2)},\alpha^{(1)})
\nonumber
\right.
\\
&
\quad
+
\left.
\sqrt{1-t}
\widetilde{A}^{(2)}
(
\sigma^{(1)}
\shortparallel
\sigma^{(2)}
,
\alpha^{(1)},
\alpha^{(2)}
)
\right]
\nonumber
\\
&
+
C
(
\eps + \delta
)
.
\end{align}
\end{lemma}
\begin{proof}
The argument is the same as in the proof of
Theorem~\ref{thm:as2:pressure-upper-bound}.
\end{proof}
Consider the family of matrices 
$
\widetilde{\mathfrak{Q}}
\equiv
\left\{
\widetilde{\mathfrak{Q}}^{(l)}
\in 
\symmetric^+(2d)
\mid
l \in [0;\mathfrak{n}+1] \cap \N
\right\}
$, 
defined as
\begin{align}
\label{eq:as2:gaussian-spins:tilide-frak-q-definition-1}
\widetilde{\mathfrak{Q}}^{(l)}
\equiv
\begin{bmatrix}
\widetilde{Q}^{(l)}
&
\widetilde{Q}^{(l)}
\\
\widetilde{Q}^{(l)}
&
\widetilde{Q}^{(l)}
\end{bmatrix}
,
\end{align}
for $l \in [0;k] \cap \N$,
and as
\begin{align}
\label{eq:as2:gaussian-spins:tilide-frak-q-definition-2}
\widetilde{\mathfrak{Q}}^{(l)}
\equiv
\begin{bmatrix}
\widetilde{Q}^{(l)}
&
\widetilde{Q}^{(k)}
\\
\widetilde{Q}^{(k)}
&
\widetilde{Q}^{(l)}
\end{bmatrix}
,
\end{align}
for $l \in [k+1;\mathfrak{n}+1] \cap \N$.
Additionally we define, for $l \in [0;\mathfrak{n}+1]$, the matrices
\begin{align*}
\widehat{\mathfrak{Q}}^{(l)}(t)
\equiv
t \mathfrak{Q}
+
(1-t) \widetilde{\mathfrak{Q}}
.
\end{align*}
Let $\widehat{Z}^{(l)} \in \R^{2d \times 2d}$, for $l \in [0;\mathfrak{x}]$, be
independent Gaussian vectors with
\begin{align*}
\cov
\left[
\widehat{Z}^{(l)}
\right]
=
2
\beta^2
\Bigl(
\widehat{\mathfrak{Q}}^{(l+1)}(t)
-
\widehat{\mathfrak{Q}}^{(l)}(t)
\Bigr)
.
\end{align*}
Given $\widehat{y} \in \R^{2d}$, $\mathfrak{L} \in \symmetric(2d)$, consider the random variable
\begin{align}
\label{eq:as2:doubled-start-of-the-recursion}
X^{(2)}_{\mathfrak{n}+1}
(
\widehat{y}
,
\mathfrak{x}
,
\widehat{\mathfrak{Q}}(t)
,
\mathfrak{L}
)
\equiv
\log
\int_{
\Sigma
}
\int_{
\Sigma
}
\exp
\Bigl(
&
\langle
\widehat{y}
,
\sigma^{(1)}
\shortparallel
\sigma^{(2)}
\rangle
+
\langle
\mathfrak{L}
(
\sigma^{(1)}
\shortparallel
\sigma^{(2)}
)
,
\sigma^{(1)}
\shortparallel
\sigma^{(2)}
\rangle
\Bigr)
\dd \mu(\sigma^{(1)})
\dd \mu(\sigma^{(2)})
.
\end{align}
Define recursively, for $l \in [\mathfrak{n};0] \cap \N$, the following
quantities
\begin{align}
\label{eq:remainder:doubled-recursive-x-term}
X^{(2)}_{l}
(
\widehat{y}
,
k
,
\mathfrak{x}
,
\widehat{\mathfrak{Q}}(t)
,
\mathfrak{L}
)
\equiv
\frac{1}{
\mathfrak{x}_l
}
\log
\E^{
\widehat{Z}^{(l)}
}
\left[
\exp
\left(
\mathfrak{x}_l
X^{(2)}_{l+1}
(
\widehat{y}
+
\widehat{Z}^{(l)}
,
k
,
\mathfrak{x}
,
\widehat{\mathfrak{Q}}^{(l)}(t)
,
\mathfrak{L}
)
\right)
\right]
.
\end{align}
\begin{lemma}
\label{lem:remainder:the-computation-of-the-doubled-a-term}
We have
\begin{align}
\label{eq:remainder:a-formula-for-the-doubled-a-term}
&
\widetilde{\Phi}^{(2),k,\mathfrak{x},\mathfrak{L}}
\left[
\sqrt{t}
A^{(2)}(\sigma^{(1)} \shortparallel \sigma^{(2)},\alpha^{(1)})
\nonumber
+
\sqrt{1-t}
\widetilde{A}^{(2)}
(
\sigma^{(1)}
\shortparallel
\sigma^{(2)}
,
\alpha^{(1)},
\alpha^{(2)}
)
\right]
\\
&
=
-
\langle
\mathfrak{L}
,
\mathfrak{U}
\rangle
+
X^{(2)}_{0}
(
0
,
\mathfrak{x}
,
\widehat{\mathfrak{Q}}^{(l)}(t)
,
\mathfrak{L}
)
.
\end{align}
\end{lemma}
\begin{proof}
This is an immediate consequence of the RPC averaging property 
\eqref{eq:remainder:averaging-property}.
\end{proof}
\begin{proposition}
\label{eq:remainder:phi-2-upper-bound}
Under the conditions of Lemma~\ref{lem:as2:rest-term-ldp-upper-bound}, we have
\begin{align*}
\varphi^{(2)}_{N}
(k,t,x,Q,\Sigma^{(2)}_N(\mathfrak{L},\mathfrak{U},\eps,\delta))
\leq
\inf_{
\mathfrak{L}
\in
\symmetric(2d)
}
\Bigl(
-
\langle
\mathfrak{L}
,
\mathfrak{U}
\rangle
+
X^{(2)}_{0}
(
0
,
\mathfrak{x}
,
\widehat{\mathfrak{Q}}(t)
,
\mathfrak{L}
)
\Bigr)
-
B^{\mathfrak{x},\mathfrak{Q}}
+
C
(
\eps + \delta
)
.
\end{align*}
\end{proposition}
\begin{remark}
Similarly to
\eqref{eq:remainder:local-free-energy-lambda-dependent-upper-bound}, 
there exists $C = C(\Sigma,\mu)>0$,
such that, for any $\mathfrak{L} \in \symmetric(2d)$,
\begin{align*}
\varphi^{(2)}_{N}
(k,t,x,Q,\Sigma^{(2)}_N(B(\mathfrak{U},\eps))
\leq
-
\langle
\mathfrak{L}
,
\mathfrak{U}
\rangle
-
B^{\mathfrak{x},\mathfrak{Q}}
+
X^{(2)}_{0}
(
0
,
\mathfrak{x}
,
\widehat{\mathfrak{Q}}(t)
,
\mathfrak{L}
)
\Bigr)
+
C
\Vert
\mathfrak{L}
\Vert_{\text{F}}
\eps
.
\end{align*}
\end{remark}
\begin{proof}
Immediately follows from Lemmata~\ref{lem:as2:rest-term-ldp-upper-bound} and
\ref{lem:remainder:the-computation-of-the-doubled-a-term}.
\end{proof}

\subsection{Adjustment of the upper bounds on $\varphi^{(2)}$}
\label{sec:remainder:adjustment-of-the-upper-bounds-on-varphi}
Proposition~\ref{thm:lecture-01-matrix-overlap-properties} implies that there
exists $r \in [1;n] \cap \N$ such that
\begin{align}
\label{eq:remainder:v-between-q}
\Vert
Q^{(r-1)}
\Vert_{\text{F}}^2
<
\Vert
V
\Vert_{\text{F}}^2
<
\Vert
Q^{(r)}
\Vert_{\text{F}}^2
.
\end{align}
Assume $r = k$. (Other cases are similar or easier as shown for 1-D in
\cite{TalagrandParisiFormula2006}.) We make the following tuning of the upper
bounds of the previous subsection. Set $\mathfrak{n} \equiv n+1$. Let $
w \in [x_{r-1}/2;x_r]
$.
Define
\begin{align}
\label{eq:remainder:mathfrak-x-case-r-equals-k}
\mathfrak{x}_l
\equiv
\mathfrak{x}_l(w)
\equiv
\begin{cases}
\frac{x_l}{2}
,
&
l \in [0; k-1] \cap \N
,
\\
w
,
&
l = k
,
\\
x_l
,
&
l \in [k+1; n+1] \cap \N
.
\end{cases}
\end{align}
Let
\begin{align*}
\widetilde{Q}^{(l)}
\equiv
\begin{cases}
Q^{(l)} 
,
&
l \in [0;k-1] \cap \N
,
\\
Q^{(l-1)} 
,
&
l \in [k; n+2] \cap \N
.
\end{cases}
\end{align*}
Moreover, suppose
$
\mathfrak{Q}
\equiv
\{
\mathfrak{Q}^{(l)}
\}_{l=0}^{n+2}
$
satisfy
\begin{align}
\label{eq:remainder:mathfrak-q-case-r-equals-k}
\Vert
\mathfrak{Q}^{(l)}
\Vert_{\text{F}}^2
=
\begin{cases}
4
\Vert
Q^{(l)}
\Vert_{\text{F}}^2
,
&
l \in [0;k-1] \cap \N
,
\\
4
\Vert
V
\Vert_{\text{F}}^2
,
&
l = k
,
\\
2
\Bigl(
\Vert
Q^{(l-1)}
\Vert_{\text{F}}^2
+
\Vert
V
\Vert_{\text{F}}^2
\Bigr)
,
&
l \in [k+1;n+2] \cap \N
.
\end{cases}
\end{align}
Such $\mathfrak{Q}$ exists due to \eqref{eq:remainder:v-between-q}. Moreover,
if $d \geq 2$, then it is obviously non-unique.
\begin{lemma}
\label{lem:remainder:replicated-b-term-in-terms-of-the-non-replicated-one}
In the above setup, we have
\begin{align*}
B^{\mathfrak{x},\mathfrak{Q}}
\equiv
t \beta^2
\Bigl\{
(w-x_{l-1})
\left(
\Vert
Q^{(k)}
\Vert_{\text{F}}^2
-
\Vert
V
\Vert_{\text{F}}^2
\right)
+
\sum_{
l=1
}^{
{n}
}
x_l
\left(
\Vert
Q^{(l+1)}
\Vert_{\text{F}}^2
-
\Vert
Q^{(l)}
\Vert_{\text{F}}^2
\right)
\Bigr\}
.
\end{align*}
\end{lemma}
\begin{proof}
The claim is a straightforward consequence of
\eqref{eq:remainder:mathfrak-x-case-r-equals-k} and \eqref{eq:remainder:mathfrak-q-case-r-equals-k}.
\end{proof}
Define the matrix
$
\mathfrak{D}^{(n+1)}
\in
\symmetric^+(2d)
$
block-wise as
\begin{align*}
\mathfrak{D}^{(n+1)}\vert_{11}
&
\equiv
\beta^2 t (U - \mathfrak{Q}^{(n+1)}\vert_{11})
+
\beta^2  (1 - t) (U - Q^{(n)})
+
\mathfrak{L}\vert_{11}
,
\\
\mathfrak{D}^{(n+1)}\vert_{12}
&
\equiv
\beta^2 t (V - \mathfrak{Q}^{(n+1)}\vert_{12})
+
\mathfrak{L}\vert_{12}
,
\\
\mathfrak{D}^{(n+1)}\vert_{21}
&
\equiv
\beta^2  t (V - \mathfrak{Q}^{(n+1)}\vert_{12})^*
+
\mathfrak{L}\vert_{12}^*
,
\\
\mathfrak{D}^{(n+1)}\vert_{22}
&
\equiv
\beta^2 t (U - \mathfrak{Q}^{(n+1)}\vert_{22})
+
\beta^2 (1 - t) (U - Q^{(n)})
+
\mathfrak{L}\vert_{22}
.
\end{align*}
Furthermore, we define
\begin{align*}
\symmetric^+(2d)
\ni
\widetilde{\mathfrak{D}}^{(n+1)}
\equiv
\begin{bmatrix}
\beta^2  (U-Q^{(n)}) + \Lambda
&
0
\\
0
&
\beta^2  (U-Q^{(n)}) + \Lambda
\end{bmatrix}
.
\end{align*}
\begin{lemma}
We have
\begin{align*}
X^{(2)}_{n+1}
(
\widehat{y}
,
k
,
\mathfrak{x}
,
\widehat{\mathfrak{Q}}(t)
,
\mathfrak{L}
)
\equiv
\log
\int_{
\Sigma
}
\int_{
\Sigma
}
&
\exp
\Bigl(
\langle
\widehat{y}
,
\sigma^{(1)}
\shortparallel
\sigma^{(2)}
\rangle
+
\langle
\mathfrak{D}^{(n+1)}
(
\sigma^{(1)}
\shortparallel
\sigma^{(2)}
)
,
\sigma^{(1)}
\shortparallel
\sigma^{(2)}
\rangle
\Bigr)
\\
&
\times
\dd \mu(\sigma^{(1)})
\dd \mu(\sigma^{(2)})
.
\end{align*}
\end{lemma}
\begin{proof}
Since $\mathfrak{x}_{n+2} = 1$, the result follows from a straightforward
calculation of the Gaussian integrals in \eqref{eq:remainder:doubled-recursive-x-term} for $l=n+1$.
\end{proof}
Define
\begin{align*}
\widetilde{\mathfrak{L}}
\equiv
\begin{bmatrix}
\Lambda
&
0
\\
0
&
\Lambda
\end{bmatrix}
,
\widetilde{\mathfrak{U}}
\equiv
\begin{bmatrix}
U
&
Q^{(k)}
\\
Q^{(k)}
&
U
\end{bmatrix}
.
\end{align*}
\begin{lemma}
\label{lem:remainder:doubled-A-term-of-tilde-mathfrak-Q}
For any 
$
y \in \R^d
$,
$
l \in [0;n+2] \cap \N
$,
we have
\begin{align*}
X^{(2)}_{l}
(
y
\shortparallel
y
,
\mathfrak{x}(w)
,
\widetilde{\mathfrak{Q}}
,
\widetilde{\mathfrak{L}}
)
\vert_{
w = x_{k-1}
}
=
\begin{cases}
2 X_{l-1}(y, x, \mathcal{Q}, U, \Lambda)
,
&
l \in [k;n+2] \cap \N
,
\\
2 X_{l}(y, x, \mathcal{Q}, U, \Lambda)
,
&
l \in [0;k-1] \cap \N
.
\end{cases}
\end{align*}
\end{lemma}
\begin{proof}
A straightforward (decreasing) induction argument on $l$ gives the result.
Indeed: for $l = n+2$, an inspection of
\eqref{eq:as2:doubled-start-of-the-recursion} and
\eqref{eq:as2:x-end-condition} immediately yields
\begin{align*}
X^{(2)}_{n+2}
(
y^{(1)}
\shortparallel
y^{(2)}
,
\mathfrak{x}(w)
,
\widetilde{\mathfrak{Q}}
,
\widetilde{\mathfrak{L}}
)
=
X_{n+1}(y^{(1)}, x, \mathcal{Q}, U, \Lambda)
+
X_{n+1}(y^{(2)}, x, \mathcal{Q}, U, \Lambda)
,
\end{align*}
where $y^{(1)}, y^{(2)} \in \R^d$.
Let 
$
\widehat{Z}^{(l)} 
$ 
be a Gaussian $2d$-dimensional vector with 
\begin{align*}
\cov
\Bigl[
\widetilde{Z}^{(l)} 
\Bigr]
=
2 \beta^2
(
\widetilde{\mathfrak{Q}}^{(l+1)}
-
\widetilde{\mathfrak{Q}}^{(l)}
)
. 
\end{align*}
Define two Gaussian $d$-dimensional vectors $\widetilde{Z}^{(l),1}$ and
$\widetilde{Z}^{(l),2}$ by demanding that 
\begin{align*}
\widetilde{Z}^{(l)} 
=
\widetilde{Z}^{(l),1}
\shortparallel
\widetilde{Z}^{(l),2}
.
\end{align*}
Due to \eqref{eq:as2:gaussian-spins:tilide-frak-q-definition-1} and
\eqref{eq:as2:gaussian-spins:tilide-frak-q-definition-2}, the vectors
$\widetilde{Z}^{(l),1}$ and $\widetilde{Z}^{(l),2}$ are independent, for $l \in [k; n+1]$. We have $
\widetilde{Z}^{(l),1}
\sim
\widetilde{Z}^{(l),2}
$,  for $l \in [0; k-1]$.
Assume that 
$
l \in [k;n+1] \cap \N
$ 
and 
\begin{align*}
X^{(2)}_{l+1}
(
y^{(1)}
\shortparallel
y^{(2)}
,
\mathfrak{x}(w)
,
\widetilde{\mathfrak{Q}}
,
\widetilde{\mathfrak{L}}
)
=
X_{l}(y^{(1)}, x, \mathcal{Q}, U, \Lambda)
+
X_{l}(y^{(2)}, x, \mathcal{Q}, U, \Lambda)
.
\end{align*}
By definition \eqref{eq:as2:doubled-start-of-the-recursion}, we have
\begin{align*}
X^{(2)}_{l}
(
y^{(1)}
\shortparallel
y^{(2)}
,
k
,
\mathfrak{x}
,
\widetilde{\mathfrak{Q}}
,
\widetilde{\mathfrak{L}}
)
&
=
\frac{1}{
\mathfrak{x}_l
}
\log
\E^{
\widetilde{Z}^{(l)}
}
\left[
\exp
\left(
\mathfrak{x}_l
X^{(2)}_{l+1}
(
y^{(1)}
\shortparallel
y^{(2)}
+
\widehat{Z}^{(l)}
,
k
,
\mathfrak{x}
,
\widetilde{\mathfrak{Q}}
,
\mathfrak{L}
)
\right)
\right]
\\
&
=
\frac{1}{
x_l
}
\log
\E^{
\widetilde{Z}^{(l)}
}
\Bigl[
\exp
\Bigl\{
x_l
\Bigl(
X_{l}(y^{(1)}+\widetilde{Z}^{(l),1}, x, \mathcal{Q}, U, \Lambda)
\\
&
\quad
+
X_{l}(y^{(2)}+\widetilde{Z}^{(l),2}, x, \mathcal{Q}, U, \Lambda)
\Bigr)
\Bigr\}
\Bigr]
\\
&
=
X_{l-1}(y^{(1)}, x, \mathcal{Q}, U, \Lambda)
+
X_{l-1}(y^{(2)}, x, \mathcal{Q}, U, \Lambda)
.
\end{align*}
By the construction and previous formula, for $l = k-1$, we have
\begin{align*}
X^{(2)}_{k-1}
(
y^{(1)}
\shortparallel
y^{(2)}
,
k
,
\mathfrak{x}
,
\widetilde{\mathfrak{Q}}
,
\widetilde{\mathfrak{L}}
)
\vert_{
w = x_{k-1}
}
&
=
X^{(2)}_{k}
(
y^{(1)}
\shortparallel
y^{(2)}
,
k
,
\mathfrak{x}
,
\widetilde{\mathfrak{Q}}
,
\widetilde{\mathfrak{L}}
)
\\
&
=
X_{k-1}(y^{(1)}, x, \mathcal{Q}, U, \Lambda)
+
X_{k-1}(y^{(2)}, x, \mathcal{Q}, U, \Lambda)
.
\end{align*}
Finally, for $l \in [0;k-2]$, we recursively obtain
\begin{align*}
X^{(2)}_{l}
(
y^{(1)}
\shortparallel
y^{(1)}
,
k
,
\mathfrak{x}
,
\widetilde{\mathfrak{Q}}
,
\widetilde{\mathfrak{L}}
)
\vert_{
w = x_{k-1}
}
&
=
\frac{1}{
\mathfrak{x}_l
}
\log
\E^{
\widetilde{Z}^{(l)}
}
\left[
\exp
\left(
\mathfrak{x}_l
X^{(2)}_{l+1}
(
y^{(1)}
\shortparallel
y^{(1)}
+
\widehat{Z}^{(l)}
,
k
,
\mathfrak{x}
,
\widetilde{\mathfrak{Q}}
,
\mathfrak{L}
)\vert_{
w = x_{k-1}
}
\right)
\right]
\\
&
=
\frac{2}{
x_l
}
\log
\E^{
\widetilde{Z}^{(l),1}
}
\Bigl[
\exp
\Bigl\{
\frac{
x_l
}{
2
}
\Bigl(
X_{l+1}(y^{(1)}+\widetilde{Z}^{(l),1}, x, \mathcal{Q}, U, \Lambda)
\\
&
\quad
+
X_{l+1}(y^{(1)}+\widetilde{Z}^{(l),1}, x, \mathcal{Q}, U, \Lambda)
\Bigr)
\Bigr\}
\Bigr]
\\
&
=
2 X_{l}(y^{(1)}, x, \mathcal{Q}, U, \Lambda)
.
\end{align*}
\end{proof}
\begin{remark}
\label{rem:remainder:the-main-thechnical-problem-of-talagrands-approach}
Motivated by 
Lemmata~\ref{thm:as2-a-sufficient-condition-for-mu-k-to-vanish}
and \ref{lem:remainder:doubled-A-term-of-tilde-mathfrak-Q} (see also
Section~\ref{sec:a-priori-estimates}), we pose the following problem. Is it
true that, as in 1-D (see
\cite{TalagrandParisiFormula2006,panchenko-free-energy-generalized-sk-2005}),
there exists $ \mathfrak{Q} \in \mathcal{Q}^\prime(\mathfrak{U},2d) $ satisfying the assumption
\eqref{eq:remainder:mathfrak-q-case-r-equals-k} such that the following
inequality holds
\begin{align}
\label{eq:remainder:the-doubling-of-the-free-energy}
&
\inf_{
\mathfrak{L}
\in
\symmetric(2d)
}
\Bigl(
-
\langle
\mathfrak{L}
,
\mathfrak{U}
\rangle
+
X^{(2)}_{0}
(
0
,
\mathfrak{x}(w)
,
\widehat{\mathfrak{Q}}(t)
,
\mathfrak{L}
)
\vert_{
w = x_{k-1}
}
\Bigr)
\nonumber
\\
&
\overset{?}{\leq}
2 
\inf_{
\Lambda
\in
\symmetric(d)
}
\Bigl(
-
\langle
\Lambda
,
U
\rangle
+
X_{0}(0, x, \mathcal{Q}, U, \Lambda)
\Bigr)
?
\end{align}
Similar problems have at first been posed in \cite{Talagrand2007}.
The resolution of the above problem seems to require more detailed information
on the behaviour of the Parisi functional \eqref{vektor-sk-pasisis-functional}
or, equivalently, of the solution of \eqref{eq:remainder:generalised-parisi-pde}
as a function of $
Q \in \mathcal{Q}(U,d)$
. 
\end{remark}

\subsection{Talagrand's a priori estimates}
\label{sec:a-priori-estimates}
We start from defining a class of the almost optimal paths for the
optimisation problem in \eqref{eq:multidim-sk:upper-bound-pde}. Recall the
following convenient definition
from \cite{panchenko-free-energy-generalized-sk-2005}.
\begin{definition}
\label{def:as2-eps-optimizer}
Given $U \in \symmetric^+(d)$, we shall call the triple
$
(n,\rho^*,\Lambda^*) 
\in
\N 
\times 
\mathcal{Q}^\prime_n(U, d)
\times 
\R^d
$
a \emph{$\theta$-optimiser} of the Parisi functional
\eqref{vektor-sk-pasisis-functional}, if it satisfies the following
two conditions
\begin{align}
\label{eq:remainder:almost-optimality-for-all-n}
\mathcal{P}(\beta, \rho^*, \Lambda^*)
\leq
\inf_{
\substack{
\rho
\in
\mathcal{Q}^\prime(U, d)
\\
\Lambda
\in
\symmetric(d)
}
} 
\mathcal{P}(\beta, \rho, \Lambda)
+
\theta
.
\\
\label{eq:remainder:optimality-for-fixed-n}
\mathcal{P}(\beta, \rho^*,\Lambda^*)
=
\inf_{
\substack{
\rho
\in
\mathcal{Q}^\prime_n(U,d)
\\
\Lambda
\in
\symmetric(d)
}
} 
\mathcal{P}(\beta, \rho, \Lambda)
.
\end{align}
\end{definition}
\begin{remark}
It is obvious that for any $\theta>0$ such a $\theta$-optimiser exists. The main
convenient feature of this definition (as pointed out in
\cite{TalagrandParisiFormula2006}) is that $n$ (the number of jumps of $\rho^*$)
is finite and fixed.
\end{remark}
Recalling \eqref{eq:remainder:guerras-sum-rule-ldp-upper-bound}, we set
\begin{align}
\label{eq:as2-psi-of-t-definition}
\phi^{(x,\mathcal{Q},\Lambda)}(t) 
\equiv 
-
\langle
U
,
\Lambda
\rangle
-
\frac{t\beta^2}{2}
\sum_{
k=1
}^{
n
}
x_k
\left(
\Vert
Q^{(k+1)}
\Vert_{\text{F}}^2
-
\Vert
Q^{(k)}
\Vert_{\text{F}}^2
\right)
+
X_0(x, \mathcal{Q}, U, \Lambda)
.
\end{align}
Under the following assumption (at first proposed in 1-D in
\cite{TalagrandParisiFormula2006}), we shall effectively prove that remainder
term almost vanishes on the $\theta$ minimisers of
\eqref{vektor-sk-pasisis-functional}, see Theorem~\ref{as2-an-a-priori-estimate}. 
\begin{assumption}
\label{as2-an-a-priori-estimate}
Let
$
\mathfrak{U} \in \symmetric^+(2d)
$
be defined by \eqref{eq:as2-mathfrak-u-definition}. We fix arbitrary $t_0 \in
[0;1)$,  $\eps > 0$ and $\delta > 0$. There exists
$K=K(t_0,\eps,\delta,\mathfrak{U})>0$, $\theta(t_0,\eps,\delta,\mathfrak{U}) > 0$, and $N_0 =
N_0(t_0,\eps,\delta,\mathfrak{U}) \in \N$ and 
$
\mathfrak{L}^* \in \symmetric(2d)
$
with the following property:

If $(n,\rho^*,\Lambda^*)$ is a $\theta$-optimiser, for some $\theta \in
(0;\theta(t_0,\eps,\delta,\mathfrak{U})]$, 
then uniformly, for all $t \in [0;t_0)$, $N > N_0$
and all $k \in [1;n] \cap \N$, we have 
\begin{align}
\label{eq:as2-a-priori-estimate}
\varphi^{(2)}_N
(
k
,
t
,
x^*,Q^*,
\Sigma^{(2)}_N(\mathfrak{L}^*,\mathfrak{U},\eps,\delta)
)
\leq
2 \phi^{(x^*,\mathcal{Q}^*,\Lambda)}(t)
-
\frac{1}{K}
\Vert
Q^{*(k)}
-
V
\Vert_{\text{F}}^2
+
C (\eps + \delta)
.
\end{align}
\end{assumption}
\begin{remark}
The validity of the above assumption for general a priori measures is an open
problem. However, in the particular case of the Gaussian a priori distribution
the assumption is indeed effectively satisfied. See
Section~\ref{sec:gaussian-spins:gaussian-spins} and
Theorem~\ref{thm:gaussian-spins:a-priori-estimate}, in particular. This gives
a complete proof of the Parisi formula for the case of Gaussian spins.
\end{remark} 
\begin{remark}
\label{rem:remainder:a-condition-for-the-assumption}
If the bound \eqref{eq:remainder:the-doubling-of-the-free-energy} holds
then
Lemma~\ref{lem:remainder:replicated-b-term-in-terms-of-the-non-replicated-one}
with $w = x_{r-1}$
would imply that
\begin{align}
\label{eq:remainder:the-doubled-free-energy-non-sharp-version}
\varphi^{(2)}_{N}
(k,t,\Sigma^{(2)}_N(\mathfrak{L}^*,\mathfrak{U},\eps,\delta))
\overset{?}{\leq}
2
\phi^{(x^*,\mathcal{Q}^*,\Lambda^*)}(t) 
+
C(\eps +\delta)
.
\end{align}
The above inequality would then be a starting point for the a priori
estimates in the spirit of Talagrand~\cite{TalagrandParisiFormula2006} which
might lead to the proof of Assumption~\ref{as2-an-a-priori-estimate}.
\end{remark}

\subsection{Gronwall's inequality and the Parisi formula} 
\begin{theorem}
\label{thm:remainder:the-conditional-parisi-formula}
Suppose Assumption~\ref{as2-an-a-priori-estimate} holds.

Then we have
\begin{align*}
\lim_{
N \uparrow +\infty
}
p_N(\beta)
=
\sup_{
U \in \symmetric^+(d)
}
\inf_{
\substack{
\rho
\in
\mathcal{Q}^\prime(U,d)
\\
\Lambda
\in
\symmetric(d)
}
}
\mathcal{P}(\beta, \rho, \Lambda)
.
\end{align*}
\end{theorem}
\begin{proof}
The proof follows the argument of \cite{TalagrandParisiFormula2006} (see also
\cite{panchenko-free-energy-generalized-sk-2005}) with the adaptations to the
case of multidimensional spins. The main ingredients are the Gronwall inequality and
Lemma~\ref{thm:as2-a-sufficient-condition-for-mu-k-to-vanish}.
Theorem~\ref{thm:multidimensional-sk:free-energy-upper-bound} implies that
\begin{align*}
\lim_{
N \uparrow +\infty
}
p_N(\beta)
\leq
\sup_{
U \in \symmetric^+(d)
}
\inf_{
\substack{
\rho
\in
\mathcal{Q}^\prime(U,d)
\\
\Lambda
\in
\symmetric(d)
}
}
\mathcal{P}(\beta, \rho, \Lambda)
.
\end{align*}
We now turn to the proof of the matching lower bound. As in the proof of
Theorem~\ref{thm:as2:pressure-lower-bound}, it is enough to show that
\begin{align}
\label{eq:remainder:lower-bound-without-remainder}
\lim_{
\eps\downarrow+0
}
\lim_{
N\uparrow+\infty
}
\varphi_N(1,x,Q,B(U,\eps)) 
\geq
\inf_{
\substack{
\rho
\in
\mathcal{Q}^\prime(U,d)
\\
\Lambda
\in
\symmetric(d)
}
}
\mathcal{P}(\beta, \rho, \Lambda)
.
\end{align}
\begin{enumerate}
\item
We fix an arbitrary $U \in \symmetric^+(d)$.
Fix also some $t_0 \in [0;1)$. By Assumption~\ref{as2-an-a-priori-estimate},
we can find the corresponding $\theta(t_0,V,U)>0$ with the properties listed in
the assumption. We pick any $\theta \in (0;\theta(t_0,V,U)]$ and let
$(n,\rho^*,\Lambda^*)$ be a correspondent $\theta$-optimiser.
Note that, by definition \eqref{eq:as2-psi-of-t-definition}, we have
\begin{align*}
\phi^{(x^*,\mathcal{Q}^*,\Lambda^*)}(1)
=
\mathcal{P}(\beta, \rho^*, U, \Lambda^*)
\end{align*}
and, by Definition \ref{def:as2-eps-optimizer},
\begin{align}
\label{eq:as2-almost-minimizer-property}
\vert
\phi^{(x^*,\mathcal{Q}^*,\Lambda^*)}(1)
-
\inf_{
\substack{
\rho
\in
\mathcal{Q}^\prime(U,d)
\\
\Lambda
\in
\symmetric(d)
}
} 
\mathcal{P}(\beta, \rho, U, \Lambda)
\vert
\leq
\theta
.
\end{align}
\item
We denote 
\begin{align*}
\Delta_N(t) 
\equiv 
\phi^{(x^*,\mathcal{Q}^*,\Lambda^*)}(t)
-
\varphi_N(t,x^*,Q^*,B(U,\eps)) 
.
\end{align*}
Note that, due to \eqref{eq:as2:sum-rule}, we obviously have 
\begin{align}
\label{eq:remainder:delta-n-lower-bound}
\Delta_N(t) 
\geq 
- 
C \eps
.
\end{align}
Define
\begin{align*}
\Delta(t) 
\equiv 
\lim_{N \uparrow +\infty}
\Delta_N(t)
.
\end{align*}
The definition
\eqref{eq:as2-psi-of-t-definition} and Theorem
\ref{thm:lectire-03-an-analytic-projection-the-v-term-without-rpc} yield
\begin{align}
\label{eq:as2-derivative-psi-phi}
\frac{\dd}{\dd t}
\Delta_N(t)
\leq
\frac{1}{2}
\sum_{k=0}^{n-1}
(
x_{k+1}
-
x_{k}
)
\mu_k
\left[
\Vert
R_N(\sigma^{(1)},\sigma^{(2)})
-
Q^{(k)}
\Vert_{\text{F}}^2
\right]
+
C \eps
.
\end{align}
\item
Let us set 
$
D 
\equiv
\sup_{
\sigma
\in
\Sigma
}
\Vert
\sigma
\Vert_2 
$. We note that, for any $\sigma^{(1)},\sigma^{(2)} \in \Sigma_N$, we have
\begin{align*}
R(\sigma^{(1)},\sigma^{(2)})
\in 
[
-D^2
;
D^2
]^{
d \times d
}
.
\end{align*}
Given the constant $K$ from \eqref{eq:as2-a-priori-estimate}, for any 
$c>0$, we define the set
\begin{align}
\label{eq:remainder:sigma-squared-k-set}
\Sigma^{(2),k}_N(U,\eps)
\equiv
\left\{
(
\sigma^{(1)}
,
\sigma^{(2)}
)
\in
\Sigma_N(B(U,\eps))^2
:
\Vert
R(\sigma^{(1)},\sigma^{(2)})
-
Q^{(k)}
\Vert_{\text{F}}^2
\geq
2 K
\left(
\Delta_N(t)
+
c
\right)
\right\}
.
\end{align}
It is easy to see that by compactness we can find a finite covering of
$\Sigma^{(2),k}_N(U,\eps)$ by the neighbourhoods
\eqref{eq:as2-sigma-2-n-v-u-delta-eps} with centres, e.g., in the corresponding
set of admissible overlap matrices
\begin{align*}
\mathcal{V}^{(k)}_N(U,\eps)
\equiv
\left\{
R(\sigma^{(1)},\sigma^{(1)})
\in 
[
-D^2
;
D^2
]^{
d \times d
}
:
(
\sigma^{(1)}
,
\sigma^{(2)}
)
\in
\Sigma^{(2),k}_N(U,\eps)
\right\}
.
\end{align*}
That is, there exists $M = M(\eps,\delta) \in \N$ and the finite collections of
matrices 
$
\left\{
V(i)
\right\}_{
i=1
}^{
M
}
\subset
\mathcal{V}^{(k)}_N(U,\eps)
$
and 
$
\left\{
U(i)
\right\}_{
i=1
}^{
M
}
\subset
B(U,\eps) \cap \symmetric^+(d)
$
such that
\begin{align}
\label{eq:remainder:finite-covering}
\Sigma^{(2),k}_N(U,\eps)
\subset
\bigcup_{
i=1
}^{
M
}
\Sigma^{(2)}_N(\mathfrak{L}^*(i),\mathfrak{U}(i),\eps,\delta)
,
\end{align}
where 
\begin{align*}
\mathfrak{U}(i)
\equiv
\begin{bmatrix}
U(i)
& 
V(i)
\\
V^*(i) 
&
U(i)
\end{bmatrix}
\in
\symmetric^+(2d)
,
\end{align*}
and $\mathfrak{L}^*(i)$ is the corresponding $\delta$-minimal Lagrange
multiplier.
\item
Given $i \in [1;M] \cap \N$, let $(n(i),x^*(i),Q^*(i),\Lambda^*(i))$ be
the corresponding to $U(i)$ $\theta(i)$-optimisers. Due to Lipschitzianity of the Parisi
functional (Proposition~\ref{eq:remainder:parisi-functional-lipschitzianity})
and the fact that $U(i) \in B(U,\eps)$ we can assume that $n(i) = n$.
Using the bound
\eqref{eq:as2-a-priori-estimate} and the definition \eqref{eq:remainder:sigma-squared-k-set}, we obtain
\begin{align*}
\varphi^{(2)}_N
(
k
,
t
,
x^*_i,Q^*_i,
\Sigma^{(2)}_N(\mathfrak{L}^*(i),\mathfrak{U}(i),\eps,\delta)
)
&
\leq
2 \phi^{(x^*(i),\mathcal{Q}^*(i),\Lambda^*(i))}(t)
-
\frac{1}{K}
\Vert
Q^{(k)}
-
V(i)
\Vert_{\text{F}}^2
+
C(\eps+\delta)
\\
&
\leq
2 
\varphi_N(t,x^*,Q^*,B(U,\eps)) 
-
c
+
C(\eps+\delta)
,
\end{align*}
where the last inequality is again due to Lipschitzianity of the Parisi
functional (Proposition~\ref{eq:remainder:parisi-functional-lipschitzianity})
which allows to approximate functional's value at
$(x^*(i),Q^*(i),\Lambda^*(i))$ by the value at $(x^*,Q^*,\Lambda^*)$ paying the cost of at most $C \eps$. Choose
$c > C(\eps+\delta)$. Then Lemma~\ref{thm:as2-a-sufficient-condition-for-mu-k-to-vanish} implies that there exists $L = L(\eps,\delta,c)>0$ such that
\begin{align*}
\mu_k
\Bigl(
\Sigma^{(2)}_N(\mathfrak{L}^*,\mathfrak{U},\eps,\delta)
\Bigr)
\leq
L
\exp
\left(
-\frac{N}{L}
\right)
.
\end{align*}
Therefore, the inclusion
\eqref{eq:remainder:finite-covering} gives
\begin{align}
\label{eq:as2-mu-k-of-sigma-u-delta-eps-vanishes}
\mu_k
\Bigl(
\Sigma^{(2),k}_N(U,\eps)
\Bigr)
\leq
L M
\exp
\left(
-\frac{N}{L}
\right)
.
\end{align}
Hence, for each $k \in [1;n]\cap\N$, we have
\begin{align}
\label{eq:as2-mu-k-upper-bound}
\mu_k
\Bigl[
\Vert
R_N(\sigma^{(1)},\sigma^{(2)})
-
Q^{(k)}
\Vert_{\text{F}}^2
\Bigr]
&
=
\mu_k
\Bigl[
\Vert
R_N(\sigma^{(1)},\sigma^{(2)})
-
Q^{(k)}
\Vert_{\text{F}}^2
\I_{
\Sigma^{(2),k}_N(U,\eps)
}
(\sigma^{(1)},\sigma^{(2)})
\Bigr]
\nonumber
\\
&
\quad
+
\mu_k
\Bigl[
\Vert
R_N(\sigma^{(1)},\sigma^{(2)})
-
Q^{(k)}
\Vert_{\text{F}}^2
\left(
1
-
\I_{
\Sigma^{(2),k}_N(U,\eps)
}
(\sigma^{(1)},\sigma^{(2)})
\right)
\Bigr]
\nonumber
\\
&
=:
\text{I}
+
\text{II}
.
\end{align}
For all 
$
(\sigma^{(1)},\sigma^{(2)}) 
\in 
\left(
\Sigma_N(B(U,\eps))^2
\setminus 
\Sigma^{(2),k}_N(U,\eps,\delta)
\right)
$,
we have by definition
\begin{align*}
\Vert
R(\sigma^{(1)},\sigma^{(2)})-Q^{(k)}
\Vert_{\text{F}}^2
<
2 K
\left(
\Delta_N(t)
+
c
\right)
.
\end{align*}
Therefore, using Remark \ref{rem:as2-mu-k-is-a-probability}, we arrive to
\begin{align}
\label{eq:gaussian:ii-leq-2-k-delta-n-t}
\text{II}
\leq
2 K
\left(
\Delta_N(t)
+
c
\right)
.
\end{align}
The bound \eqref{eq:as2-mu-k-of-sigma-u-delta-eps-vanishes} assures that
\begin{align}
\label{eq:gaussian:i-leq-l-m-exp-minus-n-over-l}
\text{I}
\leq
L M
\exp
\left(
-\frac{N}{L}
\right)
.
\end{align}
\item
Combining \eqref{eq:gaussian:ii-leq-2-k-delta-n-t} and
\eqref{eq:gaussian:i-leq-l-m-exp-minus-n-over-l} with \eqref{eq:as2-mu-k-upper-bound} and \eqref{eq:as2-derivative-psi-phi}, we obtain
\begin{align*}
\frac{\dd}{\dd t}
\Delta_N(t)
\leq
2 K
\left(
\Delta_N(t)
+
c
\right)
+
L M
\exp
\left(
-\frac{N}{L}
\right)
+
C(\eps+\delta)
.
\end{align*}
Hence,
\begin{align*}
\frac{\dd }{\dd t}
\Bigl(
(
\Delta_N(t)
+
c
)
\exp(-2 K t)
\Bigr)
&
=
\exp(-2 K t)
\Bigl(
\frac{\dd }{\dd t}
(
\Delta_N(t)
+
c
)
-
2 K
(
\Delta_N(t)
+
c
)
\Bigr)
\\
&
\leq
\exp(-2 K t)
\Bigl(
\frac{\dd }{\dd t}
(
L M
\exp
\left(
-\frac{N}{L}
\right)
+
C(\eps+\delta)
\Bigr)
.
\end{align*}
Integrating the above inequality and noting that due to \eqref{eq:as2:sum-rule}
$
\vert\Delta_N(0)\vert 
\leq
C \eps
$,
we arrive to
\begin{align*}
\Delta_N(t)
+
c
\leq
&
(
C \eps
+
c
)
\exp(-2 K t)
+
L M
\exp
\left(
-\frac{N}{L}
\right)
\\
&
+
C(\eps+\delta)
(
\exp(-2 K t) - 1
)
+
C(\eps+\delta)
.
\end{align*}
Passing consequently to the limits $N \uparrow +\infty$,
$\eps \downarrow +0$, $\delta \downarrow +0$ and finally $c \downarrow +0$ in
the above inequality, we get
\begin{align*}
\lim_{
\eps \downarrow +0
}
\Delta(t)
\leq 
0
,
\quad
\text{for all $t \in [0;t_0]$.}
\end{align*}
The existence of the $N \uparrow +\infty$ limits is guaranteed by the
general result of Guerra and Toninelli
\cite{Guerra-Toninelli-Generalized-SK-2003}. The limits $\eps \downarrow +0$,
$\delta \downarrow +0$ exist due to monotonicity. 
Finally, combining the above inequality with
\eqref{eq:remainder:delta-n-lower-bound}, we get
\begin{align}
\label{eq:as2-delta-of-t-equals-zero-below-t0}
\lim_{
\eps \downarrow +0
}
\Delta(t)
=
0
,
\quad
\text{for all $t \in [0;t_0]$.}
\end{align}
\item
Now, it is easy to extend the validity of
\eqref{eq:as2-delta-of-t-equals-zero-below-t0} onto the
whole interval $[0;1]$.
Indeed, due to the boundedness of the derivatives of $\varphi_N$ and $\phi$, we
have, for any $t \in [0;1]$,
\begin{align}
\label{eq:remainder:delta-n-of-t-leq-somethin}
\Delta_N(t)
&
\leq
\int_0^1
\frac{\dd}{\dd t}
\Delta_N(t)
\dd t
\nonumber
\\
&
=
\left(
\int_0^{t_0}
+
\int_{t_0}^1
\right)
\frac{\dd}{\dd t}
\Delta_N(t)
\dd t
\nonumber
\\
&
\leq
\left(
\Delta_N(t_0)
-
\Delta_N(0)
\right)
+
\int_{t_0}^1
\left\vert
\frac{\dd}{\dd t}
\Delta_N(t)
\right\vert
\dd t
\nonumber
\\
&
\leq
\Delta_N(t_0)
+
L(1-t_0)
.
\end{align}
Passing to the $N \uparrow +\infty$ limit, applying
\eqref{eq:as2-delta-of-t-equals-zero-below-t0}, and then to $t_0 \to 1$ limit
in \eqref{eq:remainder:delta-n-of-t-leq-somethin}, we get
\begin{align*}
\lim_{\eps \downarrow +0}
\Delta(t)=0
,
\quad
\text{
for all $t \in [0;1]$.
}
\end{align*}
\item
In particular, the previous formula yields 
\begin{align*}
0 
=
\lim_{\eps \downarrow +0}
\Delta(1)
=
\phi^{(x^*,\mathcal{Q}^*,\Lambda^*)}(1)
-
\lim_{\eps \downarrow +0}
\varphi_N(1,x^*,Q^*,B(U,\eps)) 
. 
\end{align*}
Note that $\varphi_N(1,x,Q,B(U,\eps))$ does not depend on the choice of $x$ and
$Q$. Hence, by \eqref{eq:as2-almost-minimizer-property}, we obtain
\begin{align*}
\vert
\lim_{\eps \downarrow +0}
\varphi_N(1,x^*,Q^*,B(U,\eps))
-
\inf_{
\substack{
\rho
\in
\mathcal{Q}^\prime(U,d)
\\
\Lambda
\in
\symmetric(d)
}
} 
\mathcal{P}(\beta, \rho, U, \Lambda)
\vert
\leq
\theta
.
\end{align*}
The proof of \eqref{eq:remainder:lower-bound-without-remainder} is finished by
noticing that the $\theta$ can be made arbitrary small.
\end{enumerate}
\end{proof}

\section{Proof of the local Parisi formula for the SK model with
multidimensional Gaussian spins}
\label{sec:gaussian-spins:gaussian-spins}
In this section, we prove 
Theorem~\ref{thm:gaussian:the-local-low-temperature-parisi-formula}. The rich
symmetries of the Gaussian a priori distribution allow rather explicit
computations of the $X_0$ terms (see \eqref{eq:multidimensional-sk:x-0}). This
allows us to prove that the analogon of
Assumption~\ref{as2-an-a-priori-estimate} is satisfied, implying the Parisi
formula for the local free energy (Theorem~\ref{thm:gaussian:the-local-low-temperature-parisi-formula}).

\begin{remark}
The case of Gaussian spins is very tractable due to the (unusually) good
symmetry (i.e., the rotational invariance) of the Gaussian measure.
Therefore, it is not surprising that in this case the calculus resembles the
one for the spherical SK model, cf. \cite{PanchenkoTalagrandMultipleSKModels2006,TalagrandSphericalSK}.
\end{remark}

We start from the estimates under a generic (i.e., no simultaneous
diagonalisation, cf.
Section~\ref{sec:remainder:simultaneous-diagonalisation-scenario}) scenario.

\subsection{The case of positive increments}
Let, for $k \in [0;n] \cap \N$,
\begin{align*}
\Delta Q^{(k)}
\equiv
Q^{(k+1)}-Q^{(k)}
.
\end{align*}
We define,
for 
$
\Lambda
\in
\symmetric(d)
$, 
a family of matrices
$
\left\{
D^{(l)}
\in
\R^{d \times d}
\right\}_{l=0}^{n+1}
$
as follows
\begin{align*}
D^{(n+1)}
\equiv
C
,
\end{align*}
and, further, for $k \in [0;n] \cap \N$,
\begin{align}
\label{eq:lecture-03:spherical-d-matrices-of-the-a-part}
D^{(k)}
\equiv
C
-
\Lambda
-
2 \beta^2
\sum_{
l = k 
}^n
x_l
\Delta
Q^{(l)}
.
\end{align}
We assume that the matrices $\Lambda$ and $C$ are such that, for all 
$
l 
\in
[1;n+1] \cap \N
$, we have
\begin{align*}
D^{(l)} 
\succ 0
.
\end{align*}
We need the following two small (and surely known) technical Lemmata which
exploit the symmetries of our Gaussian setting. We include their statements for
reader's convenience.
\begin{lemma}
\label{lem:lecture-03:gaussian-integral-with-linear-form}
Fix some vector $h \in \R^d$
and a Gaussian random vector $z \in \R^d$ with $\var z  = C^{-1} \in \R^{d
\times d}$.
 
Then we have
\begin{align*}
\E^z
\left[
\exp
\left(
\langle
z
,
h
\rangle
+
\langle
\Lambda
\sigma
,
\sigma
\rangle
\right)
\right]
=
&
\left(
\det
\left[
C
\left(
C-\Lambda
\right)^{-1}
\right]
\right)^{
1/2
}
\\
&
\times
\exp
\left(
\frac{1}{2}
\left\langle
(C-\Lambda)^{-1}
h
,
h
\right\rangle
\right)
.
\end{align*}
\end{lemma}
\begin{proof}
This is a standard Gaussian averaging argument.

\end{proof}
\begin{lemma}
\label{thm:lecture-03:gaussian-integral}
For a positive definite matrix $\Delta Q \in \symmetric(d)$, let $z \sim
\mathcal{N}(0,\Delta Q)$. We fix also another positive definite matrix $D \in
\symmetric(d)$ such that $\Delta Q^{-1} \succ D^{-1}$. 

Then we have 
\begin{align*}
\E^z
\left[
\exp
\left(
\frac{1}{2}
\langle
D^{-1}
(z+h)
,
z+h
\rangle
\right)
\right]
=
&
\left(
\det
\left[
D (D-\Delta Q)^{-1}
\right]
\right)^{-1/2}
\\
&
\times
\int_{
\R^d
}
\exp
\left(
\frac{1}{2}
\langle
(D-\Delta Q)^{-1}
h
,
h
\rangle
\right)
.
\end{align*}
\end{lemma}
\begin{proof}
This is a standard Gaussian averaging argument. See, e.g.,
\cite{TalagrandSphericalSK} for an argument in 1-D.

\end{proof}
Now we are ready to compute the term $X_0(x, \mathcal{Q}, U, \Lambda)$ 
(see \eqref{eq:multidimensional-sk:x-0}) corresponding to the a priori
distribution \eqref{gaussian-spins:a-priori-measure} in a rather explicit way.
\begin{lemma}
\label{lem:gaussian:a-formula-for-the-a-term}
We have
\begin{align*}
X_0(x, \mathcal{Q}, U, \Lambda)
=
\frac{1}{2}
\left(
\langle
[D^{(1)}]^{-1}
,
\Delta Q^{(0)}
\rangle
+
\langle
[D^{(1)}]^{-1}
h
,
h
\rangle
+
\sum_{
l = 1
}^{
n
}
\frac{1}{x_l}
\log
\left(
\frac{
\det D^{(l+1)}
}{
\det D^{(l)}
}
\right)
\right)
.
\end{align*}
\end{lemma}
\begin{proof}
\begin{enumerate}
\item 
We start from computing the following quantity
\begin{align}
\label{eq:lecture-03:gaussian-spins}
X_{n+1}
\equiv
\log
\int_{
\R^d
}
\exp
\left(
\sum_{l=0}^{n}
\langle
Y^{(l)}
,
\sigma
\rangle
+
\langle
\Lambda
\sigma
,
\sigma
\rangle
\right)
\dd \mu(\sigma)
,
\end{align}
where $Y^{(l)} \in \R^d$ are independent Gaussian vectors with variance 
\begin{align*}
\var
\left[
Y^{(l)}
\right]
=
2
\beta^2
\Delta
Q^{(l)}
.
\end{align*}
We denote 
\begin{align*}
\widetilde{h} 
\equiv 
h
+
\sum_{l=0}^{n}
Y^{(l)}
.
\end{align*}
Lemma \ref{lem:lecture-03:gaussian-integral-with-linear-form} gives
\begin{align*}
\int_{
\R^d
}
\exp
\left(
\sum_{l=0}^{n}
\langle
Y^{(l)}
,
\sigma
\rangle
+
\langle
\Lambda
\sigma
,
\sigma
\rangle
\right)
\dd \mu(\sigma)
=
&
\left(
\det
\left[
C
\left(
C-\Lambda
\right)^{-1}
\right]
\right)^{
1/2
}
\\
&
\times
\exp
\left(
\frac{1}{2}
\left\langle
(C-\Lambda)^{-1}
\widetilde{h}
,
\widetilde{h}
\right\rangle
\right)
.
\end{align*}
\item
Next, we define, for $l \in [0;n] \cap \N$,  recursively the following
quantities 
\begin{align*}
X_l
\equiv
\frac{1}{x_l}
\log
\E^{
Y_l
}
\left[
\exp
\left(
x_l
X_{l+1}
\right)
\right]
.
\end{align*}
Applying the Lemma \ref{thm:lecture-03:gaussian-integral} to \eqref{eq:lecture-03:gaussian-spins} 
recursively, we obtain 
\begin{align}
\label{eq:lecture-03:x-1-equals}
X_1
\equiv
\frac{1}{2}
\langle
[(D^{(1)}]^{-1}
\left(
Y^{(0)}
+
h
\right)
,
Y^{(0)}
+
h
\rangle
+
\frac{1}{2}
\sum_{
l = 1
}^{
n
}
\frac{1}{x_l}
\log
\left(
\frac{
\det D^{(l+1)}
}{
\det D^{(l)}
}
\right)
.
\end{align}
Recall that we have
\begin{align}
\label{eq:gaussian:x-0-equals-somethin}
X_0 
&
=
\lim_{
x \to +0
} 
\frac{1}{x}
\log
\E^{
Y_0
}
\left[
\exp
\left(
x
X_{1}
\right)
\right]
\nonumber
\\
&
=
\E^{
Y_0
}
\left[
X_{1}
\right]
\end{align}
and note that
\begin{align}
\label{eq:gaussian:e-y-0}
\E^{Y_0}
\left[
\langle
[D^{(1)}]^{-1}
(
Y^{(0)}
+
h
)
,
Y^{(0)}
+
h
\rangle
\right]
=
2
\beta^2
\langle
[D^{(1)}]^{-1}
,
\Delta Q^{(0)}
\rangle
+
\langle
[D^{(1)}]^{-1}
h
,
h
\rangle
.
\end{align}
Hence, combining \eqref{eq:gaussian:x-0-equals-somethin} and
\eqref{eq:gaussian:e-y-0} with \eqref{eq:lecture-03:x-1-equals}, we obtain the theorem.

\end{enumerate}
\end{proof}

\subsection{Simultaneous diagonalisation scenario}
\label{sec:gaussian:simultaneous-diagonalisation-scenario}
In what follows, we employ the simultaneous diagonalisation scenario
introduced in
Section~\ref{sec:remainder:simultaneous-diagonalisation-scenario}. Suppose
that, for $l \in [0;n+1] \cap \N$, and some matrix $O \in \mathcal{O}(d)$, we have
\begin{align*}
D^{(l)}
\equiv
O^*
d^{(l)}
O
,
\end{align*}
where the vectors $d^{(l)} \in \R^d$, for $l \in [0;n] \cap \N$, satisfy
\begin{align*}
0
\prec
d^{(l)}
\prec
d^{(l+1)}
.
\end{align*}
That is, the vectors $d^{(l)}$ are (component-wise) increasingly ordered and
non-negative.
\begin{lemma}
\label{lem:lecture-03:gaussian-spins:a-decoupling-of-coordinates-for-the-a-and-b-terms-in-the-diagonal-case}
We have
\begin{align}
\label{eq:lecture-03:decoupling-of-coordinates-for-the-a-term}
X_0(x, \mathcal{Q}, U, \Lambda)
&
=
\frac{1}{2}
\sum_{v=1}^d
\left(
\frac{
2\beta^2
q^{(1)}_v
+
h_v^2
}{
d^{(1)}_v
}
+
\sum_{
l = 1
}^{
n
}
\frac{1}{x_l}
\log
\left(
\frac{
d^{(l+1)}_v
}{
d^{(l)}_v
}
\right)
\right)
,
\\
\label{eq:lecture-03:decoupling-of-coordinates-for-the-b-term}
\frac{\beta^2}{2}
\sum_{
k=1
}^{
n
}
x_k
\left(
\Vert
Q^{(k+1)}
\Vert_\text{F}^2
-
\Vert
Q^{(k)}
\Vert_\text{F}^2
\right)
&
=
\frac{
\beta^2
}{
2
}
\sum_{
k=1
}^{
n
}
x_l
\left(
\Vert
q^{(k+1)}
\Vert_{2}^2
-
\Vert
q^{(k)}
\Vert_{2}^2
\right)
.
\end{align}
\end{lemma}
\begin{proof}
This is a standard argument which relies on the standard invariance
properties of the determinant and the matrix trace. 

\end{proof}
Define the 1-D Parisi functional for the case
\eqref{gaussian-spins:a-priori-measure} as
\begin{align}
\label{eq:lecture-03:parisi-functinal-diagonal-case}
\mathcal{P}(\rho,\lambda)
\equiv
&
-
\lambda
u
+
\frac{
2\beta^2
q^{(1)}
+
h^2
}{
d^{(1)}
}
+
\sum_{
l = 1
}^{
n
}
\frac{1}{x_l}
\log
\left(
\frac{
d^{(l+1)}
}{
d^{(l)}
}
\right)
\nonumber
\\
&
-
\beta^2
\sum_{
l=1
}^{
n
}
x_l
\left(
[q^{(l+1)}]^2
-
[q^{(l)}]^2
\right)
.
\end{align}
\begin{proposition}
\label{prp:gaussian:guerras-scheme}
There exists $C = C(\Sigma) > 0$ such that,
for all $u \in \R_+^d$ and all $\eps, \delta > 0$,
there exists an $\delta$-minimal Lagrange multiplier
$
\lambda
=
\lambda(U,\eps,\delta) 
\in 
\R^d
$
in \eqref{eq:as2:local-parisi-functional}
such that, for all $t \in [0;1]$ 
and all $(x,\mathcal{\rho})$,
we have
\begin{align}
\label{eq:gaussian:guerras-sum-rule-ldp-upper-bound}
p_N(
\Sigma_N(\mathcal{V}(\Lambda,U,\eps,\delta))
)
\leq
&
\frac{1}{2}
\inf_{
\rho,
\lambda
}
\left(
\sum_{v=1}^d
\mathcal{P}(\rho_v,\lambda_v)
\right)
+
C(\eps + \delta)
\end{align}
and
\begin{align}
\label{eq:gaussian:guerras-sum-rule-ldp-lower-bound}
\lim_{
N \uparrow +\infty
}
p_N(
\Sigma_N(B(U,\eps))
)
\geq
&
\frac{1}{2}
\inf_{
\rho,
\lambda
}
\left(
\sum_{v=1}^d
\mathcal{P}(\rho_v,\lambda_v)
+
\lim_{
N \uparrow +\infty
}
\int_{0}^{1}
\mathcal{R}(t,x,Q,\Sigma_N(B(U,\eps)))
\dd t
\right)
\nonumber
\\
&
+
C(\eps + \delta)
.
\end{align}
\end{proposition}
\begin{proof}
We combine \eqref{eq:lecture-03:decoupling-of-coordinates-for-the-a-term} and
\eqref{eq:lecture-03:decoupling-of-coordinates-for-the-b-term} and the
Proposition~\ref{eq:remainder:guerras-upper-and-lower-bounds} to get
\eqref{eq:gaussian:guerras-sum-rule-ldp-upper-bound} and
\eqref{eq:gaussian:guerras-sum-rule-ldp-lower-bound}.
\end{proof}

\subsection{The Crisanti-Sommers functional in 1-D}
In this subsection, we 
adapt the proof of \cite{TalagrandSphericalSK} to obtain the
equivalence between the (very tractable) Crisanti-Sommers functional
\cite{CrisantiSommers1992} and the Parisi one
\eqref{eq:lecture-03:parisi-functinal-diagonal-case} in the case of the Gaussian a
priori measure \eqref{gaussian-spins:a-priori-measure}. Similar ideas based on
the symmetry of the a priori measure were exploited in the case of the
spherical models by
\cite{BenArousDemboGuionnet2001,PanchenkoTalagrandMultipleSKModels2006}.

We restrict the consideration to 1-D situation for a moment. Given $u \geq 0$,
consider 
$
\rho
\in 
\mathcal{Q}^\prime_n(u,1)
$,  
$\lambda \in \R$, $h \in \R$ 
and let 
$
\{
d^{(l)}
\in \R
\}_{l=1}^{n+1}
$
be the scalars playing the role of matrices $D^{(l)}$ (cf. 
\eqref{eq:lecture-03:spherical-d-matrices-of-the-a-part}). That is,
\begin{align*}
d^{(l)}
&
\equiv
c
-
\lambda
-
2\beta^2
\sum_{k=l}^{n}
x_k
\left(
q^{(k+1)}
-
q^{(k)}
\right)
,
\\
d^{(n+1)}
&
\equiv
c
.
\end{align*}
We define, for $k \in [1;n] \cap \N$, the family of vectors 
$
\{
s^{(k)} 
\in 
\R^d 
\}_{k=0}^{n}
$ 
by
\begin{align}
\label{eq:lecture-03:spherical-quantity-s-k}
s^{(k)}
\equiv
\sum_{
l=k
}^{n}
x_l
\left(
q^{(l+1)}
-
q^{(l)}
\right)
.
\end{align}
We also define the Crisanti-Sommers functional as follows
\begin{align}
\label{eq:lecture-03:crisanti-sommers-functional}
\mathcal{CS}(\rho)
\equiv
&
1
-
c u
+
h^2
s^{(1)}
+
\frac{
q^{(1)}
}{
s^{(1)}
}
+
\sum_{l=1}^{n-1}
\frac{1}{x_l}
\log
\left(
\frac{
s^{(l)}
}{
s^{(l+1)}
}
\right)
+
\log
\left[
c
(
u
-
q^{(n)}
)
\right]
\nonumber
\\
&
+
\beta^2
\sum_{l=1}^n
x_l
\left(
[q^{(l+1)}]^2
-
[q^{(l)}]^2
\right)
.
\end{align}
\begin{lemma}
\label{lem:some-extremal-properties-of-the-parisi-functional}
If $(\rho,\lambda)$ is an optimiser for
\eqref{eq:lecture-03:parisi-functinal-diagonal-case}, that is,
\begin{align}
\label{eq:lecture-03:spherical-parisi-optimization}
\mathcal{P}(\rho,\lambda) =
\inf_{
(\rho', \lambda')
}
\mathcal{P}(\rho',\lambda')
,
\end{align}
then, for all $k \in [1;n] \cap \N$, the pair $(\rho,\lambda)$ satisfies
\begin{align}
\label{eq:lecture-03:spherical-model-q-k}
q^{(k)}
=
\frac{
h^2
+
2\beta^2
q^{(1)}
}{
[d^{(1)}]^2
}
+
\sum_{l=1}^{k-1}
\frac{1}{x_l}
\left(
\frac{1}{
d^{(l)}
}
-
\frac{1}{
d^{(l+1)}
}
\right)
.
\end{align}
Moreover,
\begin{align}
\label{eq:lecture-03:spherical-lambda-optimiality-consequence}
\lambda
=
c
-
2\beta^2
(
u
-
q^{(n)}
)
-
(
u
-
q^{(n)}
)^{-1}
,
\end{align}
and, for all $k \in [1;n] \cap \N$, we have
\begin{align}
\label{eq:lecture-03:relation-between-crisanti-sommers-and-parisi-a-parts}
\frac{1}{
s^{(k+1)}
}
-
\frac{1}{
s^{(k)}
}
=
2\beta^2
x_k
\left(
q^{(k+1)}
-
q^{(k)}
\right)
,
\end{align}
and also
\begin{align}
\label{eq:lecture-03:gausian-spins:s-k-equals-inverse-d-k}
s^{(k)}
=
\frac{1}{
d^{(k)}
}
.
\end{align}
\end{lemma}
\begin{remark}
In the formulation of the theorem (as well as elsewhere), it is implicit that 
$
d^{(k)}=d^{(k)}(\rho,\lambda)
$ 
and
$
s^{(k)}=s^{(k)}(\rho,\lambda)
$.
\end{remark}
\begin{proof}
\begin{enumerate}
\item 
Rearranging the terms in \eqref{eq:lecture-03:parisi-functinal-diagonal-case}, we observe that
\begin{align}
\label{eq:lecture-03:sperical-parisi-rearrangement}
\mathcal{P}(\rho',\lambda')
=
&
-
\lambda
u
+
\frac{
2\beta^2
q^{(1)}
+
h^2
}{
d^{(1)}
}
+
\sum_{
l = 2
}^{
n
}
\log d^{(l)}
\left(
\frac{1}{x_{l-1}}
-
\frac{1}{x_{l}}
\right)
+
\frac{1}{x_n}
\log
d^{(n+1)}
-
\frac{1}{x_1}
\log
d^{(1)}
\nonumber
\\
&
-
\beta^2
\sum_{
l=1
}^{
n
}
x_l
\left(
[q^{(l+1)}]^2
-
[q^{(l)}]^2
\right)
.
\end{align}
We compute, for $k,l \in [1;n] \cap \N$,
\begin{align}
\label{eq:gaussian:partial-d-l-partial-q-k}
\frac{
\partial
d^{(l)}
}{
\partial
q^{(k)}
}
=
\begin{cases}
0
,
& 
k < l
,
\\
2\beta^2
x_k
,
&
l = k
,
\\
2\beta^2
\left(
x_k-x_{k-1}
\right)
,
&
k > l
.
\end{cases}
\end{align}
Using  \eqref{eq:gaussian:partial-d-l-partial-q-k} and the
representation \eqref{eq:lecture-03:sperical-parisi-rearrangement}, we compute the necessary
condition for $(q,\lambda)$ satisfy
\eqref{eq:remainder:optimality-for-fixed-n}, for $k \in [2;n] \cap \N$,
\begin{align}
\label{eq:gaussian:partial-p-partial-q-k}
0
=
\frac{
\partial
}{
\partial
q^{(k)}
}
\mathcal{P}(q,\lambda)
=
&
2\beta^2
\left(
x_k
-
x_{k-1}
\right)
\left[
-
\frac{
2\beta^2
q^{(1)}
+
h^2
}{
[d^{(1)}]^2
}
+
\sum_{l=2}^{k-1}
\frac{1}{d^{(l)}}
\left(
\frac{1}{x_{l-1}}
-
\frac{1}{x_{l}}
\right)
\right.
\nonumber
\\
&
+
\left.
\frac{1}{
d^{(k)}
x_{k-1}
}
-
\frac{1}{
x_1
d^{(1)}
}
+
q_k
\right]
.
\end{align}
We also have (for $k=1$)
\begin{align}
\label{eq:gaussian:partial-p-partial-q-1}
0
=
\frac{
\partial
}{
\partial
q^{(1)}
}
\mathcal{P}(q,\lambda)
&
=
2\beta^2
\left[
\frac{
d^{(1)}
-
x_1
\left(
q^{(1)}
+
h^2
\right)
}{
[d^{(1)}]^2
}
-
\frac{x_1}{
x_1
d^{(1)}
}
+
x_1 
q^{(1)}
\right]
\nonumber
\\
&
=
2\beta^2
x_1 
\left[
q^{(1)}
-
\frac{
\left(
q^{(1)}
+
h^2
\right)
}{
[d^{(1)}]^2
}
\right]
.
\end{align}
Relations \eqref{eq:gaussian:partial-p-partial-q-k} and
\eqref{eq:gaussian:partial-p-partial-q-1} then imply \eqref{eq:lecture-03:spherical-model-q-k}.
\item
Using the fact that
\begin{align*}
\frac{
\partial
d^{(l)}
}{
\partial
\lambda
}
=
-1
,
\end{align*}
we obtain
\begin{align}
\label{eq:gaussian:partial-p-partial-lambda}
\frac{
\partial
}{
\partial
\lambda
}
\mathcal{P}(q,\lambda)
=
-
u
+
\frac{
h^2
+
2\beta^2
q^{(1)}
}{
[d^{(1)}]^2
}
+
\sum_{l=1}^{n-1}
\frac{1}{x_l}
\left(
\frac{1}{
d^{(l)}
}
-
\frac{1}{
d^{(l+1)}
}
\right)
+
\frac{1}{
d^{(n)}
}
.
\end{align}
Applying \eqref{eq:lecture-03:spherical-model-q-k} with $k=n$ in
\eqref{eq:gaussian:partial-p-partial-lambda}, we obtain that the necessary
condition for $\lambda$ to satisfy \eqref{eq:lecture-03:spherical-parisi-optimization} is as follows
\begin{align}
\label{eq:lecture-03:spherical-lambda-omptimality-condition}
0
&
=
\frac{
\partial
}{
\partial
\lambda
}
\mathcal{P}(q,\lambda)
=
-u
+
q^{(n)}
+
\frac{1}{x_n}
\left(
\frac{1}{
d^{(n)}
}
-
\frac{1}{
d^{(n+1)}
}
\right)
\nonumber
\\
&
=
-u
+
q^{(n)}
+
\frac{1}{
d^{(n)}
}
=
-u
+
q^{(n)}
+
\left(
c
-
\lambda
-
2\beta^2
(
u
-
q^{(n)}
)
\right)^{-1}
\end{align}
which implies \eqref{eq:lecture-03:spherical-lambda-optimiality-consequence}.
\item
Relation
\eqref{eq:lecture-03:relation-between-crisanti-sommers-and-parisi-a-parts} is
proved as follows. Subtracting the relations 
\eqref{eq:lecture-03:spherical-model-q-k}, we obtain, for $k \in [1;n-1] \cap
\N$,
\begin{align}
\label{eq:lecture-03:x-l-delta-q-k-equals-delta-inverse-d-k}
x_k
\left(
q^{(k+1)}
-
q^{(k)}
\right)
=
\frac{1}{
d^{(k)}
}
-
\frac{1}{
d^{(k+1)}
}
.
\end{align}
By \eqref{eq:lecture-03:spherical-lambda-omptimality-condition}, we have
\begin{align*}
x_n
\left(
q^{(n+1)}
-
q^{(n)}
\right)
=
u
-
q^{(n)}
=
\frac{1}{
d^{(n)}
}
.
\end{align*}
(That is, \eqref{eq:lecture-03:x-l-delta-q-k-equals-delta-inverse-d-k} is valid also for $k=n$.)
Combining the previous two relations, we get, for $k \in [1;n] \cap \N$,
\begin{align}
\label{eq:lecture-03:s-k-equals-inverse-d-k}
s^{(k)}
=
\frac{1}{
d^{(k)}
}
.
\end{align}
Using \eqref{eq:lecture-03:s-k-equals-inverse-d-k} and
\eqref{eq:lecture-03:x-l-delta-q-k-equals-delta-inverse-d-k}, we get 
\begin{align*}
2\beta^2
x_k
\left(
q^{(k+1)}
-
q^{(k)}
\right)
&
=
d^{(k+1)}
-
d^{(k)}
\\
\text{
(by \eqref{eq:lecture-03:x-l-delta-q-k-equals-delta-inverse-d-k})
}
&
=
d^{(k+1)}
d^{(k)}
x_k
\left(
q^{(k+1)}
-
q^{(k)}
\right)
=
d^{(k+1)}
d^{(k)}
\left(
s^{(k)}
-
s^{(k+1)}
\right)
\\
\text{
(by \eqref{eq:lecture-03:s-k-equals-inverse-d-k})
}
&
=
\frac{1}{
s^{(l+1)}
}
-
\frac{1}{
s^{(l)}
}
\end{align*}
which is
\eqref{eq:lecture-03:relation-between-crisanti-sommers-and-parisi-a-parts}.

\end{enumerate}
\end{proof}

\begin{lemma}
\label{lem:lecture-03:some-extremal-properties-of-the-crisanti-sommers-functional}
If $\rho$ is an optimiser of \eqref{eq:lecture-03:crisanti-sommers-functional},
that is,
\begin{align*}
\mathcal{CS}(\rho) =
\inf_{\rho'}
\mathcal{CS}(\rho')
,
\end{align*}
then, for all $l \in [1;n] \cap \N$,
\eqref{eq:lecture-03:relation-between-crisanti-sommers-and-parisi-a-parts} holds.
\end{lemma}
\begin{proof}
The strategy is the same as in the previous lemma. We rearrange the summands in
\eqref{eq:lecture-03:crisanti-sommers-functional} to get
\begin{align}
\label{eq:lecture-03:crisanti-sommers-regrouped-representation}
\mathcal{CS}(\rho)
=
&
h^2
s^{(1)}
+
\frac{
q^{(1)}
}{
s^{(1)}
}
+
\frac{
\log s^{(1)}
}{
x_1
}
-
\frac{
\log s^{(n)}
}{
x_{n-1}
}
+
\sum_{l=2}^{n-1}
\left(
\frac{1}{x_l}
-
\frac{1}{x_{l+1}}
\right)
\log
s^{(l)}
\nonumber
\\
&
+
\log
\left(
c
(
u
-
q^{(n)}
)
\right)
+
\beta^2
\sum_{l=1}^n
x_l
\left(
[q^{(l+1)}]^2
-
[q^{(l)}]^2
\right)
.
\end{align}
We have, for $k,l \in [1;n] \cap \N$,
\begin{align}
\label{eq:lecture-03:s-l-derivative}
\frac{
\partial
s^{(l)}
}{
\partial
q^{(k)}
}
=
\begin{cases}
0
,
&
k < l
,
\\
-x_k
,
&
k=l
,
\\
x_{k-1}
-
x_k
,
&
k > l
.
\end{cases}
\end{align}
\begin{enumerate}
\item 
Relation \eqref{eq:lecture-03:s-l-derivative} implies, for $k \in [2;n-1] \cap
\N$,
\begin{align*}
\frac{
\partial
}{
\partial
q^{(k)}
}
\mathcal{CS}(\rho)
=
&
h^2
(
x_{k-1}
-
x_k
)
-
\frac{
q^{(1)}
}{
[
s^{(1)}
]^2
}
(
x_{k-1}
-
x_k
)
+
\frac{
x_{k-1}
-
x_k
}{
x_1
s^{(1)}
}
\\
&
+
\sum_{l=2}^{k-1}
\frac{
x_{k-1}
-
x_k
}{
s^{(l)}
}
\left(
\frac{1}{
x_l
}
-
\frac{1}{
x_{l-1}
}
\right)
-
\frac{
x_k
}{
s^{(k)}
}
\left(
\frac{1}{
x_k
}
-
\frac{1}{
x_{k-1}
}
\right)
\\
&
+
2\beta^2
q^{(k)}
\left(
x_{k-1}
-
x_k
\right)
=
0
.
\end{align*}
Hence,
\begin{align}
\label{eq:lecture-03:spherical-crisanti-sommers-2-beta-quadrat-q-k}
2\beta^2
q^{(k)}
&
=
-
h^2
+
\frac{
q^{(1)}
}{
[
s^{(1)}
]^2
}
-
\frac{1}{
x_1
s^{(1)}
}
+
\frac{1}{
x_{k-1}
s^{(k)}
}
-
\sum_{l=2}^{k-1}
\frac{1}{
s^{(l)}
}
\left(
\frac{1}{
x_l
}
-
\frac{1}{
x_{l-1}
}
\right)
\nonumber
\\
&
=
-
h^2
+
\frac{
q^{(1)}
}{
[
s^{(1)}
]^2
}
-
\sum_{l=1}^{k-1}
\frac{1}{
x_l
}
\left(
\frac{1}{
s^{(l)}
}
-
\frac{1}{
s^{(l+1)}
}
\right)
.
\end{align}
\item
To handle the case $k=n$, we note that
\begin{align*}
\log 
\left(
1
+
c
(
u
-
q^{(n)}
)
\right)
=
\frac{1}{x_n}
\log
\left(
\frac{
s^{(n)}
}{
s^{(n+1)}
}
\right)
,
\end{align*}
and, hence, the argument in the previous item shows that
\eqref{eq:lecture-03:spherical-crisanti-sommers-2-beta-quadrat-q-k} is also
valid for $k=n$.
\item
Differentiating the representation
\eqref{eq:lecture-03:crisanti-sommers-regrouped-representation} with respect to
$q^{(1)}$ and using \eqref{eq:lecture-03:s-l-derivative}, we obtain
\begin{align*}
\frac{
\partial
}{
\partial
q^{(1)}
}
\mathcal{CS}(\rho)
=
-
x_1
h^2
+
\frac{1}{
s^{(1)}
}
+
\frac{
x_1
q^{(1)}
}{
[
s^{(1)}
]^2
}
-
\frac{x_1}{
x_1
s^{(1)}
}
-
2\beta^2
x_1
q^{(1)}
=
0
.
\end{align*}
Therefore,
\begin{align*}
2\beta^2
q^{(1)}
=
-h^2
+
\frac{
q^{(1)}
}{
[
s^{(1)}
]^2
}
\end{align*}
which is \eqref{eq:lecture-03:spherical-crisanti-sommers-2-beta-quadrat-q-k}, for $k=1$.
\item
Subtracting equations
\eqref{eq:lecture-03:spherical-crisanti-sommers-2-beta-quadrat-q-k}, we arrive
to \eqref{eq:lecture-03:relation-between-crisanti-sommers-and-parisi-a-parts}, for
all $k \in [1;n] \cap \N$.

\end{enumerate}
\end{proof}
\begin{proposition}
\label{thm:lecutre-03:equivalence-of-the-praisi-and-cs-functionals}
The functionals \eqref{eq:lecture-03:crisanti-sommers-functional} and
\eqref{eq:lecture-03:parisi-functinal-diagonal-case} are equivalent in the
following sense
\begin{align*}
\inf_{
\rho'
,
\lambda'
}
\mathcal{P}(
\rho'
,
\lambda'
)
=
\inf_{
\rho'
}
\mathcal{CS}(
\rho'
)
.
\end{align*}
\end{proposition}
\begin{proof}
\begin{enumerate}
\item
Let $(\rho, \lambda)$ be the solutions of equations
\eqref{eq:lecture-03:relation-between-crisanti-sommers-and-parisi-a-parts} and
\eqref{eq:lecture-03:spherical-lambda-optimiality-consequence}.
Lemma
\ref{lem:lecture-03:some-extremal-properties-of-the-crisanti-sommers-functional} 
guarantees that $\rho$ is the optimiser of the Crisnati-Sommers
functional and Lemma \ref{lem:some-extremal-properties-of-the-parisi-functional}
assures that $(\rho,\lambda)$ is the optimiser of the Parisi functional.
\item
We have
\begin{align}
\label{eq:lecture-03:parisi-minus-crisanti-sommers}
\mathcal{P}(
\rho
,
\lambda
)
-
\mathcal{CS}(
\rho
)
=
&
-
\lambda
u
+
2\beta^2
q^{(1)}
s^{(1)}
-
\frac{
q^{(1)}
}{
s^{(1)}
}
+
cu
-
1
\nonumber
\\
&
-
2\beta^2
\sum_{l=1}^{n}
x_l
\left(
[q^{(l+1)}]^2
-
[q^{(l)}]^2
\right)
.
\end{align}
We can simplify the $\Phi[B]$-like term (that is the summation) in
\eqref{eq:lecture-03:parisi-minus-crisanti-sommers}, using
\eqref{eq:lecture-03:relation-between-crisanti-sommers-and-parisi-a-parts} and
\eqref{eq:lecture-03:spherical-lambda-optimiality-consequence}. Indeed, 
\begin{align}
\label{eq:lecture-03:spherical-model-phi-b-terms-1}
2\beta^2
\sum_{l=1}^{n-1}
x_l
\left(
[q^{(l+1)}]^2
-
[q^{(l)}]^2
\right)
&
=
2\beta^2
\sum_{l=1}^{n-1}
x_l
\left(
q^{(l+1)}
[
q^{(l+1)}
-
q^{(l)}
]
+
q^{(l)}
[
q^{(l+1)}
-
q^{(l)}
]
\right)
\nonumber
\\
\text{
(by \eqref{eq:lecture-03:relation-between-crisanti-sommers-and-parisi-a-parts}
and \eqref{eq:lecture-03:spherical-quantity-s-k})
}
&
=
\sum_{l=1}^{n-1}
\left(
2\beta^2
q^{(l+1)}
\left[
s^{(l)}
-
s^{(l+1)}
\right]
+
q^{(l)}
\left[
\frac{1}{
s^{(l+1)}
}
-
\frac{1}{
s^{(l)}
}
\right]
\right)
.
\end{align}
Regrouping the summands in \eqref{eq:lecture-03:spherical-model-phi-b-terms-1},
we get
\begin{align}
\label{eq:lecture-03:spherical-model-phi-b-terms-2}
\text{
\eqref{eq:lecture-03:spherical-model-phi-b-terms-1}
}
=
&
2\beta^2
\sum_{l=1}^{n-1}
s^{(l)}
\left(
q^{(l+1)}
-
q^{(l)}
\right)
+
2\beta^2
\left(
q^{(1)}
s^{(1)}
-
q^{(n)}
s^{(n)}
\right)
\nonumber
\\
&
+
\sum_{l=1}^{n-1}
\frac{
q^{(l)}
-
q^{(l+1)}
}{
s^{(l+1)}
}
+
\left(
\frac{
q^{(n)}
}{
s^{(n)}
}
-
\frac{
q^{(1)}
}{
s^{(1)}
}
\right)
.
\end{align}
Due to
\eqref{eq:lecture-03:relation-between-crisanti-sommers-and-parisi-a-parts}, we 
have 
\begin{align*}
2\beta^2
\left(
q^{(l+1)}
-
q^{(l)}
\right)
=
\frac{
s^{(l)}
-
s^{(l+1)}
}{
x_l
s^{(l)}
s^{(l+1)}
}
=
\frac{
q^{(l+1)}
-
q^{(l)}
}{
s^{(l)}
s^{(l+1)}
}
.
\end{align*}
Applying the previous relation, we get that the
both summations in \eqref{eq:lecture-03:spherical-model-phi-b-terms-2} cancel out
and we end up with
\begin{align*}
\text{
\eqref{eq:lecture-03:spherical-model-phi-b-terms-2}
}
=
2\beta^2
\left(
q^{(1)}
s^{(1)}
-
q^{(n)}
s^{(n)}
\right)
+
\frac{
q^{(n)}
}{
s^{(n)}
}
-
\frac{
q^{(1)}
}{
s^{(1)}
}
.
\end{align*}
Now, turning back to \eqref{eq:lecture-03:parisi-minus-crisanti-sommers},
we get
\begin{align*}
\mathcal{P}(
\rho
,
\lambda
)
-
\mathcal{CS}(
\rho
)
&
=
-
\lambda
u
-
2\beta^2
\left(
u^2
-
[q^{(n)}]^2
\right)
+
2\beta^2
q^{(n)}
s^{(n)}
-
\frac{
q^{(n)}
}{
s^{(n)}
}
+
c u
-1
\\
\text{(by \eqref{eq:lecture-03:spherical-lambda-optimiality-consequence}) and \eqref{eq:lecture-03:spherical-quantity-s-k}}
&
=
-
u
\left(
c 
- 
2\beta^2 
(u-q^{(n)})
-
(u-q^{(n)})^{-1}
\right)
-
2\beta^2
\left(
u^2
-
[q^{(n)}]^2
\right)
\\
&
\quad
-
\frac{q^{(n)}}{
u-q^{(n)}
}
+
2\beta^2
q^{(n)}
\left(
u-q^{(n)}
\right)
+
c u
-
1
\\
&
=
0
.
\end{align*}

\end{enumerate}
\end{proof}

\subsection{Replica symmetric calculations}
In this subsection, we shall consider the one dimensional case of the a priori
measure \eqref{gaussian-spins:a-priori-measure} with $h = 0$. We shall also
restrict the computations to the case $n=1$ which is often referred to in
physical literature as the replica symmetric scenario. It is
indeed the right scenario under the above assumptions, as shows
Theorem~\ref{thm:gaussian:the-local-low-temperature-parisi-formula}.
\begin{lemma}
\label{lem:lecture-03:optimum-of-the-crisanti-sommers-functional}
Let $\mu$ satisfy \eqref{gaussian-spins:a-priori-measure} with $h = 0$.
Assume $d=1$, $n=1$ and $c>0$.
Given 
$
u \geq 0
$,
we have
\begin{align}
\label{eq:lecture-03:inf-spherical-cs}
\inf_{
\rho
\in
\mathcal{Q}(u,1)
}
\mathcal{CS}(\rho)
&
=
\inf_{
q \in [0;u]
}
\left(
1 
-
cu
+ 
\log 
\left(
c
(u-q)
\right)
+
\frac{q}{
u-q
}
+
\beta^2
\left(
u^2
-
q^2
\right)
\right)
=
f(c,u),
\end{align}
where $f(c,u)$ is defined in
\eqref{eq:lecture-03:inf-spherical-cs-two-clause-formula}.
\end{lemma}
\begin{proof}
Using the definitions, we obtain
\begin{align*}
\frac{
\partial
}{
\partial
q
}
\mathcal{CS}(\rho)
=
\frac{
\partial
}{
\partial
q
}
\left[
\log
\left(
u-q
\right)
+
\frac{q}{
u-q
}
+
\beta^2
\left(
u^2
-
q^2
\right)
\right]
=
\frac{q}{
(u-q)^2
}
-
2\beta^2 q
.
\end{align*}
Hence, the critical points of $q \mapsto \mathcal{CS}(q,u)$ are
\begin{align*}
q_0 
= 
0
,
q_{1,2} 
=
u
\pm
\frac{
\sqrt{2}
}{
2\beta
}
.
\end{align*}
Furthermore, we also have
\begin{align*}
\frac{
\partial^2
}{
\partial
q^2
}
\mathcal{CS}(q,u)
=
\frac{1}{
(u-q)^2
}
+
\frac{
2q
}{
(u-q)^3
}
-
2\beta^2
.
\end{align*}
Hence, as a simple calculation shows, the
infima in \eqref{eq:lecture-03:inf-spherical-cs} are attained on
\begin{align}
\label{eq:lecture-03:gaussian-spins:1d-optimizers}
q^*
=
\begin{cases}
0
, 
&
u \leq
\frac{\sqrt{2}}{2\beta}
,
\\
u
-
\frac{
\sqrt{2}
}{
2\beta
}
,
&
u > \frac{\sqrt{2}}{2\beta}
\end{cases}
\end{align}
which implies \eqref{eq:lecture-03:inf-spherical-cs}.
\end{proof}
\begin{lemma}
\label{lem:lecture-03:gaussian-spins:saddle-point-of-the-crisanti-sommers-functional}
Under the assumptions of
Lemma~\ref{lem:lecture-03:optimum-of-the-crisanti-sommers-functional}, we have
\begin{enumerate}
\item 
For $
c
\geq
2\sqrt{2}\beta
$, 
we have
\begin{align*}
\sup_{
u \geq 0
}
\inf_{
q \in [0;u]
}
\mathcal{CS}(q,u)
=
\mathcal{CS}(0,u^*)
=
\beta^2 (u^*)^2
+
\log c u^*
-
c u^*
+
1
,
\end{align*}
where
\begin{align*}
u^*
\equiv
\frac{1}{4\beta^2}
\left(
c
-
\sqrt{
c^2
-
8\beta^2
}
\right)
.
\end{align*}
\item
For $
c
<
2\sqrt{2}\beta
$, 
we have
\begin{align*}
\sup_{
u \geq 0
}
\inf_{
q \in [0;u]
}
\mathcal{CS}(q,u)
=
+\infty
.
\end{align*}
\end{enumerate}

\end{lemma}
\begin{remark}
Under the assumptions, the above theorem says that from the
point of view of the global free energy, the system can only exist in the
``high temperature'' scenario, cf.
\eqref{eq:lecture-03:inf-spherical-cs-two-clause-formula}.
The threshold at 
$
c_0
=
2\sqrt{2}\beta
$
could be easily understood from the perspective of the norms of
random matrices.
\end{remark}
\begin{proof}
\begin{enumerate}
\item 
Suppose 
$
c
\geq
2\sqrt{2}\beta
$.
Recalling \eqref{eq:lecture-03:inf-spherical-cs-two-clause-formula}, for
$
u 
\in 
(0;\frac{\sqrt{2}}{2\beta}]
$,
we
introduce the following function
\begin{align*}
f(u)
\equiv
\log(c u)
+
\beta^2 u^2
-
c u
+
1
.
\end{align*}
We have
\begin{align*}
\frac{
\partial
}{
\partial 
u
}
f(u)
=
\frac{1}{u}
+
2\beta^2 u
-
c
.
\end{align*}
Hence, the critical points of the function $f$ are
\begin{align*}
u_{1,2}
=
\frac{
c
\pm
\sqrt{
c^2
-
8\beta^2
}
}{
4\beta^2
}
.
\end{align*}
Furthermore, we have
\begin{align*}
\frac{
\partial^2
}{
\partial 
u^2
}
f(u)
=
2\beta^2
-
\frac{1}{u^2}
.
\end{align*}
We notice that
$
u^* 
\leq
\frac{\sqrt{2}}{
2\beta
}
$
and, hence, due to \eqref{eq:lecture-03:inf-spherical-cs-two-clause-formula}
\begin{align*}
\mathcal{CS}(0,u^*)
=
\beta^2 (u^*)^2
+
\log c u^*
-
c u^*
+
1
.
\end{align*}
\item
If $c < 2\sqrt{2}\beta$, then the function
\begin{align*}
u 
\mapsto 
(
2\sqrt{2} \beta 
-
c
)
u
+
\log 
\frac{c}{\beta}
-
\frac{1}{2}
\left(
1+\log 2
\right)
\end{align*}
is unbounded on 
$
(
\frac{\sqrt{2}}{2\beta}
; 
+\infty
)
$.

\end{enumerate}
\end{proof}
\subsection{The multidimensional Crisanti-Sommers functional}
Recall the definition \eqref{eq:as2-psi-of-t-definition}.
\begin{proposition}
Assume $d = 1$.
Given $u > 0$, we have
\begin{align}
\label{gaussian-spins-value-of-phi-t}
2\phi^{(x*,\mathcal{Q}*,\Lambda*)}(t)
=
\begin{cases}
\left(
3\sqrt{2}
\beta
-
c
\right)
u
+
\log 
\frac{c}{\beta}
-
1
-
\frac{
\log 2
}{
2
}
-
t
\left(
\sqrt{2} u \beta
-
\frac{1}{2}
\right)
,
&
u
>
\frac{
\sqrt{2}
}{
2\beta
}
,
\\
2 \beta^2 (u)^2
+
\log(c u)
-
c u
+
1
-
t \beta^2 (u)^2
,
&
u
\leq
\frac{
\sqrt{2}
}{
2\beta
}
.
\end{cases}
\end{align}
\end{proposition}
\begin{proof}
Combining \eqref{eq:lecture-03:parisi-functinal-diagonal-case},
\eqref{eq:lecture-03:crisanti-sommers-functional} with Lemma
\ref{lem:lecture-03:optimum-of-the-crisanti-sommers-functional} and Proposition
\ref{thm:lecutre-03:equivalence-of-the-praisi-and-cs-functionals}, we get the claim.

\end{proof}

\subsection{Talagrand's a priori estimates}
In this subsection, we prove that Assumption~\ref{as2-an-a-priori-estimate}
is satisfied in the case of the Gaussian a priori distribution
\eqref{gaussian-spins:a-priori-measure} with $h=0$.
\begin{theorem}
\label{thm:gaussian-spins:a-priori-estimate}
Let $\mu$ satisfy \eqref{gaussian-spins:a-priori-measure} with $h = 0$,
assume $U \in \symmetric^+(d)$ is such that 
$
\min_v u_v >
\frac{\sqrt{2}}{2\beta}
$ and suppose $C \succ 0$.
Let $Q = Q^*$ and $\Lambda = \Lambda^*$. 

Then, for any $t_0 \in (0;1)$ and any
$t \in (0;t_0]$, we have 
(cf. \eqref{eq:as2-a-priori-estimate} with $k=1$)
\begin{align}
\label{eq:lecture-03:gaussian-spins:a-priori-estimate}
\varphi^{(2)}_N
(
1
,
t
,
x,Q,
\Sigma^{(2)}_N(\mathfrak{L},\mathfrak{U},\eps,\delta)
)
\leq
2 \phi^{(x,\mathcal{Q},\Lambda)}(t)
-
\frac{1}{K}
\Vert
Q^{(1)}
-
V
\Vert_{\text{F}}^2
+
\mathcal{O}(\eps+\delta)
.
\end{align}
\end{theorem}
\begin{proof}
\begin{enumerate}
\item 
We employ the notations of
Section~\ref{upper-bounds-on-phi-2-guerras-scheme-revisited}. Let
$\mathfrak{n}=1$. Given $\mathfrak{U} \in \symmetric(2d)$ (cf.
\eqref{eq:as2-mathfrak-u-definition}), choose arbitrary matrices $
\left\{
\mathfrak{Q}^{(l)}
\in
\symmetric (2d)
\mid
l 
\in 
[0;2] \cap \N
\right\}
$ 
satisfying \eqref{eq:remainder:remainder-comparison-matrices-monotonicity}. Define 
$\mathfrak{x} \equiv x$ which, in particular, implies that $\zeta = \xi$. Finally,
we set, for $l \in [0;n+1] \cap \N$, 
$
\widetilde{Q}^{(l)}
\equiv
Q^{(l)}
$.
\item
Proposition~\ref{eq:remainder:phi-2-upper-bound} implies that, for any $\delta$-minimal
$
\mathfrak{L} 
\in 
\R^{2d \times 2d}
$, we have
\begin{align}
\label{eq:lecture-03:gaussian-spins:phi-2-upper-bound}
\varphi^{(2)}_N
(
1
,
t
,
x,Q,
\Sigma^{(2)}_N(\mathfrak{L},\mathfrak{U},\eps,\delta)
)
&
\leq
-
\langle
\mathfrak{L}
,
\mathfrak{U}
\rangle
-
\frac{
t \beta^2
}{2}
\left(
\Vert
\mathfrak{Q}^{(2)}
\Vert_{\text{F}}^2
-
\Vert
\mathfrak{Q}^{(1)}
\Vert_{\text{F}}^2
\right)
\nonumber
\\
&
\quad
+
X^{(2)}_{0}
(
1
,
\mathfrak{x}
,
\widehat{\mathfrak{Q}}^{(l)}(t)
,
\mathfrak{L}
)
+
\mathcal{O}
(
\eps+\delta
)
.
\end{align}
\item
We define a matrix $\mathfrak{C} \in \R^{2d \times 2d}$ as follows
\begin{align*}
\mathfrak{C}
\equiv
\begin{bmatrix}
C
& 
0
\\
0
&
C
\end{bmatrix}
.
\end{align*}
Recalling \eqref{eq:lecture-03:spherical-d-matrices-of-the-a-part}, we define
also the following matrices 
$
\mathfrak{D}^{(2)}
\equiv
\mathfrak{C}
$
and
\begin{align}
\label{eq:gaussian:mathfrak-d-1-definition}
\mathfrak{D}^{(1)}
\equiv
\mathfrak{C}
-
\mathfrak{L}
-
\left(
\widehat{\mathfrak{Q}}^{(2)}(t)
-
\widehat{\mathfrak{Q}}^{(1)}(t)
\right)
.
\end{align}
Applying Proposition \ref{prp:gaussian:guerras-scheme} to
\eqref{eq:lecture-03:gaussian-spins:phi-2-upper-bound}, we get 
\begin{align}
\label{eq:lecture-03:gaussian-spins:varphi-2-upper-bound-2}
\varphi^{(2)}_N
(
1
,
t
,
\Sigma^{(2)}_N(\mathfrak{L},\mathfrak{U},\eps,\delta)
)
&
\leq
\frac{1}{2}
\left[
-
\langle
\mathfrak{L}
,
\mathfrak{U}
\rangle
-
t \beta^2
\left(
\Vert
\mathfrak{Q}^{(2)}
\Vert_{\text{F}}^2
-
\Vert
\mathfrak{Q}^{(1)}
\Vert_{\text{F}}^2
\right)
\right.
\nonumber
\\
&
\left.
\quad
\quad
+
2\beta^2
\langle
[
\mathfrak{D}^{(1)}
]^{-1}
,
\widehat{\mathfrak{Q}}^{(1)}(t)
\rangle
+
\log
\left(
\frac{
\det 
\mathfrak{D}^{(2)}
}{
\det 
\mathfrak{D}^{(1)}
}
\right)
\right]
+
\mathcal{O}(\eps)
\nonumber
\\
&
=:
\widetilde{\Phi}^{(2),k,\mathfrak{x},\mathfrak{L}}
+
\mathcal{O}(\eps)
.
\end{align}
\item
Assume that the matrices 
\begin{align}
\label{eq:lecture-03:gaussian-spins:a-list-of-diagonilizable-matrices}
\mathfrak{Q}^{(1)}
,
\mathfrak{Q}^{(2)}
,
\mathfrak{D}^{(1)}
\in
\R^{2 d \times 2 d}
\end{align}
are simultaneously diagonalisable in the same basis which is given by
the orthogonal
matrix $\mathfrak{O} \in \R^{2d \times 2d}$. Let the vectors
\begin{align}
\label{eq:lecture-03:gaussian-spins:a-list-of-spectra}
\mathfrak{q}^{(1)}
,
\mathfrak{q}^{(2)}
,
\mathfrak{d}^{(1)}
\in
\R^{2 d}
\end{align}
be the corresponding spectra of the matrices
\eqref{eq:lecture-03:gaussian-spins:a-list-of-diagonilizable-matrices}. That
is, we assume that
\begin{align*}
\mathfrak{Q}^{(1)}
=
\mathfrak{O}^*
\diag
\mathfrak{q}^{(1)}
\mathfrak{O}
,
\mathfrak{Q}^{(2)}
=
\mathfrak{O}^*
\diag
\mathfrak{q}^{(2)}
\mathfrak{O}
,
\\
\mathfrak{D}^{(1)}
=
\mathfrak{O}^*
\mathfrak{d}^{(1)}
\mathfrak{O}
,
\widetilde{\mathfrak{Q}}^{(1)}
=
\mathfrak{O}^*
\diag
\widetilde{\mathfrak{Q}}^{\prime(1)}
\mathfrak{O}
,
\end{align*}
where we have introduced the matrix 
$
\widetilde{\mathfrak{Q}}^{\prime(1)}(t)
\in
\symmetric^+(2d)
$.
By \eqref{eq:lecture-03:gaussian-spins:1d-optimizers}, we have, 
$
Q^{(2)}-Q^{(1)}
=
\frac{\sqrt{2}}{
2\beta
}
I
$, where $I$ denotes the unit matrix of the suitable dimension.
The definitions
\eqref{eq:as2:gaussian-spins:tilide-frak-q-definition-1} and 
\eqref{eq:as2:gaussian-spins:tilide-frak-q-definition-2}
then imply
\begin{align}
\label{eq:gaussian:delta-q-tilde}
\widetilde{\mathfrak{Q}}^{(2)}
-
\widetilde{\mathfrak{Q}}^{(1)}
=
\frac{\sqrt{2}}{
2\beta
}
I
.
\end{align}
Using the definitions and the above relation, we obtain
\begin{align}
\label{eq:gaussian:diagonal-representation-of-the-main-matrices}
\widehat{\mathfrak{Q}}^{(1)}_v(t)
&
=
\mathfrak{O}^*
\Bigl(
t
\diag
\mathfrak{q}^{(1)}
+
(1-t)
\widetilde{\mathfrak{Q}}^{\prime(1)}
\Bigr)
\mathfrak{O}
,
\nonumber
\\
\widehat{\mathfrak{Q}}^{(2)}(t)
-
\widehat{\mathfrak{Q}}^{(1)}(t)
&
=
\mathfrak{O}^*
\Bigl(
t
\diag
(
\mathfrak{q}^{(2)}
-
\mathfrak{q}^{(1)}
)
+
(1-t)
\frac{\sqrt{2}}{
2\beta
}
I
\Bigr)
\mathfrak{O}
.
\end{align}
Motivated by \eqref{eq:lecture-03:gausian-spins:s-k-equals-inverse-d-k}, we
set 
\begin{align}
\label{eq:gaussian:mathfrak-d-1-choice}
\mathfrak{d}^{(1)}_v
\equiv
\left(
\mathfrak{u}_v
-
\mathfrak{q}^{(1)}_v
\right)^{-1}
.
\end{align}
In view of \eqref{eq:gaussian:mathfrak-d-1-definition}, the above choice
necessarily yields (cf.
\eqref{eq:lecture-03:spherical-lambda-optimiality-consequence})
\begin{align}
\label{eq:gaussian:mathfrak-l-formula}
\mathfrak{L}
&
=
\mathfrak{C}
-
\mathfrak{O}^*
\diag
(
\mathfrak{u}_v
-
\mathfrak{q}^{(1)}_v
)^{-1}
\mathfrak{O}
-
\Bigl(
\widehat{\mathfrak{Q}}^{(2)}(t)
-
\widehat{\mathfrak{Q}}^{(1)}(t)
\Bigr)
\nonumber
\\
&
=
\mathfrak{C}
-
\mathfrak{O}^*
\Bigl(
\diag
(
\mathfrak{u}_v
-
\mathfrak{q}^{(1)}_v
)^{-1}
+
t
\diag
(
\mathfrak{q}^{(2)}
-
\mathfrak{q}^{(1)}
)
+
(1-t)
\frac{\sqrt{2}}{
2\beta
}
I
\Bigr)
\mathfrak{O}
.
\end{align}
Applying
Lemma~\ref{lem:lecture-03:gaussian-spins:a-decoupling-of-coordinates-for-the-a-and-b-terms-in-the-diagonal-case}
to
\eqref{eq:lecture-03:gaussian-spins:varphi-2-upper-bound-2} and using
\eqref{eq:gaussian:mathfrak-l-formula},
\eqref{eq:gaussian:mathfrak-d-1-choice},
\eqref{eq:gaussian:diagonal-representation-of-the-main-matrices},
we get the following diagonalised representation of
\eqref{eq:lecture-03:gaussian-spins:phi-2-upper-bound}
\begin{align}
\label{eq:gaussian:first-diagonalised-bound-on-phi-2}
\varphi^{(2)}_N
(
1
,
t
,
x,Q,
\Sigma^{(2)}_N(\mathfrak{L},\mathfrak{U},\eps,\delta)
)
\leq
&
\frac{1}{2}
\log\det\mathfrak{C}
-
\frac{1}{2}
\langle
\mathfrak{C}
,
\mathfrak{U}
\rangle
\nonumber
\\
&
+
\frac{1}{2}
\sum_{
v=1
}^{2 d}
\Bigl\{
\mathfrak{u}_v
\Bigl[
(\mathfrak{u}_v
-
\mathfrak{q}^{(1)}_v
)^{-1}
+
2\beta^2
\Bigl(
t
(
\mathfrak{q}^{(2)}_v
-
\mathfrak{q}^{(1)}_v
)
+
(1-t)
\frac{\sqrt{2}}{
2\beta
}
\Bigr)
\Bigr]
\nonumber
\\
&
\quad\quad
+
2\beta^2
(
\mathfrak{u}_v
-
\mathfrak{q}^{(1)}_v
)
\left(
t
\mathfrak{q}^{(1)}_v
+
(1-t)
\widetilde{\mathfrak{q}}^{(1)}_v
\right)
+
\log
(
\mathfrak{u}_v
-
\widetilde{\mathfrak{Q}}^{\prime(1)}_{v,v}
)
\nonumber
\\
&
\quad\quad
-
t \beta^2
\Bigl(
(
\mathfrak{q}^{(2)}_v
)^2
-
(
\mathfrak{q}^{(1)}_v
)^2
\Bigr)
\Bigr\}
+
\mathcal{O}(\eps)
.
\end{align}
Using the definitions, we get 
\begin{align}
\label{eq:gaussian:c-part-tracial-relations}
\langle
\mathfrak{C}
,
\mathfrak{U}
\rangle
&
=
2
\langle
C
,
U
\rangle
=
2
\sum_{v=1}^d
c_v
u_{v}
,
\nonumber
\\
\log\det\mathfrak{C}
&
=
2
\log \det C
=
2
\sum_{v=1}^d
\log
c_v
.
\end{align}
Motivated by \eqref{eq:gaussian:delta-q-tilde} (or by
\eqref{eq:lecture-03:gaussian-spins:1d-optimizers}), 
we define
\begin{align}
\label{eq:gaussian:q-1-assumption}
\mathfrak{q}^{(1)}_v
:=
\mathfrak{u}_v
-
\frac{
\sqrt{2}
}{
2\beta
}
.
\end{align}
In this case, as a straightforward calculation shows, the expression in the
curly brackets in \eqref{eq:gaussian:first-diagonalised-bound-on-phi-2} equals
\begin{align}
\label{eq:gaussian:simplified-curly-brackets}
2\sqrt{2} \beta \mathfrak{u}_v
+
\beta\sqrt{2}
\widetilde{\mathfrak{Q}}^{\prime(1)}_{v,v}
(1-t)
-
\log \beta
-
\frac{1}{2}
(
\log 2
-
t
)
.
\end{align}
By the definitions and the general properties of matrix trace, we have
\begin{align}
\label{eq:gaussian:tracial-relations}
\sum_{v=1}^{2d}
\widetilde{\mathfrak{Q}}^{\prime(1)}_{v,v}
&
=
\sum_{v=1}^{2d}
\widetilde{\mathfrak{Q}}^{(1)}_{v,v}
=
2
\sum_{v=1}^d
Q^{(1)}_{v,v}
,
\nonumber
\\
\sum_{v=1}^{2d}
\mathfrak{u}_v
&
=
2
\sum_{v=1}^{d}
U_{v,v}
.
\end{align}
Combining
\eqref{eq:gaussian:first-diagonalised-bound-on-phi-2} with 
\eqref{eq:gaussian:simplified-curly-brackets},
\eqref{eq:gaussian:tracial-relations} and
\eqref{eq:gaussian:c-part-tracial-relations},
we obtain
\begin{align}
\label{eq:gaussian:non-sharp-a-priori-bound}
\varphi^{(2)}_N
(
1
,
t
,
x,Q,
\Sigma^{(2)}_N(\mathfrak{L},\mathfrak{U},\eps,\delta)
)
&
\leq
\sum_{v=1}^{d}
\Bigl(
-
c_v u_{v}
+
\log c_v
+
3\sqrt{2}u\beta
\nonumber
\\
&
\quad\quad
-
\frac{1}{2}
(
\log 2
-
t
)
-
\sqrt{2} \beta t u
-
\log\beta
-
1
\Bigr)
+
\mathcal{O}(\eps)
\nonumber
\\
&
=
2
\sum_{v=1}^{d}
\phi(t)\vert_{
\substack{
c=c_v,
\\
u=u_v
}
}
+
\mathcal{O}(\eps)
,
\end{align}
where in the last line we have used the relation
\eqref{gaussian-spins-value-of-phi-t}.
\item
To get the version of the a priori
bound \eqref{eq:gaussian:non-sharp-a-priori-bound} with the quadratic correction
term as stated in \eqref{eq:lecture-03:gaussian-spins:a-priori-estimate}, we
perturb the r.h.s of \eqref{eq:lecture-03:gaussian-spins:phi-2-upper-bound}
around our choice of $\mathfrak{D}^{(1)}$ in
\eqref{eq:gaussian:mathfrak-d-1-choice}, i.e.,
\begin{align*}
\mathfrak{D}^{(1)} 
=
\left(
\mathfrak{U}_v
-
\mathfrak{Q}^{(1)}_v
\right)^{-1}
=
\sqrt{2}\beta
I
,
\end{align*}
where in the last equality we used \eqref{eq:gaussian:q-1-assumption}.
\end{enumerate}
\end{proof}

\subsection{The local low temperature Parisi formula}

\begin{proof}[Proof of
Theorem~\ref{thm:gaussian:the-local-low-temperature-parisi-formula}]
The result follows from Theorem~\ref{thm:gaussian-spins:a-priori-estimate} and
Theorem~\ref{thm:remainder:the-conditional-parisi-formula}. Note that the proof
of Theorem~\ref{thm:remainder:the-conditional-parisi-formula} requires a minor modification to cope with the fact that the a priori distribution
\eqref{gaussian-spins:a-priori-measure} is unbounded. This minor problem
can be fixed by considering the pruned Gaussian distribution and using the
elementary estimates to bound the tiny Gaussian tails.
\end{proof}

\appendix
\section{}
\label{sec:thermodynamic-limit}
The general result of Guerra and Toninelli
\cite{Guerra-Toninelli-Generalized-SK-2003} implies that the thermodynamic
limit of the local free energy \eqref{eq:as2:local-limiting-free-energy} exists
almost surely and in $L^1$. The following existence of the limiting average overlap is an immediate consequence of this.
\begin{proposition}
We have
\begin{align*}
\E
\left[
\mathcal{G}_N(\beta) \otimes \mathcal{G}_N(\beta)
\left[
\var H_N(\sigma)
-
\E
\left[
H_N(\sigma)
H_N(\sigma^\prime)
\right]
\right]
\right]
\xrightarrow[N\uparrow+\infty]{}
C(\beta)
\geq 0
,
\end{align*}
where $C: \R_+ \to \R_+$.
\end{proposition}
\begin{proof}
The free energy is a convex function of $\beta$ (a consequence of the Hölder
inequality). Hence, by a result in \cite{Griffiths1964} the following holds
\begin{align*}
\lim_{N\uparrow\infty}
\frac{\dd}{\dd \beta}
\E
\left[
p_N(\beta)
\right] 
=
\frac{\dd}{\dd \beta}
\E
\left[
p(\beta)
\right] 
.
\end{align*}
Proposition \ref{free-energy-beta-derivative} implies
\begin{align*}
\frac{\dd}{\dd \beta} 
\E
\left[
p_N(\beta)
\right]
=
\beta
\E
\left[
\mathcal{G}_N(\beta) \otimes \mathcal{G}_N(\beta)
\left[
\var H_N(\sigma)
-
\E
\left[
H_N(\sigma)
H_N(\sigma^\prime)
\right]
\right]
\right]
.
\end{align*}
\end{proof}
The following super-additivity result is an application of the Gaussian
comparison inequalities obtained in
Subection~\ref{sec:gaussian-comparison-inequalities-for-free-energy-like-functionals}.
Note that the result does not provide enough information for the
cavity-like argument of \cite{AizenmanSimsStarr2003}.
\begin{proposition}
\label{lem:as2:local-almost-subadditivity}
For any 
$
\mathcal{V} 
\equiv
B(U,\eps)
\subset
\mathcal{U}
$,
we have
\begin{align*}
N 
\E
\left[
p_{N}(\mathcal{V})
\right]
+
M 
\E
\left[
p_{M}(\mathcal{V})
\right]
\leq
(N+M)
\E
\left[
p_{N+M}(\mathcal{V})
\right]
+
(N+M)\mathcal{O}(\eps)
,
\end{align*}
as 
$
\eps \downarrow +0
$.
\end{proposition}
\begin{proof}
Define the process 
$
Y_{N,M}
\equiv
\{
Y(\sigma)
:
\sigma
=
\alpha
\shortparallel
\tau
;
\alpha \in \Sigma_N
,
\tau \in \Sigma_M
\}
$
as follows
\begin{align*}
Y(
\alpha
\shortparallel
\tau
)
\equiv
\left(
\frac{N}{N+M}
\right)^{1/2}
X^{(1)}_N(\alpha)
+
\left(
\frac{M}{N+M}
\right)^{1/2}
X^{(2)}_M(\tau)
,
\end{align*}
where $X^{(1)}$ and $X^{(2)}$ are two independent copies of the process $X$.
Given some Gaussian process
$
\{
C(\sigma)
\}_{
\sigma \in \Sigma_N
}
$,
let us introduce the functional 
$
\Phi_N(\beta)[C]
$
as follows
\begin{align*}
\Phi_{N,M}(\beta)[C]
\equiv
\E
\left[
\log
\mu^{\otimes (N+M)}
\left[
\I_{
\Sigma_N(\mathcal{V})
}
\I_{
\Sigma_M(\mathcal{V})
}
\exp
(
\beta
\sqrt{N+M}
C
)
\right]
\right]
.
\end{align*}
Now, set
$
\varphi(t)
\equiv
\Phi_{N+M}(\beta)
\left[
\sqrt{t}
X_{N+M}
+
\sqrt{1-t}
Y_{N,M}
\right]
$
.
Applying Proposition~\ref{prp:as2:gaussian-comparison-of-free-energy}, we get
\begin{align}
\label{eq:as2:local-subadditivity-interpolation-speed}
\frac{\dd}{\dd t}
\varphi(t)
&=
\frac{
\beta^2
(N+M)
}{
2
}
\E
\left[
\mathcal{G}(t)
\otimes 
\mathcal{G}(t)
\left[
\right.
\right.
\nonumber
\\
&
\quad
\left(
\var X_{N+M}(\sigma^{(1)}) 
- 
\var Y_{N,M}(\sigma^{(1)}) 
\right)
\nonumber
\\
&
\quad
\left.
\left.
-
\left(
\cov
\left[
X_{N+M}(\sigma^{(1)})
,
X_{N+M}(\sigma^{(2)})
\right]
- 
\cov 
\left[
Y_{N,M}(\sigma^{(1)})
,
Y_{N,M}(\sigma^{(2)})
\right]
\right)
\right]
\right]
.
\end{align}
Note that we have 
\begin{align}
\label{eq:as2:local-subadditivity-interpolation-start-end}
\varphi(0) 
&
= 
N 
\E
\left[
p_N(\mathcal{V})
\right] 
+ 
M 
\E
\left[
p_M(\mathcal{V})
\right]
,
\nonumber
\\
\varphi(1) 
&
\leq
(N+M)
\E
\left[
p_{N+M}(\mathcal{V})
\right] 
,
\end{align}
where the last inequality is due to the fact that, for all 
$
\alpha 
\in
\Sigma_N(\mathcal{V})
$ 
and all 
$
\tau
\in
\Sigma_N(\mathcal{V})
$,
we have
\begin{align*}
\alpha \shortparallel \tau 
\in 
\Sigma_{N+M}(\mathcal{V})
.
\end{align*}
Moreover, for 
$
\sigma = \alpha \shortparallel \tau
$
with 
$
\alpha \in
\Sigma_N(\mathcal{V})
$
and 
$
\sigma \in
\Sigma_M(\mathcal{V})
$
we have
\begin{align*}
\var X_{N+M}(\sigma) 
- 
\var Y_{N,M}(\sigma) 
&
=
\left\Vert
\frac{N}{N+M}
R_{N}(
\alpha
,
\alpha
)
+
\frac{M}{N+M}
R_{M}(
\tau
,
\tau
)
\right\Vert_2^2
-
\frac{N}{N+M}
\Vert
R_{N}(
\alpha
,
\alpha
)
\Vert_2^2
\\
&
\quad
-
\frac{M}{N+M}
\Vert
R_{M}(
\tau
,
\tau
)
\Vert_2^2
=
\mathcal{O}(\eps)
.
\end{align*}
Also, due to convexity of the norm, we have
\begin{align*}
&
\cov
\left[
X_{N+M}(\sigma^{(1)})
,
X_{N+M}(\sigma^{(2)})
\right]
- 
\cov 
\left[
Y_{N,M}(\sigma^{(1)})
,
Y_{N,M}(\sigma^{(2)})
\right]
\\
&
=
\left\Vert
\frac{N}{N+M}
R_{N}(
\alpha^{(1)}
,
\alpha^{(2)}
)
+
\frac{M}{N+M}
R_{M}(
\tau^{(1)}
,
\tau^{(2)}
)
\right\Vert_2^2
-
\frac{N}{N+M}
\Vert
R_{N}(
\alpha^{(1)}
,
\alpha^{(2)}
)
\Vert_2^2
\\
&
\quad
-
\frac{M}{N+M}
\Vert
R_{M}(
\tau^{(1)}
,
\tau^{(2)}
)
\Vert_2^2
\leq
0
.
\end{align*}
Applying 
$
\int_0^1 \dd t
$
to 
\eqref{eq:as2:local-subadditivity-interpolation-speed} and using
the previous two formulae, we get the claim.
\end{proof}

\bibliographystyle{plain}
\bibliography{bibliography}

\begin{thebibliography}{10}

\bibitem{AizenmanSimsStarr2003}
Michael Aizenman, Robert Sims, and Shannon~L. Starr.
\newblock An {E}xtended {V}ariational {P}rinciple for the {SK} {S}pin-{G}lass
  {M}odel.
\newblock {\em Phys. Rev. B}, 68:214403, 2003.

\bibitem{AizenmanSimsStarr2006}
Michael Aizenman, Robert Sims, and Shannon~L. Starr.
\newblock Mean-field spin glass models from the cavity-{ROS}t perspective.
\newblock In {\em Prospects in mathematical physics}, volume 437 of {\em
  Contemp. Math.}, pages 1--30. Amer. Math. Soc., Providence, RI, 2007.
\newblock arXiv:math-ph/0607060.

\bibitem{Arguin2006}
Louis-Pierre Arguin.
\newblock Spin glass computations and {R}uelle's probability cascades.
\newblock {\em J. Stat. Phys.}, 126(4-5):951--976, 2007.
\newblock arXiv:math-ph/0608045v1.

\bibitem{BenArousDemboGuionnet2001}
G{\'e}rard {Ben Arous}, Amir Dembo, and Alice Guionnet.
\newblock Aging of spherical spin glasses.
\newblock {\em Probab. Theory Related Fields}, 120(1):1--67, 2001.

\bibitem{Bogachev1998}
Vladimir~I. Bogachev.
\newblock {\em Gaussian measures}, volume~62 of {\em Mathematical Surveys and
  Monographs}.
\newblock American Mathematical Society, Providence, RI, 1998.

\bibitem{BolthausenSznitman1998}
Erwin Bolthausen and Alain-Sol Sznitman.
\newblock On {R}uelle's probability cascades and an abstract cavity method.
\newblock {\em Comm. Math. Phys.}, 197(2):247--276, 1998.

\bibitem{BovierBook2006}
Anton Bovier.
\newblock {\em Statistical mechanics of disordered systems}.
\newblock Cambridge Series in Statistical and Probabilistic Mathematics.
  Cambridge University Press, Cambridge, 2006.
\newblock A mathematical perspective.

\bibitem{BovierKurkova2007}
Anton Bovier and Irina Kurkova.
\newblock Much ado about {D}errida's {GREM}.
\newblock In {\em Spin glasses}, volume 1900 of {\em Lecture Notes in Math.},
  pages 81--115. Springer, Berlin, 2007.

\bibitem{Briand2007}
Philippe Briand and Ying Hu.
\newblock {Quadratic BSDEs with convex generators and unbounded terminal
  conditions}.
\newblock {\em Probab. Theory Related Fields}, pages~--, 2007.
\newblock arXiv:math/0703423.

\bibitem{Comets1989}
Francis Comets.
\newblock Large deviation estimates for a conditional probability distribution.
  {A}pplications to random interaction {G}ibbs measures.
\newblock {\em Probab. Theory Related Fields}, 80(3):407--432, 1989.

\bibitem{CrisantiSommers1992}
Andrea Crisanti and Hans-J{\"u}rgen Sommers.
\newblock {The spherical p-spin interaction spin glass model: the statics}.
\newblock {\em Zeitschrift f{\"u}r Physik B Condensed Matter}, 87(3):341--354,
  1992.

\bibitem{DaLioLey2006}
Francesca {Da~Lio} and Olivier Ley.
\newblock {Uniqueness Results for Second-Order Bellman--Isaacs Equations under
  Quadratic Growth Assumptions and Applications}.
\newblock {\em SIAM J. Control Optim.}, 45(1):74--106, 2006.

\bibitem{denHollanderLDPBook2000}
Frank den Hollander.
\newblock {\em Large deviations}, volume~14 of {\em Fields Institute
  Monographs}.
\newblock American Mathematical Society, Providence, RI, 2000.

\bibitem{Evans1998}
Lawrence~C. Evans.
\newblock {\em Partial differential equations}, volume~19 of {\em Graduate
  Studies in Mathematics}.
\newblock American Mathematical Society, Providence, RI, 1998.

\bibitem{Froehlich-Zegarlinski-Generalized-SK-1987}
J{\"u}rg Fr{\"o}hlich and Bogus{\l}aw Zegarli{\'n}ski.
\newblock Some comments on the {S}herrington-{K}irkpatrick model of spin
  glasses.
\newblock {\em Comm. Math. Phys.}, 112(4):553--566, 1987.

\bibitem{Griffiths1964}
Robert~B. Griffiths.
\newblock A proof that the free energy of a spin system is extensive.
\newblock {\em J. Mathematical Phys.}, 5:1215--1222, 1964.

\bibitem{Guerra2003a}
Francesco Guerra.
\newblock Broken replica symmetry bounds in the mean field spin glass model.
\newblock {\em Comm. Math. Phys.}, 233(1):1--12, 2003.

\bibitem{GuerraReview2005}
Francesco Guerra.
\newblock Spin glasses.
\newblock Preprint, 2005.
\newblock arXiv:cond-mat/0507581v1.

\bibitem{Guerra-Toninelli-Generalized-SK-2003}
Francesco Guerra and Fabio~Lucio Toninelli.
\newblock The infinite volume limit in generalized mean field disordered
  models.
\newblock {\em Markov Process. Related Fields}, 9(2):195--207, 2003.

\bibitem{PanchenkoParisiMeasures2004}
Dmitry Panchenko.
\newblock {A question about the Parisi functional.}
\newblock {\em Electron. Commun. Probab.}, 10:155--166, 2005.

\bibitem{panchenko-free-energy-generalized-sk-2005}
Dmitry Panchenko.
\newblock Free energy in the generalized {S}herrington-{K}irkpatrick mean field
  model.
\newblock {\em Rev. Math. Phys.}, 17(7):793--857, 2005.

\bibitem{Panchenko2007}
Dmitry Panchenko and Michel Talagrand.
\newblock {Guerra's interpolation using Derrida-Ruelle cascades}.
\newblock {\em Preprint}, 2007.
\newblock arXiv:0708.3641v2 [math.PR].

\bibitem{PanchenkoTalagrandMultipleSKModels2006}
Dmitry Panchenko and Michel Talagrand.
\newblock {On the overlap in the multiple spherical SK models}.
\newblock {\em Ann. Probab.}, 35(6):2321--2355, 2007.

\bibitem{Ruelle1987}
David Ruelle.
\newblock A mathematical reformulation of {D}errida's {REM} and {GREM}.
\newblock {\em Comm. Math. Phys.}, 108(2):225--239, 1987.

\bibitem{Sherrington2007}
David Sherrington.
\newblock Spin glasses: a perspective.
\newblock In {\em Spin glasses}, volume 1900 of {\em Lecture Notes in Math.},
  pages 45--62. Springer, Berlin, 2007.

\bibitem{SherringtonKirkpatrick1975}
David Sherrington and Scott Kirkpatrick.
\newblock {Solvable Model of a Spin-Glass}.
\newblock {\em Physical Review Letters}, 35(26):1792--1796, 1975.

\bibitem{Talagrand2000}
Michel Talagrand.
\newblock Large deviation principles and generalized
  {S}herrington-{K}irkpatrick models.
\newblock {\em Ann. Fac. Sci. Toulouse Math. (6)}, 9(2):203--244, 2000.

\bibitem{TalagrandSpinGlassesBook2003}
Michel Talagrand.
\newblock {\em Spin glasses: a challenge for mathematicians}, volume~46 of {\em
  Ergebnisse der Mathematik und ihrer Grenzgebiete. 3. Folge. A Series of
  Modern Surveys in Mathematics [Results in Mathematics and Related Areas. 3rd
  Series. A Series of Modern Surveys in Mathematics]}.
\newblock Springer-Verlag, Berlin, 2003.
\newblock Cavity and mean field models.

\bibitem{TalagrandSphericalSK}
Michel Talagrand.
\newblock Free energy of the spherical mean field model.
\newblock {\em Probab. Theory Related Fields}, 134(3):339--382, 2006.

\bibitem{TalagrandParisiFormula2006}
Michel Talagrand.
\newblock The {P}arisi formula.
\newblock {\em Ann. of Math. (2)}, 163(1):221--263, 2006.

\bibitem{TalagrandParisiMeasures2006a}
Michel Talagrand.
\newblock Parisi measures.
\newblock {\em J. Funct. Anal.}, 231(2):269--286, 2006.

\bibitem{Talagrand2007a}
Michel Talagrand.
\newblock {Large Deviations, Guerra's and A.S.S. Schemes, and the Parisi
  Hypothesis}.
\newblock {\em Journal of Statistical Physics}, 126(4):837--894, 2007.

\bibitem{Talagrand2007}
Michel Talagrand.
\newblock Mean field models for spin glasses: some obnoxious problems.
\newblock In {\em Spin glasses}, volume 1900 of {\em Lecture Notes in Math.},
  pages 63--80. Springer, Berlin, 2007.

\bibitem{Toubol1998}
Alain Toubol.
\newblock High temperature regime for a multidimensional
  {S}herrington-{K}irkpatrick model of spin glass.
\newblock {\em Probab. Theory Related Fields}, 110(4):497--534, 1998.

\end{thebibliography}

\end{document}